\documentclass[a4paper,12pt]{article}

\usepackage{amsmath}
\usepackage{amssymb}
\usepackage{amsthm}

\newtheorem{theorem}{Theorem}[section]
\newtheorem{lemma}[theorem]{Lemma}
\newtheorem{coro}[theorem]{Corollary}
\newtheorem{proposition}[theorem]{Proposition}
\theoremstyle{definition}
\newtheorem{example}[theorem]{Example}
\newtheorem{remark}[theorem]{Remark}
\newtheorem{rules}[theorem]{Rule}

\def\sln{\mathfrak{sl}_n}
\def\slnh{\widehat{\mathfrak{sl}}_n}
\def\geh{\mathfrak{g}}
\def\La{\Lambda}
\def\la{\lambda}
\def\ol#1{\overline{#1}}

\def\wt{\mbox{\sl wt}\,}
\def\cd{\cdots}
\def\ot{\otimes}

\def\veps{\varepsilon}
\def\vphi{\varphi}

\def\P{{\mathbb P}}
\def\Q{{\mathbb Q}}
\def\T{{\mathbb T}}
\def\bbt{\hat{\T}}
\def\Z{{\mathbb Z}}
\def\Zn{\Z_{\ge0}}
\def\etd{\tilde{e}}
\def\ftd{\tilde{f}}
\def\et#1{\tilde{e}_{#1}}
\def\ft#1{\tilde{f}_{#1}}

\newcommand{\fbx}[1]{\fbox{$\vphantom{\overline{1}} #1$}}
\newcommand{\bbx}[2]{\overbrace{\fbox{$\vphantom{\ol{1}} #1 \cd #1$}}^{#2}}
\newcommand{\mapright}[1]{%
  \smash{\mathop{%
   \hbox to 1cm{\rightarrowfill}}\limits^{#1}}}
\newcommand{\mapdown}[1]{\Big\downarrow
  \llap{$\vcenter{\hbox{$\scriptstyle#1\,$}}$ }}
\makeatletter
\@addtoreset{equation}{section}
\makeatother
\renewcommand{\theequation}{\arabic{section}.\arabic{equation}}
\title{Combinatorial $R$ matrices for a family \\
of crystals : $C^{(1)}_n$ and $A^{(2)}_{2n-1}$ cases}
\author{
Goro Hatayama\thanks{
Institute of Physics, University of Tokyo, Komaba, Tokyo 153-8902, Japan},
Atsuo Kuniba,$\hspace{-1.2mm}^*$\\
Masato Okado\thanks{Department of Informatics and Mathematical Science,
Graduate School of Engineering Science,
Osaka University,
Toyonaka, Osaka 560-8531,
Japan}
and Taichiro Takagi\thanks{
Department of Mathematics and Physics, National Defense Academy,
Yokosuka 239-8686, Japan}}
\date{}
\begin{document}
\maketitle
\begin{abstract}
The combinatorial $R$ matrices are obtained for
a family $\{B_l \}$
of crystals for $U'_q(C^{(1)}_n)$ and $U'_q(A^{(2)}_{2n-1})$,
where $B_l$ is the crystal of the irreducible module
corresponding to  the one-row Young diagram of length $l$.
The isomorphism $B_l\ot B_k \simeq B_k \ot B_l$ and the
energy function are described explicitly in terms of a $C_n$-analogue of the
Robinson-Schensted-Knuth type insertion algorithm.
As an application a $C^{(1)}_n$-analogue of
the Kostka  polynomials
is calculated for several cases.
\end{abstract}
\section{Introduction}
\label{sec:intro}
\subsection{Background}
{\em Physical combinatorics} might be defined na\"{\i}vely as
combinatorics guided by ideas or insights from physics.
A distinguished example can be given by the Kostka polynomial.
It is a polynomial $K_{\la\mu}(q)$ in $q$ depending on two
partitions $\la,\mu$ with the same number of nodes. Although
there are several aspects in this polynomial, one can regard
it as a $q$-analogue of the multiplicity of the irreducible
$\sln$-module $V_\la$ in the $m$-fold tensor product
$V_{(\mu_1)}\ot V_{(\mu_2)}\ot\cd\ot V_{(\mu_m)}$
($\mu=(\mu_1,\cd,\mu_m)$). Here for $\la=(\la_1,\cd,\la_n)$
($\la_1\ge\cd\ge\la_n\ge0$) $V_\la$ denotes the irreducible
$\sln$-module with highest weight
$\sum_{i=1}^{n-1}(\la_i-\la_{i+1})\La_i$, $\La_i$ being the
fundamental weight of $\sln$. In particular, $V_{(\mu_i)}$ is the
symmetric tensor representation of degree $\mu_i$.

In \cite{KR} Kirillov and Reshetikhin presented the following expression
for the Kostka polynomial\footnote{This expression differs from the
conventional definition of $K_{\la\mu}(q)$ by an overall power of $q$.}:
\begin{eqnarray}
K_{\la\mu}(q)&=&\sum_{\{m\}}q^{c(\{m\})}
\prod_{{\scriptstyle 1 \le a \le n-1} \atop
   {\scriptstyle i \ge 1}}
\left[ \begin{array}{c} p^{(a)}_i +  m^{(a)}_i
 \\   m^{(a)}_i \end{array} \right], \label{eq:KostkaFF}\\
c(\{ m  \}) & = & \frac{1}{2}\sum_{1 \le a, b \le n-1} C_{a b}
\sum_{i, j \ge 1} \mbox{min}(i, j) m^{(a)}_i m^{(b)}_j \nonumber\\
&& \qquad\quad - \sum_{i, j \ge 1} \mbox{min}(i, \mu_j)
m^{(1)}_i, \nonumber\\
p^{(a)}_i & = & \delta_{a 1} \sum_{j \ge 1} \mbox{min}(i, \mu_j)
- \sum _{1 \le b \le n-1} C_{a b} \sum_{j \ge 1} \mbox{min}(i,j) m^{(b)}_j,
\nonumber
\end{eqnarray}
where the sum $\sum_{\{m\}}$ is taken over $\{m^{(a)}_i\in\Zn\mid
1\le a\le n-1,i\ge1\}$ satisfying $p^{(a)}_i\ge0$ for $1\le a \le n-1, i\ge1$,
and $\sum_{i\ge1}im^{(a)}_i=\la_{a+1}+\la_{a+2}+\cd+\la_n$ for $1\le a\le n-1$.
$(C_{ab})_{1\le a,b\le n-1}$ is the Cartan matrix of $\sln$, and
$\left[M\atop N\right]$ is the $q$-binomial coefficient or Gaussian polynomial.
An intriguing point is that this expression was obtained through the string
hypothesis of the Bethe ansatz \cite{Be} for the $\sln$-invariant Heisenberg
chain, which is certainly in the field of physics.

Another important idea comes from Baxter's corner transfer matrix (CTM)
\cite{Ba,ABF}. In the course of the study of CTM eigenvalues, the notion of
one dimensional sum (1dsum) has appeared \cite{DJKMO},
and it was recognized that 1dsums give affine Lie algebra characters.
Such phenomena were clarified by the theory of perfect crystals
\cite{KMN1,KMN2}. As far as the Kostka polynomial is concerned, Nakayashiki
and Yamada obtained the following expression \cite{NY}:
\begin{equation} \label{eq:Kostka1dsum}
K_{\la\mu}(q)=\sum_p q^{E(p)},
\end{equation}
where $p$ ranges over the elements $p=b_1\ot\cd\ot b_m$ of
$B_{(\mu_1)}\ot\cd\ot B_{(\mu_m)}$ satisfying $\et{i}p=0$
($i=1,\cd,n-1$) and $\wt p=\sum_{i=1}^{n-1}(\la_i-\la_{i+1})\La_i$.
$B_{(\mu_i)}$ is the crystal base of the irreducible $U_q(\sln)$-module
with highest weight corresponding to $(\mu_i)$ and $\et{i}$ is the
so-called Kashiwara operator. $E(p)$ is called the energy of $p$
and calculated by using the energy function $H$ as
\[
E(p)=\sum_{1\le i<j\le m}H(b_i\ot b^{(i+1)}_j),
\]
where $b^{(i)}_j$ is defined through the crystal isomorphism:
\begin{eqnarray*}
B_{(\mu_i)}\ot B_{(\mu_{i+1})}\ot\cd\ot B_{(\mu_j)}
&\simeq&
B_{(\mu_j)}\ot B_{(\mu_i)}\ot\cd\ot B_{(\mu_{j-1})}\\
b_i\ot b_{i+1}\ot\cd\ot b_j
&\mapsto&
b^{(i)}_j\ot b'_i\ot\cd\ot b'_{j-1}.
\end{eqnarray*}
In the two-fold tensor case, the crystal isomorphism
$B_{(\mu_i)}\ot B_{(\mu_j)}\simeq B_{(\mu_j)}\ot B_{(\mu_i)}:
b_i\ot b_j\mapsto b'_j\ot b'_i$ combined with the value $H(b_i\ot b_j)$
is called the combinatorial $R$ matrix. The crystal base $B_{(l)}$ has
a generalization to the rectangular shape $B_{(l^k)}$, and the
corresponding generalization of the Kostka polynomial is considered in
\cite{SW,S}.

In view of the equality (\ref{eq:KostkaFF})=(\ref{eq:Kostka1dsum}), one
is led to an application of the perfect crystal theory of $B_{(l)}$.
Define a branching function $b^V_\la(q)$ for an $\slnh$-module $V$ by
\begin{eqnarray*}
b^V_\la(q)&=&\mbox{tr}\,_{{\cal H}(V,\la)}\;q^{-d},\\
{\cal H}(V,\la)&=&\{v\in V\mid e_iv=0\,(i=1,\cd,n-1),\wt v=\la\}.
\end{eqnarray*}
Here $d$ is the degree operator. Let $V(l\La_0)$ be the integrable
$\slnh$-module with affine highest weight $l\La_0$. Then
(\ref{eq:KostkaFF})=(\ref{eq:Kostka1dsum}) implies the spinon character
formula:
\begin{equation} \label{eq:spinon}
b^{V(l\La_0)}_\la(q)=\sum_\eta
\frac{K_{\xi\eta}(q)F^{(l)'}_\eta(q)}{(q)_{\zeta_1}\cd(q)_{\zeta_{n-1}}}.
\end{equation}
For the definitions of $\xi,(\zeta_1,\cd,\zeta_{n-1}),F^{(l)'}_\eta(q)$
along with the summing range of $\eta$, see Proposition 4.12 of \cite{HKKOTY}.

A key to the derivation of (\ref{eq:spinon}) is the fact that a suitable
subset of the semi-infinite tensor product $\cd\ot B_{(l)}\ot\cd\ot B_{(l)}$
can be identified with the crystal base $B(l\La_0)$ of the integrable
$U_q(\slnh)$-module with highest weight $l\La_0$. Since all components are
the same, such a case is called homogeneous. Recently, a generalization of
such results to inhomogeneous cases is obtained \cite{HKKOT}. For example,
a suitable subset of
\[
\cd\ot B_{(l_1+l_2)}\ot B_{(l_2)}\ot\cd\ot B_{(l_1+l_2)}\ot B_{(l_2)}
\ot B_{(l_1+l_2)}\ot B_{(l_2)}
\]
can be identified with $B(l_1\La_0)\ot B(l_2\La_0)$. Taking the corresponding
limit of $\mu$ in the equality (\ref{eq:KostkaFF})=(\ref{eq:Kostka1dsum}),
one obtains an expression for the branching function
$b^{V(l_1\La_0)\ot V(l_2\La_0)}_\la(q)$.

Another important application of the inhomogeneous case is found
in soliton cellular automata.
Recently several such automata have been related to
known soliton equations
through a limiting procedure called ultra-discretization \cite{TS,TTMS}.
Although they seem to have nothing to do with the theory of crystals
at first view, recent studies revealed their underlying crystal structure
\cite{HKT,FOY,HHIKTT}.
Namely, the combinatorial $R$ matrix appears as
the scattering rule of solitons as well as the time evolution operator
for the automaton.

\subsection{Present work}
In the $\slnh$ case, a typical example of
the isomorphism $B_{(3)} \ot B_{(2)} \simeq B_{(2)}\ot B_{(3)}$  is
\begin{displaymath}
112 \otimes 23 \mapsto 12 \otimes 123.
\end{displaymath}
It may be viewed as a scattering process of
two composite particles 112 and 23.
Through the collision the constituent particles are re-shuffled and
then recombined into two other composite particles 12 and 123.

In this paper we study the combinatorial $R$ matrices  for a family of
$U_q'(C^{(1)}_n)$ and $U_q'(A^{(2)}_{2n-1})$ crystals.
This includes a new type of examples as
\begin{align*}
123 \otimes \bar{2}\bar{1} &\mapsto 23 \otimes 0\bar{2}\bar{0}\quad
\text{for $U_q'(C^{(1)}_n)$ case},\\
&\mapsto 13 \otimes 1\bar{1}\bar{1}\quad
\text{for $U_q'(A^{(2)}_{2n-1})$ case}.
\end{align*}
Here we observe ``anti-particles", which undergo
a pair annihilation and a pair creation:
$(1) + (\bar{1}) \longrightarrow (0) + (\bar{0})$ or
$(2) + (\bar{2}) \longrightarrow (1) + (\bar{1})$.

We shall consider a family $\{B_l\mid l \in \Z_{\ge 1} \}$
of crystals for $U_q'(C^{(1)}_n)$ and $U_q'(A^{(2)}_{2n-1})$.
The above example corresponds to $B_3\ot B_2 \simeq B_2\ot B_3$.
Here  $B_l$ is the crystal of the irreducible $U_q'$-module
corresponding to the $l$-fold symmetric ``fusion" of the vector representation.
For $U_q'(A^{(2)}_{2n-1})$ it was constructed in \cite{KKM}.
For $U_q'(C^{(1)}_n)$,  $B_l$ in this paper
denotes $B_{l/2}$ in \cite{KKM}
($B_l$ in \cite{HKKOT}) when $l$ is even (odd).
Our main result is the explicit description of the
isomorphism  $B_l \ot B_k \simeq B_k\ot B_l$  and
the associated energy function for any $l$ and $k$.
It will be done through a slight modification of the
insertion algorithm for the $C$-tableaux
introduced by T.~H.~Baker \cite{B}.
Since the two affine algebras $C^{(1)}_n$ and $A^{(2)}_{2n-1}$
share the common classical part $C_n$,
they allow a parallel treatment and the results are similar in many respects.
Let us sketch them along the content of the paper.

In Section 2, we recall some basic facts about crystals.
As a $U_q(C_n)$ crystal, it is known that $U_q'(C^{(1)}_n)$ crystal
$B_l$ decomposes into the disjoint union of $B(l\Lambda_1), B((l-2)\Lambda_1), \ldots$,
where $B(\lambda)$ denotes the crystal of the irreducible representation
with highest weight $\lambda$.
Within each $B(l'\Lambda_1)$ it is natural \cite{KN,B} to parametrize
the elements by length $l'$ one-row semistandard tableaux
with  letters $1 < \cdots n < \ol{n} < \cdots \ol{2} < \ol{1}$.
Instead of doing so we will represent elements in $B_l$ uniformly via
length $l$  one-row  semistandard tableaux with letters
$0 < 1 < \cdots n < \ol{n} < \cdots \ol{2} < \ol{1} < \ol{0}$.
Here the number $x_0$ of $0$ and $\ol{x}_0$ of $\ol{0}$ must be the same,
according to which the elements belong to $B((l-2x_0)\Lambda_1)$.
Thus the number of letters in the tableaux has increased from $2n$ to  $2(n+1)$.
In fact, under the insertion scheme in later sections, these tableaux will
behave like those for $U_q(C_{n+1})$ \cite{B} in some sense.

In Section 3 we first define an insertion algorithm for the tableaux
introduced in Section 2.
When there is no $(x,\ol{x})$ pair, it is the same as the well known
$\sln$ case \cite{F}.
In general, our algorithm is essentially Baker's one \cite{B}
for $U_q(C_{n+1})$ if
$0 < \cdots < n < \ol{n} < \cdots < \ol{0}$
is regarded as $1 < \cdots < n+1 < \ol{n+1} < \cdots < \ol{1}$.
See Remark \ref{re:bakertonokankei}.
We describe it only for those tableaux with depth at most two, which
suffices for our aim.
We then state a main theorem, which describes the combinatorial $R$ matrix
of $U'_q(C^{(1)}_n)$  explicitly in terms of the insertion scheme.

In Section 4 we prove the main theorem.
As a $U_q(C_n)$ crystal, $B_l\ot B_k$ decomposes into connected components
which are isomorphic to the crystals of irreducible $U_q(C_n)$-modules.
Within each component the general elements are obtained by
applying $\tilde{f}_i$'s $(1 \le i \le n)$ to the $U_q(C_n)$ highest elements.
Our strategy is first to verify the theorem directly for the
highest elements.
For general elements the theorem follows from the
fact due to Baker that our insertion algorithm on letters $0,1, \ldots, \ol{1},\ol{0}$
can be regarded as the isomorphism of $U_q(C_{n+1})$ crystals.
It turns out that
$B_l \supset B_{l-2} \supset B_{l-4} \supset \cdots$
as the sets according to the number of $(0,\ol{0})$
pairs contained in the tableaux.
We shall utilize this fact to remove the $(0,\ol{0})$ pairs before the
insertion so as to avoid the pair annihilation of the boxes under the insertions
and the resulting bumping-sliding transition in \cite{B}.

In Section 5, a parallel treatment is done for $U'_q(A^{(2)}_{2n-1})$.
This case is simpler in that $B_l$ coincides with $B(l\Lambda_1)$
as a $U_q(C_n)$ crystal.
Consequently we do not have letters $0$ and $\ol{0}$ in the tableaux.
The main difference from $U'_q(C^{(1)})$ case is to remove $1$ and
$\ol{1}$ appropriately before the insertion.

In Appendix A, we detail the calculation for the proof of
Proposition \ref{th:main2}.

In Appendix B, another  rule for finding the image under
$B_l \ot B_k \simeq B_k \ot B_l$ is given for $U'_q(C_n)$ case.
In practical calculations it is often more efficient than the one based on
the insertion scheme in the main text.

In Appendix C, the $C^{(1)}_n$-analogue $X_{\lambda, \mu}(t)$
of the Kostka polynomials
in the sense of Section 1.1
is listed up to $\vert \mu \vert = 6$.
They coincide with the Kostka polynomial if $\vert \lambda \vert =
\vert \mu \vert$.

We remark that
the isomorphism $B_l\ot B_k \simeq B_k \ot B_l$
for $U'_q(C^{(1)}_n)$ in this paper has been identified
with the two body scattering rule in the soliton cellular automaton \cite{HKT}.

\vspace{0.4cm}
\noindent
{\bf Acknowledgements} \hspace{0.1cm}
The authors thank T. H. Baker for useful discussions.
They also thank
M. Kashiwara and  T. Miwa for
organizing the conference
``Physical Combinatorics''
in Kyoto during January 28 - February 2, 1999.

\section{Definitions}
\label{sec:defini}

%
\subsection{Brief summary of crystals}
Let $I$ be an index set.
A crystal $B$ is a set $B$ with the maps
\[
\et{i},\ft{i}: B\sqcup\{0\}\longrightarrow B\sqcup\{0\}
\quad (i \in I)
\]
satisfying the following properties:
\begin{itemize}
\item[] $\et{i}0=\ft{i}0=0$,
\item[] for any $b$ and $i$, there exists $n>0$ such that
$\et{i}^nb=\ft{i}^nb=0$,
\item[] for $b,b'\in B$ and $i\in I$, $\ft{i}b=b'$ if and only if
$b=\et{i}b'$.
\end{itemize}
For an element $b$ of $B$ we set
\[
\veps_i(b)=\max\{n\in\Zn\mid\et{i}^nb\neq0\},\quad
\vphi_i(b)=\max\{n\in\Zn\mid\ft{i}^nb\neq0\}.
\]
For two crystals $B$ and $B'$, the tensor product $B\ot B'$
is defined.
\[
B\ot B'=\{b\ot b'\mid b\in B,b'\in B'\}.
\]
The actions of $\et{i}$ and $\ft{i}$ are defined by
\begin{eqnarray}
\et{i}(b\ot b')&=&\left\{
\begin{array}{ll}
\et{i}b\ot b'&\mbox{ if }\vphi_i(b)\ge\veps_i(b')\\
b\ot \et{i}b'&\mbox{ if }\vphi_i(b) < \veps_i(b'),
\end{array}\right. \label{eq:ot-e}\\
\ft{i}(b\ot b')&=&\left\{
\begin{array}{ll}
\ft{i}b\ot b'&\mbox{ if }\vphi_i(b) > \veps_i(b')\\
b\ot \ft{i}b'&\mbox{ if }\vphi_i(b)\le\veps_i(b').
\end{array}\right. \label{eq:ot-f}
\end{eqnarray}
Here $0\ot b$ and $b\ot0$ are understood to be $0$.

%
\subsection{\mathversion{bold}Energy function and combinatorial $R$ matrix}

Let $\geh$ be an affine Lie algebra and
let $B$ and $B'$ be two $U_q'(\geh)$ crystals.
We assume that $B$ and $B'$ are finite sets,
and that $B \ot B'$ is connected.
The algebra $U_q'(\geh)$ is a subalgebra of
$U_q(\geh)$.
Their definitions are given in Section 2.1 (resp. 3.2) of \cite{KMN1}
for $U_q(\geh)$ (resp. $U_q'(\geh)$).

Suppose $b\ot b'\in B\ot B'$ is mapped to $\tilde{b'}\ot \tilde{b}
\in B'\ot B$ under the isomorphism
$B \ot B' \simeq B' \ot B$
of $U_q'(\geh)$ crystals.
A $\Z$-valued function
$H$ on $B\ot B'$ is called an {\em energy function} if for any $i$
and $b\ot b'\in B\ot B'$ such that $\et{i}(b\ot b')\neq0$,
 it satisfies
\begin{equation}\label{eq:e-func}
H(\et{i}(b\ot b'))=
\begin{cases}
	H(b\ot b')+1
	&\text{ if }i=0,\vphi_0(b)\geq\veps_0(b'),
	\vphi_0(\tilde{b'})\geq\veps_0(\tilde{b}),\\
	H(b\ot b')-1
	&\text{ if }i=0,\vphi_0(b)<\veps_0(b'),
	\vphi_0(\tilde{b'})<\veps_0(\tilde{b}),\\
	H(b\ot b')
	&\text{ otherwise}.
\end{cases}
\end{equation}
When we want to emphasize $B\ot B'$, we write $H_{BB'}$ for $H$.
This definition of the energy function
is due to (3.4.e) of \cite{NY}, which is a
generalization of the definition for $B=B'$ case
in \cite{KMN1}.
The energy function is unique up to an additive constant, since
$B\ot B'$ is connected. 
By definition, $H_{BB'}(b\ot b') - H_{B'B}(\tilde{b'}\ot \tilde{b})$ 
is a constant independent of 
$b\ot b'$. 
In this paper we choose the constant to be $0$.
We call the isomorphism
$B \ot B' \simeq B' \ot B$
endowed with the energy function $H_{BB'}$ the {\em combinatorial $R$-matrix}.

\subsection{\mathversion{bold}$C_n^{(1)}$ crystals}

Given a non-negative integer $l$,
we consider a $U_q'(C_n^{(1)})$ crystal denoted by $B_l$.
If $l$ is even, $B_{l}$ is the same as that defined
in \cite{KKM}.
(Their $B_l$ is identical to our $B_{2l}$.)
If $l$ is odd,  $B_{l}$ is defined in \cite{HKKOT}.
$B_l$'s are the crystals associated with
the crystal bases
of the irreducible finite dimensional representation
of the quantum affine algebra $U_q'(C_n^{(1)})$.
As a set $B_{l}$ reads
$$
B_{l} = \left\{(
x_1,\ldots, x_n,\overline{x}_n,\ldots,\overline{x}_1) \Biggm|
x_i, \overline{x}_i \in \Z_{\ge 0},
\sum_{i=1}^n(x_i + \overline{x}_i) \in \{l,l-2,\ldots \} \right\}.
$$
The crystal structure is given by (\ref{eq:structure}).

$B_{l}$ is isomorphic to
$\bigoplus_{0 \leq j \leq l,\, j \equiv l \pmod{2}} B(j \La_1)$ as
crystals for $U_q(C_n)$, where $B(j \La_1)$ is
the one associated with the irreducible representation of
with highest weight $j \Lambda_1$.
As a special case of the more general family of $U_q(C_n)$ crystals
\cite{KN},
the crystal $B(j \La_1)$ has a description with the
semistandard $C$-tableaux.
The entries are $1,\ldots ,n$ and $\ol{1}, \ldots ,\ol{n}$, with the
total order:
\begin{displaymath}
1 < 2 < \cd < n < \ol{n} < \cd < \ol{2} < \ol{1}.
\end{displaymath}
In this description
$b= (x_1, \ldots, x_n, \overline{x}_n,\ldots,\overline{x}_1) \in
B(j \La_1)$
is depicted by
\begin{equation}
b=\overbrace{\fbox{$\vphantom{\ol{1}} 1 \cd 1$}}^{x_1}\!
\fbox{$\vphantom{\ol{1}}\cd$}\!
\overbrace{\fbox{$\vphantom{\ol{1}}n \cd n$}}^{x_n}\!
\overbrace{\fbox{$\vphantom{\ol{1}}\ol{n} \cd \ol{n}$}}^{\ol{x}_n}\!
\fbox{$\vphantom{\ol{1}}\cd$}\!
\overbrace{\fbox{$\ol{1} \cd \ol{1}$}}^{\ol{x}_1}.
\end{equation}
The length of this one-row tableau is equal to $j$, namely
$\sum_{i=1}^n(x_i + \overline{x}_i) =j$.
Here and in the remaining part of this paper we denote
$\overbrace{\fbx{i} \fbx{i}  \fbx{\cd}  \fbx{i}}^{x}$ by
\par\noindent
\setlength{\unitlength}{5mm}
\begin{picture}(22,3)(-6,0)
\put(0,0){\makebox(10,3)
{$\overbrace{\fbox{$\vphantom{\ol{1}} i \cd i$}}^{x}$ or
more simply by }}
\put(10,0.5){\line(1,0){3}}
\put(10,1.5){\line(1,0){3}}
\put(10,0.5){\line(0,1){1}}
\put(13,0.5){\line(0,1){1}}
\put(10,0.5){\makebox(3,1){$i$}}
\put(10,1.5){\makebox(3,1){${\scriptstyle x}$}}
\put(13,0){\makebox(1,1){.}}
\end{picture}
\par\noindent

We shall depict the elements of
$B_{l}$ 
by one-row tableaux with length $l$, by supplying
pairs of \fbx{0} and \fbx{\ol{0}}.
Adding $0$ and $\ol{0}$ into the set of the entries of the
tableaux, we assume
the total order $0 < 1 < \cd <  \ol{1} < \ol{0}$.
Thus we depict
$b= (x_1, \ldots, x_n, \overline{x}_n,\ldots,\overline{x}_1) \in
B_{l}$ by
\begin{equation}\label{eq:boxx}
\T (b)=\bbx{0}{x_0}\!
\overbrace{\fbox{$\vphantom{\ol{1}} 1 \cd 1$}}^{x_1}\!
\fbox{$\vphantom{\ol{1}}\cd$}\!
\overbrace{\fbox{$\vphantom{\ol{1}}n \cd n$}}^{x_n}\!
\overbrace{\fbox{$\vphantom{\ol{1}}\ol{n} \cd \ol{n}$}}^{\ol{x}_n}\!
\fbox{$\vphantom{\ol{1}}\cd$}\!
\overbrace{\fbox{$\ol{1} \cd \ol{1}$}}^{\ol{x}_1}\!
\overbrace{\fbox{$\ol{0} \cd \ol{0}$}}^{\ol{x}_0},
\end{equation}
where $x_0 = \overline{x}_0 = (l-\sum_{i=1}^n (x_i + \overline{x}_i))/2$.
If $x_0 = 0$ we say that $\T (b)$ has no \fbx{0}.
Sometimes we identify  $\T (b)$ with $b$, and
omit the frame of $\T (b)$, e.g., 
$0\bar{2}\bar{0} = (0,\ldots, 0,1,0)  \in B_3$.

This description means that we have embedded $B_{l}$, as a set, into
$U_q(C_{n+1})$ crystal $B(l \Lambda_1)$.
Let us denote by $\varsigma$ this embedding:
\begin{equation}
\varsigma :
\mbox{$U_q(C_n)$ crystal $B_{l}$ as a set} \hookrightarrow
\mbox{$U_q(C_{n+1})$ crystal $B(l \Lambda_1)$}.
\label{eq:embedding}
\end{equation}
It shifts the entries of the tableaux as
$\varsigma (i) = i+1$ and $\varsigma (\ol{i}) = \ol{i+1}$ for
$i=0,1,\ldots,n$.
For example $\varsigma(0\bar{2}\bar{0}) = 1\bar{3}\bar{1}$.
For $b \in B_l$ and $i=1,\ldots,n$ one has $\varsigma (\et{i}b) = \et{i+1} \varsigma (b)$
and $\varsigma (\ft{i}b) = \ft{i+1} \varsigma (b)$.
The crystal structure of $B_l$ is give by
{\allowdisplaybreaks
\begin{eqnarray}
\et{0} b &=&
\begin{cases}
	(x_1-2,x_2,\ldots,\overline{x}_2,\overline{x}_1)  &
	\mbox{ if }x_1 \ge \overline{x}_1+2, \\
	(x_1-1,x_2,\ldots,\overline{x}_2,\overline{x}_1+1) &
	\mbox{ if } x_1 = \overline{x}_1+1,\\
	(x_1,x_2,\ldots,\overline{x}_2,\overline{x}_1+2) &
	\mbox{ if } x_1 \le \overline{x}_1,\\
\end{cases}\nonumber\\
\et{n} b &=& (x_1,\ldots,x_n+1,\overline{x}_n-1,\ldots,\overline{x}_1),\nonumber\\
\et{i} b &=&
\begin{cases}
	(x_1,\ldots,x_i+1,x_{i+1}-1,\ldots,\overline{x}_1) &
	\mbox{ if } x_{i+1} > \overline{x}_{i+1},\\
	(x_1,\ldots,\overline{x}_{i+1}+1,\overline{x}_i-1,\ldots,\overline{x}_1) &
	\mbox{ if } x_{i+1} \le \overline{x}_{i+1},
\end{cases}\nonumber\\[4pt]
\ft{0} b &=&
\begin{cases}
	 (x_1+2,x_2,\ldots,\overline{x}_2,\overline{x}_1) &
	\mbox{ if }x_1 \ge \overline{x}_1,\\
	(x_1+1,x_2,\ldots,\overline{x}_2,\overline{x}_1-1) &
	\mbox{ if } x_1 = \overline{x}_1-1,\\
	(x_1,x_2,\ldots,\overline{x}_2,\overline{x}_1-2) &
	\mbox{ if } x_1 \le \overline{x}_1-2,
\end{cases}\nonumber\\
\ft{n} b &=&
	(x_1,\ldots,x_n-1,\overline{x}_n+1,\ldots,\overline{x}_1),\nonumber\\
\ft{i} b &=&
\begin{cases}
	(x_1,\ldots,x_i-1,x_{i+1}+1,\ldots,\overline{x}_1) &
	\mbox{ if } x_{i+1} \ge \overline{x}_{i+1},\\
	(x_1,\ldots,\overline{x}_{i+1}-1,\overline{x}_i+1,\ldots,\overline{x}_1) &
	\mbox{ if } x_{i+1} < \overline{x}_{i+1},
\end{cases}
\label{eq:structure}
\end{eqnarray}
}
where $b= (x_1, \ldots, x_n, \overline{x}_n,\ldots,\overline{x}_1)$
and $i = 1, \ldots, n-1$.
For this $b$ we have
\begin{eqnarray}
\vphi_i (b) &=& x_i + (\ol{x}_{i+1}-x_{i+1})_+
\quad \mbox{for} \, i=0,1,\ldots,n-1, \nonumber\\
\veps_i (b) &=& \ol{x}_i + (x_{i+1}-\ol{x}_{i+1})_+
\quad \mbox{for} \, i=0,1,\ldots,n-1, \nonumber\\
\vphi_n (b) &=& x_n , \quad \veps_n (b) = \ol{x}_n. \label{eq:eps}
\end{eqnarray}
Here $(x)_+ := \max (x,0)$.

Except for Section \ref{sec:discuss} 
concerning $A^{(2)}_{2n-1}$,  we normalize the energy function for $C^{(1)}_n$ case as
\begin{equation*}
H_{B_lB_k}((l,0,\ldots,0) \ot (0,k,0,\ldots,0)) =0,
\end{equation*}
irrespective of $l < k$ or $l \ge k$.
(For $A^{(2)}_{2n-1}$ we will employ a different normalization.
See (\ref{eq:kikaku}).)

\section{Explicit description of isomorphism and energy function}
\label{sec:descri}
\subsection{The algorithm of column insertions}

Set an alphabet $\mathcal{X}=\mathcal{A} \sqcup \bar{\mathcal{A}},\,
\mathcal{A}=\{ 0,1,\dots,n\}$ and
$\bar{\mathcal{A}}=\{\ol{0},\ol{1},\dots,\ol{n}\}$,
with the total order
$0 < 1 < \dots < n < \ol{n} < \dots < \ol{1} < \ol{0}$.
Unless otherwise stated, a tableau means a (column-strict) semistandard  one
with entries taken from $\mathcal{X}$
in Section \ref{sec:descri} and \ref{sec:theor}.
For the alphabet $\mathcal{X}$, we follow the convention
that Greek letters $ \alpha, \beta, \ldots $ belong to
$\mathcal{A} \sqcup \bar{\mathcal{A}}$ while Latin
letters $x,y,\ldots$ (resp. $\ol{x},\ol{y},\ldots$) belong to
$\mathcal{A}$ (resp. $\bar{\mathcal{A}}$).

Given a letter $\alpha \in \mathcal{X}$ and
the tableau $T$ that have at most two rows,
we define a tableau denoted by $(\fbx{\alpha} \rightarrow T)$, and
call such an algorithm a ``column insertion of a
letter $\alpha$ into a tableau $T$~".
(We sometimes identify a letter $\alpha$ with a box \fbx{\alpha}.)
Let us begin with such $T$'s that have at most one column.
The procedure of the column insertion $(\fbx{\alpha} \rightarrow T)$
can be summarized as follows:
\par\noindent
\setlength{\unitlength}{5.5mm}
\begin{picture}(22,1.5)(-3,0)
\put(0,0){\makebox(0,1)[r]{case 1a}}
\put(0,0){\makebox(1,1){$\Bigl($}}
\put(1,0){\line(1,0){1}}
\put(1,1){\line(1,0){1}}
\put(1,0){\line(0,1){1}}
\put(2,0){\line(0,1){1}}
\put(1,0){\makebox(1,1){$\alpha$}}
\put(2,0){\makebox(1,1){$\to$}}
\put(3,0){\makebox(1,1){$\emptyset$}}
\put(4,0){\makebox(1,1){$\Bigr)$}}
\put(4.8,0){\makebox(1,1){$=$}}
\put(6,0){\line(1,0){1}}
\put(6,1){\line(1,0){1}}
\put(6,0){\line(0,1){1}}
\put(7,0){\line(0,1){1}}
\put(6,0){\makebox(1,1){$\alpha$}}
\put(9,0){\makebox(10,1)[l]{,}}
\end{picture}
\par\noindent
\setlength{\unitlength}{5.5mm}
\begin{picture}(22,2.5)(-3,0)
\put(0,0){\makebox(0,2)[r]{case 2a}}
\put(0,0.5){\makebox(1,1){$\Bigl($}}
\put(1,0.5){\line(1,0){1}}
\put(1,1.5){\line(1,0){1}}
\put(1,0.5){\line(0,1){1}}
\put(2,0.5){\line(0,1){1}}
\put(1,0.5){\makebox(1,1){$\beta$}}
\put(2,0.5){\makebox(1,1){$\to$}}
\put(3,0.5){\line(1,0){1}}
\put(3,1.5){\line(1,0){1}}
\put(3,0.5){\line(0,1){1}}
\put(4,0.5){\line(0,1){1}}
\put(3,0.5){\makebox(1,1){$\alpha$}}
\put(4,0.5){\makebox(1,1){$\Bigr)$}}
\put(4.8,0.5){\makebox(1,1){$=$}}
\put(6,0){\line(1,0){1}}
\put(6,1){\line(1,0){1}}
\put(6,2){\line(1,0){1}}
\put(6,0){\line(0,1){2}}
\put(7,0){\line(0,1){2}}
\put(6,1){\makebox(1,1){$\alpha$}}
\put(6,0){\makebox(1,1){$\beta$}}
\put(9,0){\makebox(10,2)[l]{if $\alpha < \beta$,}}
\end{picture}
\par\noindent
\setlength{\unitlength}{5.5mm}
\begin{picture}(22,1.5)(-3,0)
\put(0,0){\makebox(0,1)[r]{case 1b}}
\put(0,0){\makebox(1,1){$\Bigl($}}
\put(1,0){\line(1,0){1}}
\put(1,1){\line(1,0){1}}
\put(1,0){\line(0,1){1}}
\put(2,0){\line(0,1){1}}
\put(1,0){\makebox(1,1){$\alpha$}}
\put(2,0){\makebox(1,1){$\to$}}
\put(3,0){\line(1,0){1}}
\put(3,1){\line(1,0){1}}
\put(3,0){\line(0,1){1}}
\put(4,0){\line(0,1){1}}
\put(3,0){\makebox(1,1){$\beta$}}
\put(4,0){\makebox(1,1){$\Bigr)$}}
\put(4.8,0){\makebox(1,1){$=$}}
\put(6,0){\line(1,0){2}}
\put(6,1){\line(1,0){2}}
\put(6,0){\line(0,1){1}}
\put(7,0){\line(0,1){1}}
\put(8,0){\line(0,1){1}}
\put(6,0){\makebox(1,1){$\alpha$}}
\put(7,0){\makebox(1,1){$\beta$}}
\put(9,0){\makebox(10,1)[l]{if $\alpha \le \beta$,}}
\end{picture}
\par\noindent
\setlength{\unitlength}{5.5mm}
\begin{picture}(22,2.5)(-3,0)
\put(0,0){\makebox(0,2)[r]{case 2b}}
\put(0,0){\makebox(1,2){$\Biggl($}}
\put(1,0){\line(1,0){1}}
\put(1,1){\line(1,0){1}}
\put(1,0){\line(0,1){1}}
\put(2,0){\line(0,1){1}}
\put(1,0){\makebox(1,1){$\beta$}}
\put(2,0){\makebox(1,1){$\to$}}
\put(3,0){\line(1,0){1}}
\put(3,1){\line(1,0){1}}
\put(3,2){\line(1,0){1}}
\put(3,0){\line(0,1){2}}
\put(4,0){\line(0,1){2}}
\put(3,0){\makebox(1,1){$\gamma$}}
\put(3,1){\makebox(1,1){$\alpha$}}
\put(4,0){\makebox(1,2){$\Biggr)$}}
\put(4.8,0){\makebox(1,2){$=$}}
\put(6,0){\line(1,0){1}}
\put(6,1){\line(1,0){2}}
\put(6,2){\line(1,0){2}}
\put(6,0){\line(0,1){2}}
\put(7,0){\line(0,1){2}}
\put(8,1){\line(0,1){1}}
\put(6,0){\makebox(1,1){$\beta$}}
\put(6,1){\makebox(1,1){$\alpha$}}
\put(7,1){\makebox(1,1){$\gamma$}}
\put(9,0){\makebox(10,2)[l]{if
$\alpha < \beta \leq \gamma$ and $(\alpha,\gamma) \ne (x,\ol{x})$,}}
\end{picture}
\par\noindent
\setlength{\unitlength}{5.5mm}
\begin{picture}(22,2.5)(-3,0)
\put(0,0){\makebox(0,2)[r]{case 3b}}
\put(0,0){\makebox(1,2){$\Biggl($}}
\put(1,0){\line(1,0){1}}
\put(1,1){\line(1,0){1}}
\put(1,0){\line(0,1){1}}
\put(2,0){\line(0,1){1}}
\put(1,0){\makebox(1,1){$\alpha$}}
\put(2,0){\makebox(1,1){$\to$}}
\put(3,0){\line(1,0){1}}
\put(3,1){\line(1,0){1}}
\put(3,2){\line(1,0){1}}
\put(3,0){\line(0,1){2}}
\put(4,0){\line(0,1){2}}
\put(3,0){\makebox(1,1){$\gamma$}}
\put(3,1){\makebox(1,1){$\beta$}}
\put(4,0){\makebox(1,2){$\Biggr)$}}
\put(4.8,0){\makebox(1,2){$=$}}
\put(6,0){\line(1,0){1}}
\put(6,1){\line(1,0){2}}
\put(6,2){\line(1,0){2}}
\put(6,0){\line(0,1){2}}
\put(7,0){\line(0,1){2}}
\put(8,1){\line(0,1){1}}
\put(6,0){\makebox(1,1){$\gamma$}}
\put(6,1){\makebox(1,1){$\alpha$}}
\put(7,1){\makebox(1,1){$\beta$}}
\put(9,0){\makebox(10,2)[l]{if
$\alpha \leq \beta < \gamma$ and
$(\alpha, \gamma) \ne (x,\ol{x})$,}}
\end{picture}
\par\noindent
\setlength{\unitlength}{5.5mm}
\begin{picture}(22,2.5)(-3,0)
\put(0,0){\makebox(0,2)[r]{case 4b}}
\put(0,0){\makebox(1,2){$\Biggl($}}
\put(1,0){\line(1,0){1}}
\put(1,1){\line(1,0){1}}
\put(1,0){\line(0,1){1}}
\put(2,0){\line(0,1){1}}
\put(1,0){\makebox(1,1){$\beta$}}
\put(2,0){\makebox(1,1){$\to$}}
\put(3,0){\line(1,0){1}}
\put(3,1){\line(1,0){1}}
\put(3,2){\line(1,0){1}}
\put(3,0){\line(0,1){2}}
\put(4,0){\line(0,1){2}}
\put(3,0){\makebox(1,1){$\ol{x}$}}
\put(3,1){\makebox(1,1){$x$}}
\put(4,0){\makebox(1,2){$\Biggr)$}}
\put(4.8,0){\makebox(1,2){$=$}}
\put(6,0){\line(1,0){2}}
\put(6,1){\line(1,0){4}}
\put(6,2){\line(1,0){4}}
\put(6,0){\line(0,1){2}}
\put(8,0){\line(0,1){2}}
\put(10,1){\line(0,1){1}}
\put(6,0){\makebox(2,1){$\beta$}}
\put(6,1){\makebox(2,1){$x\!-\!1$}}
\put(8,1){\makebox(2,1){$\ol{x\!-\!1}$}}
\put(11,0){\makebox(8,2)[l]{if
$x \le \beta \le \ol{x}$ and $ x \ne 0$,}}
\end{picture}
\par\noindent
\setlength{\unitlength}{5.5mm}
\begin{picture}(22,2.5)(-3,0)
\put(0,0){\makebox(0,2)[r]{case 5b}}
\put(0,0){\makebox(1,2){$\Biggl($}}
\put(1,0){\line(1,0){1}}
\put(1,1){\line(1,0){1}}
\put(1,0){\line(0,1){1}}
\put(2,0){\line(0,1){1}}
\put(1,0){\makebox(1,1){$x$}}
\put(2,0){\makebox(1,1){$\to$}}
\put(3,0){\line(1,0){1}}
\put(3,1){\line(1,0){1}}
\put(3,2){\line(1,0){1}}
\put(3,0){\line(0,1){2}}
\put(4,0){\line(0,1){2}}
\put(3,0){\makebox(1,1){$\ol{x}$}}
\put(3,1){\makebox(1,1){$\beta$}}
\put(4,0){\makebox(1,2){$\Biggr)$}}
\put(4.8,0){\makebox(1,2){$=$}}
\put(6,0){\line(1,0){2}}
\put(6,1){\line(1,0){3}}
\put(6,2){\line(1,0){3}}
\put(6,0){\line(0,1){2}}
\put(8,0){\line(0,1){2}}
\put(9,1){\line(0,1){1}}
\put(6,0){\makebox(2,1){$\ol{x\!+\!1}$}}
\put(6,1){\makebox(2,1){$x\!+\!1$}}
\put(8,1){\makebox(1,1){$\beta$}}
\put(11,0){\makebox(8,2)[l]{if
$x < \beta < \ol{x}$ and $x \ne n$.}}
\end{picture}
\par\noindent
The cases 2b - 5b do not cover all the tableaux with two rows,
but we only deal with these situations in this paper.
In particular, the tableaux generated by these insertions
have at most two rows.
Note that the algorithm except for the cases  4b and 5b
agrees with the Knuth-type column insertion.
We call  the cases 1b - 5b  the ``bumping cases".

When  $T$ is a general tableau with at most two rows, we repeat the above procedure:
we insert a box into the leftmost
column of $T$ according to the above formula.
If it is not a bumping case,
replace the column by the right hand side of the formula.
Otherwise,
replace the column by the right hand side of the formula without
the right box. We regard that this right box is bumped.
We insert it into the second column of $T$ from the left
and repeat the procedure above until we come to
a non-bumping case 1a or 2a.

\begin{example}
$n=4$.\par\noindent
\setlength{\unitlength}{4.8mm}
\begin{picture}(29,4.5)(0,-2)
\put(0,0){\makebox(1,2){$\Biggr($}}
\put(1,0){\line(1,0){1}}
\put(1,1){\line(1,0){1}}
\put(1,0){\line(0,1){1}}
\put(2,0){\line(0,1){1}}
\put(1,0){\makebox(1,1){$2$}}
\put(2,0){\makebox(1,1){$\to$}}
\put(3,0){\line(1,0){3}}
\put(3,1){\line(1,0){4}}
\put(3,2){\line(1,0){4}}
\multiput(3,0)(1,0){4}{\line(0,1){2}}
\put(7,1){\line(0,1){1}}
\put(3,0){\makebox(1,1){$4$}}
\put(4,0){\makebox(1,1){$\ol{3}$}}
\put(5,0){\makebox(1,1){$\ol{1}$}}
\put(3,1){\makebox(1,1){$0$}}
\put(4,1){\makebox(1,1){$3$}}
\put(5,1){\makebox(1,1){$4$}}
\put(6,1){\makebox(1,1){$\bar{0}$}}
\put(7,0){\makebox(1,2){$\Biggr)$}}
\put(7.8,0){\makebox(1,2){$=$}}
\multiput(0,0)(5,0){3}{
	\put(9,0){\line(1,0){3}}
	\put(9,1){\line(1,0){4}}
	\put(9,2){\line(1,0){4}}
	\multiput(9,0)(1,0){4}{\line(0,1){2}}
	\put(13,1){\line(0,1){1}}
}
\put( 9,0){\makebox(1,1){$2$}}
\put(10,0){\makebox(1,1){$\ol{3}$}}
\put(11,0){\makebox(1,1){$\ol{1}$}}
\put(9,1){\makebox(1,1){$0$}}
\put(10,1){\makebox(1,1){$3$}}
\put(11,1){\makebox(1,1){$4$}}
\put(12,1){\makebox(1,1){$\bar{0}$}}
\put(13,0){\makebox(1,2){$=$}}
\put(10,-1){\makebox(1,1){$\uparrow$}}
\multiput(10,-2)(1,0){2}{\line(0,1){1}}
\multiput(10,-2)(0,1){2}{\line(1,0){1}}
\put(10,-2){\makebox(1,1){$4$}}
\put(14,0){\makebox(1,1){$2$}}
\put(15,0){\makebox(1,1){$4$}}
\put(16,0){\makebox(1,1){$\ol{1}$}}
\put(14,1){\makebox(1,1){$0$}}
\put(15,1){\makebox(1,1){$2$}}
\put(16,1){\makebox(1,1){$4$}}
\put(17,1){\makebox(1,1){$\bar{0}$}}
\put(18,0){\makebox(1,2){$=$}}
\put(16,-1){\makebox(1,1){$\uparrow$}}
\multiput(16,-2)(1,0){2}{\line(0,1){1}}
\multiput(16,-2)(0,1){2}{\line(1,0){1}}
\put(16,-2){\makebox(1,1){$\ol{2}$}}
\put(19,0){\makebox(1,1){$2$}}
\put(20,0){\makebox(1,1){$4$}}
\put(21,0){\makebox(1,1){$\ol{2}$}}
\put(19,1){\makebox(1,1){$0$}}
\put(20,1){\makebox(1,1){$2$}}
\put(21,1){\makebox(1,1){$4$}}
\put(22,1){\makebox(1,1){$\bar{0}$}}
\put(23,0){\makebox(1,2){$=$}}
\put(22,-1){\makebox(1,1){$\uparrow$}}
\multiput(22,-2)(1,0){2}{\line(0,1){1}}
\multiput(22,-2)(0,1){2}{\line(1,0){1}}
\put(22,-2){\makebox(1,1){$\ol{1}$}}
\put(24,0){\line(1,0){3}}
\put(24,1){\line(1,0){5}}
\put(24,2){\line(1,0){5}}
\multiput(24,0)(1,0){4}{\line(0,1){2}}
\multiput(28,1)(1,0){2}{\line(0,1){1}}
\put(24,0){\makebox(1,1){$2$}}
\put(25,0){\makebox(1,1){$4$}}
\put(26,0){\makebox(1,1){$\ol{2}$}}
\put(24,1){\makebox(1,1){$0$}}
\put(25,1){\makebox(1,1){$2$}}
\put(26,1){\makebox(1,1){$4$}}
\put(27,1){\makebox(1,1){$\bar{1}$}}
\put(28,1){\makebox(1,1){$\bar{0}$}}
\end{picture}
\end{example}

For a tableau $T$ we denote by $w(T)$ the
Japanese reading word of $T$.
The $w(T)$ is a sequence of letters that is created by
reading all letters on $T$ from the rightmost column
to the leftmost column,
and in each column from the top to the bottom.
For instance,
\par\noindent
\setlength{\unitlength}{5mm}
\begin{picture}(22,1.5)(-5,0)
\put(0,0){\makebox(1,1){$w($}}
\put(1,0){\line(1,0){5}}
\put(1,1){\line(1,0){5}}
\put(1,0){\line(0,1){1}}
\put(2,0){\line(0,1){1}}
\put(3,0){\line(0,1){1}}
\put(5,0){\line(0,1){1}}
\put(6,0){\line(0,1){1}}
\put(1,0){\makebox(1,1){$\alpha_1$}}
\put(2,0){\makebox(1,1){$\alpha_2$}}
\put(3,0){\makebox(2,1){$\cd$}}
\put(5,0){\makebox(1,1){$\alpha_j$}}
\put(6,0){\makebox(6,1){$)=\alpha_j \cd \alpha_2 \alpha_1,$}}
\end{picture}
\par\noindent
and
\par\noindent
\setlength{\unitlength}{5mm}
\begin{picture}(22,2.5)(-2,0)
\put(0,0){\makebox(2,2){$w \Biggl($}}
\put(2,0){\line(1,0){5}}
\put(2,1){\line(1,0){10}}
\put(2,2){\line(1,0){10}}
\put(2,0){\line(0,1){2}}
\put(3,0){\line(0,1){2}}
\put(4,0){\line(0,1){2}}
\put(6,0){\line(0,1){2}}
\put(7,0){\line(0,1){2}}
\put(9,1){\line(0,1){1}}
\put(11,1){\line(0,1){1}}
\put(12,1){\line(0,1){1}}
\put(2,0){\makebox(1,1){$\beta_1$}}
\put(2,1){\makebox(1,1){$\alpha_1$}}
\put(3,0){\makebox(1,1){$\beta_2$}}
\put(3,1){\makebox(1,1){$\alpha_2$}}
\put(4,0){\makebox(2,1){$\cd$}}
\put(4,1){\makebox(2,1){$\cd$}}
\put(6,0){\makebox(1,1){$\beta_i$}}
\put(6,1){\makebox(1,1){$\alpha_i$}}
\put(7,1){\makebox(2,1){$\alpha_{i+1}$}}
\put(9,1){\makebox(2,1){$\cd$}}
\put(11,1){\makebox(1,1){$\alpha_j$}}
\put(12,0){\makebox(14,2){$\Biggr) =\alpha_j \cd \alpha_{i+1}
\alpha_i \beta_i \cd
\alpha_2 \beta_2 \alpha_1 \beta_1 ,$}}
\end{picture}
\par\noindent
and so on.
Let $T$ and $T'$ be one-row tableaux.
By abuse of  notation we denote by $T' \to T$
the tableau constructed by
successive
column insertions of the letters of the word $w(T')$ into $T$.
Namely if
\begin{displaymath}
w(T') = \tau_j \tau_{j-1} \cd \tau_1,
\end{displaymath}
then we write
\begin{displaymath}
(T' \to T) = (\tau_1 \to \cd (\tau_{j-1} \to ( \tau_j \to T ) ) \cd ).
\end{displaymath}
(Following a usual convention in type $A$ \cite{F},
it might be written as a {\em product tableau}, $T' \cdot T$.)
In particular $(T \to \emptyset) = T$ for any $T$.

Throughout this paper we let $T_1$ be the length of the
first row of a tableau $T$.

\begin{remark}\label{re:bakertonokankei}
Our algorithm is a specialization of
the column insertion for $C_n$-case \cite{B}.
Let $b_i$ be an element of $U_q(C_{n+1})$ crystal $B(l_i \Lambda_1)\,
(i=1,2)$.
Denote by $b_2 \stackrel{\ast}{\longrightarrow} b_1$
the tableau obtained by  successive
column insertions with the original definition \cite{B} of
$w(\T(b_2))$ into $b_1$.
(In \cite{B}, $b_2 \stackrel{\ast}{\longrightarrow} b_1$ is
denoted by $b_1 * b_2$.)
Then, for any element $b_i$ of $U'_q(C_{n}^{(1)})$ crystal $B_{l_i}$
such that  $\T (b_1)$ or $\T (b_2)$ has no \fbx{0},
our $\bigl( \T (b_2) \rightarrow \T (b_1) \bigr)$ has been determined so that
\[
\varsigma \Bigl( \bigl( \T (b_2) \rightarrow \T (b_1) \bigr) \Bigr)
=
\varsigma\bigl( \T(b_2) \bigr) \stackrel{\ast}{\longrightarrow}
	\varsigma\bigl( \T(b_1) \bigr).
\]
We will  calculate $\bigl( \T (b_2) \rightarrow \T (b_1) \bigr)$
only when  $\T(b_1)$ or $\T(b_2)$ has no \fbx{0} in this paper.
Under such a situation no pair annihilation takes place  during  the insertion
in the right hand side.
See also Remark \ref{rem:avoid}.
\end{remark}

\vskip3pt

We shall also use the reverse bumping algorithm \cite{B}.
In our case where the tableau has at most two rows,
the algorithm is rather simple.
We only use them in the following five cases.
\par\noindent
\setlength{\unitlength}{5.5mm}
\begin{picture}(22,3.5)(-3,0)
\put(-0.5,0){\makebox(0,1)[r]{case 1c}}
\put(0,0){
	\multiput(0,0)(0,1){2}{\line(1,0){1}}
	\multiput(0,0)(1,0){2}{\line(0,1){1}}
	\put(0,0){\makebox(1,1){$\alpha$}}
	}
\put(4,0){
	\multiput(0,0)(0,1){2}{\line(1,0){1}}
	\multiput(0,0)(1,0){2}{\line(0,1){1}}
	\put(0,0){\makebox(1,1){$\beta$}}
	}
\put(0,2){
	\multiput(0,0)(0,1){2}{\line(1,0){1}}
	\multiput(0,0)(1,0){2}{\line(0,1){1}}
	\put(0,0){\makebox(1,1){$\beta$}}
	}
\put(3,2){
	\multiput(0,0)(0,1){2}{\line(1,0){1}}
	\multiput(0,0)(1,0){2}{\line(0,1){1}}
	\put(0,0){\makebox(1,1){$\alpha$}}
	}
\put(1.5,0){\makebox(1,1){$=$}}
\put(0,1){\makebox(1,1){$\downarrow$}}
\put(3,1){\makebox(1,1){$\downarrow$}}
\put(6,0){\makebox(0,1)[l]{if $\alpha \le \beta$,}}
\end{picture}
\par\noindent
\setlength{\unitlength}{5.5mm}
\begin{picture}(22,5.5)(-3,0)
\put(-0.5,0){\makebox(0,2)[r]{case 2c}}
\put(0,0){
	\multiput(0,0)(0,1){3}{\line(1,0){1}}
	\multiput(0,0)(1,0){2}{\line(0,1){2}}
	\put(0,1){\makebox(1,1){$\alpha$}}
	\put(0,0){\makebox(1,1){$\beta$}}
	}
\put(4,0){
	\multiput(0,0)(0,1){3}{\line(1,0){1}}
	\multiput(0,0)(1,0){2}{\line(0,1){2}}
	\put(0,1){\makebox(1,1){$\alpha$}}
	\put(0,0){\makebox(1,1){$\gamma$}}
	}
\put(0,3){
	\multiput(0,0)(0,1){2}{\line(1,0){1}}
	\multiput(0,0)(1,0){2}{\line(0,1){1}}
	\put(0,0){\makebox(1,1){$\gamma$}}
	}
\put(3,3){
	\multiput(0,0)(0,1){2}{\line(1,0){1}}
	\multiput(0,0)(1,0){2}{\line(0,1){1}}
	\put(0,0){\makebox(1,1){$\beta$}}
	}
\put(1.5,0){\makebox(1,2){$=$}}
\put(0,2){\makebox(1,1){$\downarrow$}}
\put(3,2){\makebox(1,1){$\downarrow$}}
\put(6,0){\makebox(0,1)[l]{if
$\alpha < \beta \leq \gamma$ and $(\alpha,\gamma) \ne (x,\ol{x})$,}}
\end{picture}
\par\noindent
\setlength{\unitlength}{5.5mm}
\begin{picture}(22,5.5)(-3,0)
\put(-0.5,0){\makebox(0,2)[r]{case 3c}}
\put(0,0){
	\multiput(0,0)(0,1){3}{\line(1,0){1}}
	\multiput(0,0)(1,0){2}{\line(0,1){2}}
	\put(0,1){\makebox(1,1){$\alpha$}}
	\put(0,0){\makebox(1,1){$\gamma$}}
	}
\put(4,0){
	\multiput(0,0)(0,1){3}{\line(1,0){1}}
	\multiput(0,0)(1,0){2}{\line(0,1){2}}
	\put(0,1){\makebox(1,1){$\beta$}}
	\put(0,0){\makebox(1,1){$\gamma$}}
	}
\put(0,3){
	\multiput(0,0)(0,1){2}{\line(1,0){1}}
	\multiput(0,0)(1,0){2}{\line(0,1){1}}
	\put(0,0){\makebox(1,1){$\beta$}}
	}
\put(3,3){
	\multiput(0,0)(0,1){2}{\line(1,0){1}}
	\multiput(0,0)(1,0){2}{\line(0,1){1}}
	\put(0,0){\makebox(1,1){$\alpha$}}
	}
\put(1.5,0){\makebox(1,2){$=$}}
\put(0,2){\makebox(1,1){$\downarrow$}}
\put(3,2){\makebox(1,1){$\downarrow$}}
\put(6,0){\makebox(0,1)[l]{if
$\alpha \leq \beta < \gamma$ and $(\alpha,\gamma) \ne (x,\ol{x})$,}}
\end{picture}
\par\noindent
\setlength{\unitlength}{5.5mm}
\begin{picture}(22,5.5)(-3,0)
\put(-0.5,0){\makebox(0,2)[r]{case 4c}}
\put(0,0){
	\multiput(0,0)(0,1){3}{\line(1,0){1}}
	\multiput(0,0)(1,0){2}{\line(0,1){2}}
	\put(0,1){\makebox(1,1){$x$}}
	\put(0,0){\makebox(1,1){$\beta$}}
	}
\put(4,0){
	\multiput(0,0)(0,1){3}{\line(1,0){2}}
	\multiput(0,0)(2,0){2}{\line(0,1){2}}
	\put(0,1){\makebox(2,1){$x\!+\!1$}}
	\put(0,0){\makebox(2,1){$\ol{x\!+\!1}$}}
	}
\put(0,3){
	\multiput(0,0)(0,1){2}{\line(1,0){1}}
	\multiput(0,0)(1,0){2}{\line(0,1){1}}
	\put(0,0){\makebox(1,1){$\ol{x}$}}
	}
\put(3,3){
	\multiput(0,0)(0,1){2}{\line(1,0){1}}
	\multiput(0,0)(1,0){2}{\line(0,1){1}}
	\put(0,0){\makebox(1,1){$\beta$}}
	}
\put(1.5,0){\makebox(1,2){$=$}}
\put(0,2){\makebox(1,1){$\downarrow$}}
\put(3,2){\makebox(1,1){$\downarrow$}}
\put(8,0){\makebox(0,1)[l]{if
$x < \beta < \ol{x}$ and $x \ne n$,}}
\end{picture}
\par\noindent
\setlength{\unitlength}{5.5mm}
\begin{picture}(22,5.5)(-3,0)
\put(-0.5,0){\makebox(0,2)[r]{case 5c}}
\put(0,0){
	\multiput(0,0)(0,1){3}{\line(1,0){1}}
	\multiput(0,0)(1,0){2}{\line(0,1){2}}
	\put(0,1){\makebox(1,1){$x$}}
	\put(0,0){\makebox(1,1){$\ol{x}$}}
	}
\put(5,0){
	\multiput(0,0)(0,1){3}{\line(1,0){2}}
	\multiput(0,0)(2,0){2}{\line(0,1){2}}
	\put(0,1){\makebox(2,1){$\beta$}}
	\put(0,0){\makebox(2,1){$\ol{x\!-\!1}$}}
	}
\put(0,3){
	\multiput(0,0)(0,1){2}{\line(1,0){1}}
	\multiput(0,0)(1,0){2}{\line(0,1){1}}
	\put(0,0){\makebox(1,1){$\beta$}}
	}
\put(3,3){
	\multiput(0,0)(0,1){2}{\line(1,0){2}}
	\multiput(0,0)(2,0){2}{\line(0,1){1}}
	\put(0,0){\makebox(2,1){$x\!-\!1$}}
	}
\put(1.5,0){\makebox(1,2){$=$}}
\put(0,2){\makebox(1,1){$\downarrow$}}
\put(3,2){\makebox(2,1){$\downarrow$}}
\put(8,0){\makebox(0,1)[l]{if
$x \le \beta \le \ol{x}$ and $ x \ne 0$,}}
\end{picture}
\par\noindent
where
\begin{picture}(6.5,4.5)(-1,0)
\put(0,0){\makebox(1,1){$C$}}
\put(0,2){
	\multiput(0,0)(0,1){2}{\line(1,0){1}}
	\multiput(0,0)(1,0){2}{\line(0,1){1}}
	\put(0,0){\makebox(1,1){$\alpha$}}
	}
\put(3,2){
	\multiput(0,0)(0,1){2}{\line(1,0){1}}
	\multiput(0,0)(1,0){2}{\line(0,1){1}}
	\put(0,0){\makebox(1,1){$\beta$}}
	}
\put(4,0){\makebox(1,1){$C'$}}
\put(1.5,0){\makebox(1,1){$=$}}
\put(0,1){\makebox(1,1){$\downarrow$}}
\put(3,1){\makebox(1,1){$\downarrow$}}
\end{picture}
means that
if a letter $\beta$ is column inserted into a column $C'$ then
the column is changed to $C$ and a letter $\alpha$ is bumped out.

\subsection{\mathversion{bold}Main theorem : $C^{(1)}_n$ case}
\label{subsec:isorule}
Fix $l, k \in \Z_{\ge 1}$.
Given $b_1 \otimes b_2 \in B_{l} \otimes B_{k}$,
we define the element
$b'_2 \otimes b'_1 \in B_{k} \otimes B_{l}$
and $l',k', m \in \Z_{\ge 0}$ by the following rule.

\begin{rules}\label{rule:C}
\hfill\par\noindent
Set $z = \min(\sharp\,\fbx{0} \text{ in }\T(b_1),\, \sharp\,\fbx{0} \text{ in }\T(b_2))
=$
$\min(\sharp\,\fbx{\ol{0}} \text{ in }\T(b_1),\, \sharp\,\fbx{\ol{0}} \text{ in }\T(b_2))$.
Remove ( \fbx{0} , \fbx{\ol{0}} ) pairs simultaneously from $\T (b_1)$
and $\T (b_2)$ $z$ times.
Denote the resulting tableaux by $\hat{\T} (b_1)$ and $\hat{\T} (b_2)$, and set
$l' = \hat{\T} (b_1)_{1} = l-2z$ and $k' = \hat{\T} (b_2)_{1}=k-2z$.
($T_{1}$ is the length of the first row of a tableau $T$.)
Operate the column insertion and set
$\hat{\P} (b_2 \to b_1) = \left(\hat{\T} (b_2) \longrightarrow \hat{\T} (b_1)\right)$.
It has the form:

\setlength{\unitlength}{5mm}
\begin{picture}(20,4)
\put(8,1){\line(1,0){3.5}}
\put(8,2){\line(1,0){9}}
\put(8,3){\line(1,0){9}}
\put(8,1){\line(0,1){2}}
\put(11.5,1){\line(0,1){1}} 
\put(12.5,2){\line(0,1){1}} 
\put(17,2){\line(0,1){1}}
\put(12.5,2){\makebox(4.5,1){$i_{m+1} \;\cdots\; i_{l'}$}}
\put(8,1){\makebox(3,1){$\;\;i_1 \cdots i_m$}}
\put(8.5,2){\makebox(3,1){$\;\;j_1 \cdots\cdots j_{k'}$}}
\end{picture}

\noindent
where $m$ is the length of the second row, hence that of the first
row is $l'+k'-m$. ($0 \le m \le k'$.)

Next we  bump out  $l'$ letters from
the tableau $T^{(0)} = \hat{\P} (b_2 \to b_1)$ by the reverse bumping
algorithm.
For the boxes containing $i_{l'}, i_{l'-1}, \ldots, i_1$ in the above
tableau, we do it first for $i_{l'}$ then $i_{l'-1}$ and so on.
Correspondingly, let $w_{1}$ be the first letter that is  bumped out from
the leftmost column and $w_2$ be the second and so on.
Denote by $T^{(i)}$  the resulting tableau when $w_i$ is bumped out
($1 \le i \le l'$).
Note that $w_1 \le w_2 \le \cdots \le w_{l'}$.
Now $b'_1 \in B_l$ and $b'_2 \in B_k$ are uniquely specified by
\begin{eqnarray*}
\T(b'_2) &=&
\setlength{\unitlength}{5mm}
\begin{picture}(9.5,1.4)(0,0.3)
\multiput(0,0)(0,1){2}{\line(1,0){9}}
\multiput(0,0)(3,0){4}{\line(0,1){1}}
\put(0,0){\makebox(3,1){$0\cdots 0$}}
\put(3,0){\makebox(3,1){$T^{(l')}$}}
\put(6,0){\makebox(3,1){$\ol{0}\cdots \ol{0}$}}
\multiput(0,0.9)(6,0){2}{\put(0,0){\makebox(3,1){$z$}}}
\end{picture},\\
\T(b'_1) &=&
\begin{picture}(10.5,2)(0,0.3)
\multiput(0,0)(0,1){2}{\line(1,0){10}}
\multiput(0,0)(3,0){2}{\line(0,1){1}}
\multiput(4.25,0)(1.5,0){2}{\line(0,1){1}}
\multiput(7,0)(3,0){2}{\line(0,1){1}}
\put(0,0){\makebox(3,1){$0\cdots 0$}}
\put(3,0){\makebox(1.25,1){$w_{1}$}}
\put(4.25,0){\makebox(1.5,1){$\cdots$}}
\put(5.75,0){\makebox(1.25,1){$w_{l'}$}}
\put(7,0){\makebox(3,1){$\ol{0}\cdots \ol{0}$}}
\multiput(0,0.9)(7,0){2}{\put(0,0){\makebox(3,1){$z$}}}
\end{picture}.
\end{eqnarray*}
\end{rules}

\vskip3ex
Our main result for $U'_q(C^{(1)}_n)$ is
\begin{theorem}\label{th:main1}
Given $b_1 \ot b_2 \in B_l \ot B_k$, specify $b'_2 \ot b'_1 \in B_k \ot B_l$
and $l', k', m$ by Rule \ref{rule:C}.
Let $\iota: B_l \ot B_k \stackrel{\sim}{\rightarrow} B_k \ot B_l$ be the isomorphism of
$U'_q(C^{(1)}_n)$ crystal.
Then we have
\begin{align*}
\iota(b_1\otimes b_2)& = b'_2 \otimes b'_1,\\
H_{B_l B_k}(b_1 \otimes b_2) &= \min(l',k')- m.
\end{align*}
\end{theorem}

In  Appendix \ref{sec:anoth} we give an alternative algorithm
equivalent to Rule \ref{rule:C}, which is analogous
to the type $A$ case (Rule 3.11 of \cite{NY}).
In practical calculations it is often
more efficient than Rule \ref{rule:C}  based on the insertion algorithm.

\par\noindent
\begin{remark}
Associated with the tableau
$\hat{\P} (b_2 \to b_1)=
\left( \bbt (b_2) \longrightarrow \bbt (b_1) \right)$,
we have the recording tableau $\hat{\Q} (b_2 \to b_1)$,
as in the Robinson-Schensted-Knuth correspondence \cite{F}.
$\hat{\Q} (b_2 \to b_1)$ has a common shape with
$\hat{\P} (b_2 \to b_1)$, and its entries are the consecutive
integers from $1$ to $l'+k'$.
(Integers from $l'+1$ to $l'+ m$ are in the second row.)
With $\hat{\Q} (b_2 \to b_1)$, we can reverse the column insertion
procedure and recover $\bbt (b_1)$ and
$\bbt (b_2)$ from $\hat{\P} (b_2 \to b_1)$.
The recording tableau $\hat{\Q} (b'_1 \to b'_2)$ for the column insertion
$\left( \bbt (b'_1) \longrightarrow \bbt (b'_2) \right) $ is
similarly defined.
It has a common shape with
$\hat{\P} (b_2 \to b_1)$, and its entries are the consecutive
integers from $1$ to $l'+k'$.
(Integers from $k'+1$ to $k'+m$ are in the second row.)
In Rule \ref{rule:C} we have constructed
$T^{(l')}$  and \fbox{$w_{1} \cdots w_{l'}$}
from $\hat{\P} (b_2 \to b_1)$
with the help of the recording tableau $\hat{\Q} (b'_1 \to b'_2)$.
\end{remark}

\begin{example}
Let us assume  $n \ge 4$ and take
$b_1 = (2,0,1,1,0,\ldots,0,1,1,1) \in B_9$ and
$b_2 = (0,\ldots,0,3,0,0,2) \in B_7$.
Then we have $z = 1, l' = 7, k' = 5$ and
\begin{displaymath}
\bbt (b_1) = 1134\bar{3}\bar{2}\bar{1},\quad
\bbt (b_2) = \bar{4} \bar{4} \bar{4} \bar{1} \bar{1}.
\end{displaymath}
In this example we have
\begin{displaymath}
\hat{\P} (b_2 \to b_1)=\begin{array}{l}
	0034\bar{3}\bar{2}\bar{1} \\
	\bar{4}\bar{4}\bar{4}\bar{0}\bar{0}
	\end{array}, \>
\hat{\Q} (b_2 \to b_1)=\begin{array}{l}
	1234567 \\
	89\acute{0}\acute{1}\acute{2}
	\end{array}, \>
\hat{\Q} (b_1' \to b_2')=\begin{array}{l}
	12345\acute{1}\acute{2} \\
	6789\acute{0}
	\end{array}.
\end{displaymath}
Here we have written  $10,11,12$ as $\acute{0},\acute{1},\acute{2}$.
The column insertion
$\left( \bbt (b_2) \longrightarrow \bbt (b_1) \right)$
goes as
\begin{displaymath}
	\begin{array}{l}
	1134\bar{3}\bar{2}\bar{1} \\
	\bar{1}
	\end{array},\,
	\begin{array}{l}
	0134\bar{3}\bar{2}\bar{1} \\
	\bar{1}\bar{0}
	\end{array},\,
	\begin{array}{l}
	0134\bar{3}\bar{2}\bar{1} \\
	\bar{4}\bar{1}\bar{0}
	\end{array},\,
	\begin{array}{l}
	0034\bar{3}\bar{2}\bar{1} \\
	\bar{4}\bar{4}\bar{0}\bar{0}
	\end{array},\,
    \begin{array}{l}
	0034\bar{3}\bar{2}\bar{1} \\
	\bar{4}\bar{4}\bar{4}\bar{0}\bar{0}
	\end{array}.
\end{displaymath}
The reverse bumping according to the recording tableau
$\hat{\Q} (b_1' \to b_2')$ goes as
\begin{displaymath}
    \begin{array}{l}
	0034\bar{3}\bar{2}\bar{1} \\
	\bar{4}\bar{4}\bar{4}\bar{0}\bar{0}
	\end{array},\,
    \begin{array}{l}
	034\bar{3}\bar{2}\bar{1} \\
	\bar{4}\bar{4}\bar{4}\bar{0}\bar{0}
	\end{array},\,
    \begin{array}{l}
	034\bar{2}\bar{1} \\
	\bar{4}\bar{4}\bar{3}\bar{0}\bar{0}
	\end{array},\,
    \begin{array}{l}
	044\bar{2}\bar{1} \\
	\bar{4}\bar{4}\bar{0}\bar{0}
	\end{array},\,
    \begin{array}{l}
	044\bar{2}\bar{1} \\
	\bar{4}\bar{0}\bar{0}
	\end{array},\,
    \begin{array}{l}
	144\bar{2}\bar{1} \\
	\bar{1}\bar{0}
	\end{array},\,
    \begin{array}{l}
	144\bar{2}\bar{1} \\
	\bar{0}
	\end{array}.
\end{displaymath}
Adding the ( \fbx{0} , \fbx{\ol{0}} ) pair $z=1$ time, we get
$$
\T (b_2') = 0144\bar{2}\bar{1}\bar{0},\qquad
\T (b_1') =  00 \bar{4} \bar{4} \bar{4} \bar{4} \bar{1} \bar{0} \bar{0}.
$$
Therefore we obtain
$$b'_1 = (0,\ldots,0,4,0,0,1) \in B_9, \quad
b'_2 = (1,0,0,2,0,\ldots,0,1,1) \in B_7.$$
\end{example}

\section{\mathversion{bold}Proof : $C^{(1)}_n$ case}
\label{sec:theor}

We call an element $b$ of a $U'_q(C_n^{(1)})$ crystal
a {\itshape $\mathit{U_q(C_n)}$ highest  element} if it satisfies
$\etd_i b = 0$ for $i=1,2,\ldots ,n$.
Let $b'_2 \otimes b'_1 = \iota(b_1 \otimes b_2)$
under the isomorphism of $U_q'(C_n^{(1)})$ crystals
$\iota : B_{l} \otimes B_{k} \stackrel{\sim}{\rightarrow} B_{k} \otimes B_{l} $.
By definition if $b_1 \otimes b_2$ is a
$U_q(C_n)$ highest element so is $b'_2 \otimes b'_1$.
In Section \ref{subsec:thpartx} we prove Proposition \ref{th:main2}.
It verifies Theorem \ref{th:main1} when $b_1\ot b_2$ is a
$U_q(C_n)$ highest element and either $\T(b_1)$ or $\T(b_2)$ is free of
\fbx{0}.
In Section \ref{subsec:thpartxx} we prove that
if both $\T (b_1)$ and $\T (b_2)$ have at least one \fbx{0},
then the combinatorial $R$ on $B_l \ot B_k$ is reduced to
the combinatorial $R$ on $B_{l-2} \ot B_{k-2}$ by removing
a ( \fbx{0} , \fbx{\ol{0}} ) pair.
In Section \ref{subsec:thpartxxx} we quote a proposition \cite{B}
that assures the compatibility of the column insertion algorithm with a
$U_q(C_n)$ crystal isomorphism.
Based on these preparations we complete the proof of Theorem \ref{th:main1}
for general elements in Section \ref{subsec:thpartxxxx}.

\subsection{\mathversion{bold}Combinatorial $R$ for a class of highest elements}
\label{subsec:thpartx}

\begin{proposition}
\label{th:main2}
Given  $b_1 \ot b_2 \in B_l \ot B_k$, let
$b_2'\ot b'_1 = \iota(b_1 \ot b_2) \in B_k \ot B_l$
be the image under the isomorphism.
Suppose that $b_1 \ot b_2$ is a $U_q(C_n)$ highest element,
and
$\T (b_1)$ or $\T (b_2)$ has no \fbx{0}.
Then $\T (b'_2)$ or $\T (b'_1)$ also has no \fbx{0}, and
their column insertions give a common tableau:
\begin{equation}
\left( \T (b_2) \longrightarrow \T (b_1) \right) =
\left( \T (b'_1) \longrightarrow \T (b'_2) \right) .
\end{equation}
The value of the energy function is given by
$$H(b_1 \ot b_2) =
\left( \T (b_2) \longrightarrow \T (b_1) \right)_1 -\max (l,k).$$
\end{proposition}
\noindent
We give a proof of Proposition \ref{th:main2} by a case checking
in Appendix \ref{ap:proofiso}.
In this subsection we only list up all the
$U_q(C_n)$ highest elements of the above type.
We also list up the values of their energy functions.
Let $(x_1,x_2, \mbox{---},\ol{x}_1)$ stand for
$(x_1,x_2,0,\ldots ,0,\ol{x}_1)$ .
\begin{lemma}\label{lem:choohigh}
We have
$$
\iota : (l,0,\mbox{---},0) \ot (k,0,\mbox{---},0) \mapsto
(k,0,\mbox{---},0) \ot (l,0,\mbox{---},0)
$$
under the isomorphism
$\iota : B_{l} \otimes B_{k} \stackrel{\sim}{\rightarrow} B_{k} \otimes B_{l} $.
\end{lemma}
\begin{proof}
They are the unique elements in $B_l \ot B_k$ and $B_k \ot B_l$
respectively that do not vanish when $(\et{0})^{l+k}$ is applied.
\end{proof}

\begin{lemma}
Let $b_1 \ot b_2 \in B_l \ot B_k \, (l \geq k)$.
Suppose that $b_1 \ot b_2$ is a $U_q(C_n)$ highest element,
and
$\T (b_1)$ or $\T (b_2)$ has no \fbx{0}.
Then it has either the form:
$$(l,0,\mbox{---},0) \otimes (x_1,x_2,\mbox{---},\overline{x}_1)$$
with $x_1,x_2,\ol{x}_1 \in \Z_{\geq 0}$ and $x_1+x_2+\ol{x}_1 \leq k$,
or the form:
$$(l-2y_0,0,\mbox{---},0) \otimes (x_1,x_2,\mbox{---},k-x_1-x_2)$$
with $y_0(\neq 0),x_1,x_2 \in \Z_{\geq 0}$, $l-k \geq 2y_0-x_1$ and 
$x_1 + x_2 \le k$.
\end{lemma}
We call the former a {\itshape type I},
and the latter a {\itshape type II} $U_q(C_n)$ highest element.
They are exclusive.
\begin{proof}
Let $b_1 = (y_1,\ldots,y_n,\ol{y}_n,\ldots,\ol{y}_1)$
and $b_2 = (x_1,\ldots,x_n,\ol{x}_n,\ldots,\ol{x}_1)$.
Since $b_1 \ot b_2$ is a $U_q(C_n)$ highest element,
$\veps_i(b_1 \ot b_2) = \max (\veps_i(b_1),
\veps_i(b_1)+\veps_i(b_2)-\vphi_i(b_1) ) = 0$ for
$i=1,\ldots,n$.
It means that $\veps_i(b_1)=0$ and
$\vphi_i(b_1) \geq \veps_i(b_2)$ for $i=1,\ldots,n$.
Thus we have $\veps_n(b_1) = \ol{y}_n =0$, and then
$\veps_{n-1}(b_1) = \ol{y}_{n-1} +(y_n - \ol{y}_n )_+=0$,
i.e.~$\ol{y}_{n-1}=y_n=0$.
Repeating the same process we come to
$\veps_{1}(b_1) = \ol{y}_{1} +(y_{2} - \ol{y}_{2} )_+=0$,
i.e.~$\ol{y}_{1}=y_{2}=0$.
Thus $b_1 = (l-2y_0,0,\mbox{---},0)$ and $\vphi_i(b_1)=
(l-2y_0)\delta_{i,1}$.
Thus we have $\veps_n(b_2) = \ol{x}_n =0$, and then
$\veps_{n-1}(b_2) = \ol{x}_{n-1} +(x_n - \ol{x}_n )_+=0$,
i.e.~$\ol{x}_{n-1}=x_n=0$.
Repeating the same process we come to
$\veps_{2}(b_2) = \ol{x}_{2} +(x_{3} - \ol{x}_{3} )_+=0$,
i.e.~$\ol{x}_{2}=x_{3}=0$.
Thus $b_2 = (x_1,x_2,\mbox{---},\ol{x}_1)$ and
$\veps_1(b_2)=
x_2+\ol{x}_1$.
Therefore we have a condition
$\vphi_1(b_1)=l-2y_0 \geq x_2+\ol{x}_1$.
If $\T (b_1)$ has no \fbx{0} then $y_0=0$ and this condition
certainly holds.
If $\T (b_2)$ has no \fbx{0} then $x_2+\ol{x}_1=k-x_1$, thus
we impose the condition $l-k \geq 2y_0-x_1$.
\end{proof}
\begin{lemma}
\label{lem:3}
Under the isomorphism of $U_q'(C_n^{(1)})$ crystals
\begin{displaymath}
\iota : B_{l} \otimes B_{k} \stackrel{\sim}{\rightarrow} B_{k} \otimes B_{l}
         \quad (l \geq k),
\end{displaymath}
the {\itshape type I} $U_q(C_n)$ highest element
is mapped as
\begin{eqnarray*}
&&  (l,0,\mbox{---},0) \otimes (x_1,x_2,\mbox{---},\overline{x}_1) \\
&&           \mapsto (k,0,\mbox{---},0) \otimes
           (x_1+l-k-y,x_2,\mbox{---},\overline{x}_1-y),
\end{eqnarray*}
where $y=\min [l-k, (\overline{x}_1 - x_1)_+ ]$.
The value of the energy function for this element is
$x_0 + (x_1 - \overline{x}_1)_+$ with $x_0 =
(k-x_1-x_2-\overline{x}_1)/2$.
\end{lemma}
\begin{proof}
For a set of operators ${\mathcal O}_1 ,{\mathcal O}_2, \ldots$
on the $U_q'(C_n^{(1)})$ crystals,
we define ${\displaystyle
\prod_{i}^{2 \swarrow \nwarrow 1} {\mathcal O}_i}$ by
${\mathcal O}_2 {\mathcal O}_1$ for $n=2$,
${\mathcal O}_2 {\mathcal O}_3 {\mathcal O}_2 {\mathcal O}_1$
for $n=3$,
${\mathcal O}_2 {\mathcal O}_3 {\mathcal O}_4 {\mathcal O}_3
{\mathcal O}_2 {\mathcal O}_1$ for $n=4$ and so on.
Let $l=2m$ or $l=2m-1$.
The lemma can be proved by applying the following sequence of operators
\begin{equation}
\ft{0}^{m+x_0+(x_1-\ol{x}_1)_+}
\et{1}^{\min (x_1,\ol{x}_1)}
\left(
\prod_{i}^{2 \swarrow \nwarrow 1} \et{i}^{x_2 + \min (x_1,\ol{x}_1)}
\right)
\et{0}^{k+m}
\end{equation}
to the  both sides of Lemma \ref{lem:choohigh}.

In the sequel we will show
$$
H((l,0,\mbox{---},0)\ot(x_1,x_2,\mbox{---},\ol{x}_1)) =
H((l,0,\mbox{---},0)\ot(k,0,\mbox{---},0)) - k + x_0 + (x_1-\ol{x}_1)_+.
$$
In the case $l=2m$ and $m \geq k$, the value of the energy function
was lowered by $k$ when the first to the $k$-th $\et{0}$'s
were applied, and raised by $x_0 + (x_1 - \overline{x}_1)_+$
when the $(m+1)$-th to the last $\ft{0}$'s were applied.
In the case $l=2m$ and $m < k$, in addition to the same change
as in the previous case,
the value of the energy function
was raised by $k-m$ when the $(2m+1)$-th to the last $\et{0}$'s
were applied, and lowered by the same amount
when the first to the $(k-m)$-th $\ft{0}$'s were applied.

In the case $l=2m-1$ and $m-1 > k$, the value of the energy function
was lowered by $k$ when the first to the $k$-th $\et{0}$'s
were applied, and raised by $x_0 + (x_1 - \overline{x}_1)_+$
when the $(m+1)$-th to the last $\ft{0}$'s were applied.
In the case $l=2m-1$ and $m-1 \leq k$, in addition to the same change
as the previous case,
the value of the energy function
was raised by $k-m+1$ when the $2m$-th to the last $\et{0}$'s
were applied, and lowered by the same amount
when the first to the $(k-m+1)$-th $\ft{0}$'s were applied.

Recall that we have normalized the energy function as
$H_{B_lB_k}((l,0,\mbox{---},0) \ot (0,k,\mbox{---},0))=0$.
Thus we have 
$H((l,0,\mbox{---},0)\ot(x_1,x_2,\mbox{---},\ol{x}_1)) =
x_0 + (x_1-\ol{x}_1)_+$.
\end{proof}
\noindent

\begin{coro}\label{cor:hvalue}
For any $l,k \in \Zn$ we have
\begin{displaymath}
H_{B_lB_k}((l,0,\mbox{---},0) \ot (k,0,\mbox{---},0))=\min (l,k).
\end{displaymath}
\end{coro}

\begin{lemma}
\label{lem:4}
Under the isomorphism of $U_q'(C_n^{(1)})$ crystals
\begin{displaymath}
\iota : B_{l} \otimes B_{k} \stackrel{\sim}{\rightarrow} B_{k} \otimes B_{l}
         \quad (l \geq k),
\end{displaymath}
the {\itshape type II} $U_q(C_n)$ highest element
is mapped as
\begin{displaymath}
b_1 \ot b_2 := (l-2y_0,0,\mbox{---},0) \otimes (x_1,x_2,\mbox{---},k-x_1-x_2)
     \mapsto b'_2 \otimes b'_1,
\end{displaymath}
where $b'_2 \otimes b'_1$ is given by the following 1 and 2.
Let $\overline{x}_1 = k-x_1-x_2$.

\begin{enumerate}
	\item If $l-k > y_0 \geq x_1-\overline{x}_1$,
$$
b'_2 \otimes b'_1
 = (k,0,\mbox{---},0) \otimes
(x_1+l-k-y_0-z,x_2,\mbox{---},\overline{x}_1+y_0-z),
$$
where $z=\min [y_0+\overline{x}_1 - x_1,l-k-y_0 ]$.
$H(b_1\ot b_2) = 0$.
	\item If $l-k \leq y_0$ or $y_0 < x_1-\overline{x}_1$,
$$
b'_2 \otimes b'_1
= (k-2y_0+2w,0,\mbox{---},0) \otimes
(x_1+l-k-w,x_2,\mbox{---},\overline{x}_1+w),
$$
where $w=\min [l-k, (2y_0-x_1+\overline{x}_1)_+ ]$.
$H(b_1\ot b_2) = \max [y_0-l+k,x_1-\overline{x}_1-y_0]$.
\end{enumerate}
\end{lemma}

\begin{proof}
If $y_0 \geq x_1-\ol{x}_1$,
let $l=2m$ or $l=2m-1$.
The lemma can be proved by applying the following sequence of operators
\begin{equation}
\ft{0}^{m-y_0}
\et{1}^{x_1}
\left(
\prod_{i}^{2 \swarrow \nwarrow 1} \et{i}^{x_2 + x_1}
\right)
\et{0}^{k+m}
\end{equation}
to the both sides of
Lemma \ref{lem:choohigh}.
In the case $l=2m$ and $m \geq k$, the value of the energy function
was lowered by $k$ when the first to the $k$-th $\et{0}$'s
were applied.
In the case $l=2m$, $m < k$ and $2m-k>y_0$ (resp.~$2m-k \leq y_0$),
in addition to the same change in the previous case,
the value of the energy function
was raised by $k-m$ when the $(2m+1)$-th to the last $\et{0}$'s
were applied, and lowered by $k-m$ (resp.~$m-y_0$)
when the first to the $(k-m)$-th (resp.~the last) $\ft{0}$'s were applied.
In the case $l=2m-1$ and $m-1 > k$, the value of the energy function
was lowered by $k$ when the first to the $k$-th $\et{0}$'s
were applied.
In the case $l=2m-1$, $m-1 \leq k$ and $2m-1-k>y_0$
(resp.~$2m-1-k \leq y_0$),
in addition to the same change in the previous case,
the value of the energy function
was raised by $k-m+1$ when the $2m$-th to the last $\et{0}$'s
were applied, and lowered by $k-m+1$ (resp.~$m-y_0$)
when the first to the $(k-m+1)$-th (resp.~the last) $\ft{0}$'s were applied.

If $y_0 < x_1-\ol{x}_1 \leq 2 y_0$, one can check that
$\et{0}^{x_1-\ol{x}_1-y_0}(b_1 \ot b_2) =
b_1 \ot \et{0}^{x_1-\ol{x}_1-y_0}b_2$.
Lemma \ref{lem:1} and the previous case of the present lemma
enable us to obtain its image
under the map $\iota$.
They also tell us that now the value of the energy function is
equal to $(2y_0-x_1+\ol{x}_1+k-l)_+$.
Then apply $\ft{0}^{x_1-\ol{x}_1-y_0}$.
Since it again turns out to hit the right component of the tensor
product,
the value of the energy function is raised by $x_1-\ol{x}_1-y_0$.

If $2 y_0 < x_1-\ol{x}_1$, one can check that
$\et{0}^{x_1-\ol{x}_1-y_0}(b_1 \ot b_2) =
b_1 \ot \et{0}^{x_1-\ol{x}_1-y_0}b_2$.
Lemma \ref{lem:1} and \ref{lem:3} enable us to obtain its image
under the map $\iota$.
They also tell us that now the value of the energy function is
equal to $0$.
Then apply $\ft{0}^{x_1-\ol{x}_1-y_0}$.
Since it again hits the right component of the tensor
product,
the value of the energy function is raised by $x_1-\ol{x}_1-y_0$.
\end{proof}

\subsection{\mathversion{bold}Relation of $R$ on $B_l \ot B_k$ and $B_{l-2} \ot B_{k-2}$}
\label{subsec:thpartxx}
Let $l \ge 3$.
For any $b = (x_1, \ldots, \bar{x}_1) \in B_{l-2}$
we define $\tau^l_{l-2}(b) \in B_l$ to be
the unique element such that the tableau $\T (\tau^l_{l-2} (b))$ is made
from the tableau $\T (b)$ by adding a ( \fbx{0} , \fbx{\ol{0}} ) pair.
Note that $\tau_{l-2}^l(b)$ also has the same presentation
$(x_1,\ldots, \bar{x}_1)$ in $B_{l}$.
The map $\tau^l_{l-2} : B_{l-2} \rightarrow B_l$ is injective and has the property:
\begin{align}
\tilde{f}_i\tau^l_{l-2}(b) &= \tau^l_{l-2}(\tilde{f}_ib)\quad(0 \le i \le n)
\qquad \text{ if } \tilde{f}_ib \neq 0, \notag\\
\tilde{e}_i\tau^l_{l-2}(b) &= \tau^l_{l-2}(\tilde{e}_ib)\quad(0 \le i \le n)
\qquad \text{ if } \tilde{e}_ib \neq 0. \label{eq:tauandf}
\end{align}
\begin{lemma}
\label{lem:1}
We have
$\tau^l_{l-2}(b_1) \ot \tau^k_{k-2}(b_2) \simeq
\tau^k_{k-2}(b_2') \ot \tau^l_{l-2}(b_1')$ under the isomorphism
$B_{l} \ot B_{k} \simeq B_{k} \ot B_{l}$,
if and only if
$b_1 \ot b_2 \simeq b_2' \ot b_1'$
under  $B_{l-2} \ot B_{k-2} \simeq B_{k-2} \ot B_{l-2}$.
We also have
$H_{B_lB_k}(\tau^l_{l-2}(b_1) \ot \tau^k_{k-2}(b_2)) =
H_{B_{l-2}B_{k-2}}(b_1 \ot b_2)$.
\end{lemma}
\begin{proof}
Since $\tau^l_{l-2}$ and $\tau^k_{k-2}$ are injective,
the \textit{only if} part of the statement follows
immediately after when the \textit{if} part is proved.
Without loss of generality we assume $l \ge k$.
Set
$$b^{(l)}=(l,0,\mbox{---},0) \in B_l.$$
First consider the case $b_1=b^{(l-2)} \in B_{l-2}$ and
$b_2=b^{(k-2)} \in B_{k-2}$.
Then $b_2'=b^{(k-2)}$,
$b_1'=b^{(l-2)}$ and
$H_{B_{l-2}B_{k-2}}(b^{(l-2)} \ot b^{(k-2)}) = k-2$ by Corollary \ref{cor:hvalue}.
On the other hand we have
\begin{align*}
\psi\left(b^{(l)}\ot b^{(k)}\right) &=
\tau^{l}_{l-2}(b^{(l-2)}) \ot \tau^k_{k-2}(b^{(k-2)}) \in B_l \ot B_k,\\
\psi\left(b^{(k)}\ot b^{(l)}\right) &=
\tau^{k}_{k-2}(b^{(k-2)}) \ot \tau^l_{l-2}(b^{(l-2)}) \in B_k \ot B_l,
\end{align*}
where
\begin{eqnarray*}
\psi &=& \etd_0
(\etd_1)^{l+k-2}(\etd_2)^{l+k-2}\cdots(\etd_{n-1})^{l+k-2}
(\etd_{n})^{l+k-2} \nonumber\\
&& \times (\etd_{n-1})^{l+k-2}\cdots
(\etd_2)^{l+k-2}(\etd_1)^{l+k-2}(\etd_0)^{l+k-1}.
\end{eqnarray*}
By Lemma \ref{lem:choohigh} one has
\begin{equation}\label{eq:tautau}
\tau^{l}_{l-2}(b^{(l-2)}) \ot \tau^k_{k-2}(b^{(k-2)}) \simeq
\tau^{k}_{k-2}(b^{(k-2)}) \ot \tau^l_{l-2}(b^{(l-2)})
\end{equation}
under the isomorphism $B_l \ot B_k \simeq B_k \ot B_l$.
The energy was lowered by $k$ when the first to the $k$-th
$\etd_0$'s were
applied, and raised by $k-1$ when the $(l+1)$-th to the
$(l+k-1)$-th
$\etd_0$'s were applied.
Then it was lowered by $1$ when the leftmost $\etd_0$ was
applied.
Thus we have
\begin{align}
H_{B_lB_k}(\tau^l_{l-2}(b^{(l-2)}) \ot \tau^k_{k-2}(b^{(k-2)})) 
&= H_{B_lB_k}(b^{(l)}\ot b^{(k)}) - 2 = k-2 \notag \\
& = H_{B_{l-2}B_{k-2}}(b^{(l-2)} \ot b^{(k-2)}).\label{eq:hhh}
\end{align}
The proof is finished in this special case from
Corollary \ref{cor:hvalue}.

Now we consider the general elements
$b_1 \otimes b_2 \in B_{l-2} \otimes B_{k-2}$ and
$b'_2 \otimes b'_1 \in B_{k-2} \otimes B_{l-2}$
that are mapped to each other under the isomorphism.
Take any finite sequence $\psi'$  made of $\etd_i$'s
and $\ftd_i$'s
$(i=0,1,\ldots,n)$ such that
\begin{equation}\label{eq:psiprime}
b_1\ot b_2 = \psi'(b^{(l-2)}\ot b^{(k-2)}),
\end{equation}
which is equivalent to
\begin{equation}\label{eq:psiprime2}
b'_2\ot b'_1 = \psi'(b^{(k-2)}\ot b^{(l-2)}).
\end{equation}
For any operator in $\psi'$, the rules
(\ref{eq:ot-e})-(\ref{eq:ot-f}) determine
whether it should hit  the left or the right component of the
tensor product.
For any $c_1 \ot c_2 \in B_{l-2} \otimes B_{k-2}$
we have
$\vphi_i (\tau^l_{l-2} (c_1)) = \vphi_i (c_1) +
\delta_{i,0}$
and
$\veps_i (\tau^k_{k-2} (c_2)) = \veps_i (c_2) +
\delta_{i,0}$ from (\ref{eq:eps}).
Thus the alternatives in
(\ref{eq:ot-e})-(\ref{eq:ot-f})
are not changed by $\tau^l_{l-2} \ot \tau^k_{k-2}$.
{}From (\ref{eq:tauandf}) it follows that
$(\tau^l_{l-2}\ot \tau^k_{k-2})(\psi'(c_1\ot c_2))  =
\psi'(\tau^l_{l-2}(c_1) \ot \tau^k_{k-2}(c_2))$.
Applying $\tau^l_{l-2}\ot \tau^k_{k-2}$
(resp.  $\tau^k_{k-2}\ot \tau^l_{l-2}$ )
to (\ref{eq:psiprime}) (resp. (\ref{eq:psiprime2}))
we thus get
\begin{align}
\tau^l_{l-2}(b_1)\ot \tau^k_{k-2}(b_2) &= \psi'(
\tau^l_{l-2}(b^{(l-2)})\ot \tau^k_{k-2}(b^{(k-2)})),
\label{eq:ppp1}\\
\tau^k_{k-2}(b'_2)\ot \tau^l_{l-2}(b'_1) &= \psi'(
\tau^k_{k-2}(b^{(k-2)})\ot \tau^l_{l-2}(b^{(l-2)})).
\label{ppp2}
\end{align}
{}From (\ref{eq:tautau}) it follows that
\[
\tau^l_{l-2}(b_1)\ot \tau^k_{k-2}(b_2) \simeq
\tau^k_{k-2}(b'_2)\ot \tau^l_{l-2}(b'_1)
\]
under the isomorphism $B_{l} \otimes B_{k} \simeq B_{k}
\otimes B_{l}$.
When comparing (\ref{eq:psiprime}) and (\ref{eq:ppp1})
change of the value of the
energy function caused by $\psi'$
is not affected by $\tau^l_{l-2} \ot \tau^k_{k-2}$.
Therefore from (\ref{eq:hhh})  we have $H_{B_lB_k}(\tau^l_{l-2}(b_1) \ot
\tau^k_{k-2}(b_2)) =
H_{B_{l-2}B_{k-2}}(b_1 \ot b_2)$.
\end{proof}
\subsection{\mathversion{bold}Column insertion and $U_q(C_n)$ crystal morphism}
\label{subsec:thpartxxx}
The next proposition is due to Baker (Proposition 7.1 of \cite{B}).
For a dominant integral weight $\lambda$ of the $C_n$ root system,
let $B(\lambda)$ be the $U_q(C_n)$ crystal
associated with the irreducible highest weight representation $V(\lambda)$ \cite{KN}.
The elements of $B(\lambda)$ can be represented by
the semistandard $C$-tableaux of shape $\lambda$ \cite{KN}.
\begin{proposition}
\label{lem:2}
Let $B(\mu) \ot B(\nu) \simeq \bigoplus_j B(\lambda_j)^{\oplus m_j}$
be the tensor product decomposition of crystals. Here $\la_j$'s are
distinct highest weights and $m_j(\ge1)$ is the multiplicity of $B(\la_j)$.
Forgetting the multiplicities we have the canonical morphism from
$B(\mu) \ot B(\nu)$ to $\bigoplus_j B(\lambda_j)$.
Define $\psi_C$ by
\begin{displaymath}
\psi_C (b_1 \ot b_2) = \left( b_2 \stackrel{\ast}{\longrightarrow}
b_1 \right).
\end{displaymath}
Then $\psi_C$ gives the unique crystal morphism from
$B(\mu) \ot B(\nu)$ to $\bigoplus_j B(\lambda_j)$.
\end{proposition}
\par\noindent
Here $b_2 \stackrel{\ast}{\longrightarrow} b_1$ is
the tableau obtained from successive column insertions of letters of
the Japanese reading word of $b_2$ into $b_1$ by the original definition
in \cite{B}
(In \cite{B}, $b_2 \stackrel{\ast}{\longrightarrow} b_1$ is
denoted by $b_1 * b_2$.)
\par\noindent
\begin{remark}\label{rem:avoid}
The insertion $b_2 \stackrel{\ast}{\longrightarrow} b_1$ may
include such a process that \fbx{x} and \fbx{\ol{x}} annihilate
pairwise and an empty box thereby produced slides out.
In \cite{B} this process was called a \textit{bumping-sliding transition}.
Consider the case that both $b_1$ and $b_2$ are one-row tableaux.
(We shall omit the symbol $\T$ here.)
In this case the bumping-sliding transition can occur only when
\fbx{\ol{1}} is inserted into \fbx{1}.
We defined our column insertion
$(\longrightarrow)$ so that
$\varsigma((b_2 \longrightarrow b_1))$ is equivalent to
$(\varsigma(b_2) \stackrel{\ast}{\longrightarrow} \varsigma(b_1))$.
In Rule \ref{rule:C} for $U_q'(C_n^{(1)})$ combinatorial $R$ matrix
we have removed ( \fbx{0} , \fbx{\ol{0}} ) pairs beforehand which
become ( \fbx{1} , \fbx{\ol{1}} ) under $\varsigma$.
Thus we have avoided the bumping-sliding
transition to occur.
\end{remark}

\subsection{Proof of Theorem \ref{th:main1}}
\label{subsec:thpartxxxx}

With no loss of generality we assume $l \geq k$.
Let $b'_2 \otimes b'_1$ be the image of $b_1 \otimes b_2$
under the isomorphism of $U_q'(C_n^{(1)})$ crystals
$\iota : B_{l} \otimes B_{k} \stackrel{\sim}{\rightarrow} B_{k} \otimes B_{l}$.
In order to prove Theorem \ref{th:main1} we are to show the claims:
\par\noindent
\begin{enumerate}
\item
Let
$$z_0 = \min( \sharp \fbx{0} \text{ in } \T (b_1),
\sharp \fbx{0} \text{ in } \T (b_2)),\;
z'_0 = \min( \sharp \fbx{0} \text{ in } \T (b'_1),
\sharp \fbx{0} \text{ in } \T (b'_2)).
$$
Then $z'_0 = z_0$.
\item
Remove $(\, \fbx{0}\, , \fbx{\ol{0}}\, )$ pairs $z_0$ times from
$\T (b_1)$, $\T (b_2)$, $\T (b'_1)$ and $\T (b'_2)$.
Call the resulting tableaux
$\bbt (b_1)$, $\bbt (b_2)$, $\bbt (b'_1)$ and $\bbt (b'_2)$, respectively.
Then we have
\begin{equation}
\left( \bbt (b_2) \longrightarrow \bbt (b_1) \right) =
\left( \bbt (b'_1) \longrightarrow \bbt (b'_2) \right) .
\end{equation}
\item
$
H_{B_l B_k}(b_1 \ot b_2) = \left(
\bbt (b_2) \longrightarrow \bbt (b_1) \right)_1 - \bbt (b_1)_1.
$
\end{enumerate}

\begin{proof}
Thanks to Lemma \ref{lem:1} it suffices to verify the above
claims only when $\T(b_1)$ or $\T(b_2)$ has no \fbx{0}.
Such a case can be reduced to
Proposition \ref{th:main2} by the argument as follows.

Let $b_1 \ot b_2$ be an element of $B_{l} \ot B_{k}$,
which is not necessarily a $U_q(C_n)$ highest element and
either $\T (b_1)$ or $\T (b_2)$ is free of \fbx{0}.
Let $b'_2 \ot b'_1 = \iota (b_1 \ot b_2)$
under the isomorphism
$\iota : B_{l} \ot B_{k} \to B_{k} \ot B_{l}$.
There exists a sequence
$i_1,i_2,\ldots,i_s \,(1 \leq i_\alpha \leq n,\,\alpha = 1,2,\ldots,s)$
such that
$\dot{b}_1 \ot \dot{b}_2 := \et{i_s} \cd \et{i_1} (b_1 \ot b_2)$
is a $U_q(C_n)$ highest element.
Then we have $b'_2 \ot b'_1 = \ft{i_1} \cd \ft{i_s} \circ \iota \circ
\et{i_s} \cd \et{i_1} (b_1 \ot b_2)$.
Let
$\dot{b}'_2 \ot \dot{b}'_1 = \iota (\dot{b}_1 \ot \dot{b}_2)$.
Since $\et{i} \, (1 \leq i \leq n)$ does not change
$\sharp$\fbx{0},
$\T (\dot{b}_1)$ or $\T (\dot{b}_2)$ also has no \fbx{0}.
Therefore from Proposition \ref{th:main2} we have
\begin{equation}\label{eq:46}
\left( \T (\dot{b}_2) \longrightarrow \T (\dot{b}_1) \right) =
\left( \T (\dot{b}'_1) \longrightarrow \T (\dot{b}'_2) \right),
\end{equation}
where $\T (\dot{b}'_1)$ or $\T (\dot{b}'_2)$ has no \fbx{0}.
Since $b'_2 \ot b'_1 = \ft{i_1} \cd \ft{i_s}(\dot{b}'_2 \ot \dot{b}'_1)$
and $1 \le i_\alpha \le n$, we conclude that
$\T(b'_1)$ or $\T(b'_2)$ has no \fbx{0}.
Thus Claim 1 is indeed valid as $z_0 = z'_0 = 0$.
By Remark \ref{re:bakertonokankei}, (\ref{eq:46}) is equivalent to
\begin{equation}
\left( \varsigma (\T (\dot{b}_2)) \stackrel{\ast}{\longrightarrow}
\varsigma (\T (\dot{b}_1)) \right) =
\left( \varsigma (\T (\dot{b}'_1)) \stackrel{\ast}{\longrightarrow}
\varsigma (\T (\dot{b}'_2)) \right) =: \P.
\end{equation}
Here  $\varsigma$ is defined in (\ref{eq:embedding}).
Regarding $\P$ as an element of a $U_q(C_{n+1})$ crystal, we
apply  Proposition \ref{lem:2}, to get
\begin{equation}
\left( \varsigma (\T (b_2)) \stackrel{\ast}{\longrightarrow}
\varsigma (\T (b_1)) \right) = \ft{i_1+1} \cd \ft{i_s+1} (\P)
= \left( \varsigma (\T (b'_1)) \stackrel{\ast}{\longrightarrow}
\varsigma (\T (b'_2)) \right),
\end{equation}
which is equivalent to
\begin{equation}
\left( \T (b_2) \longrightarrow \T (b_1) \right) =
\left( \T (b'_1) \longrightarrow \T (b'_2) \right),
\end{equation}
showing Claim 2.
Since $\ft{i} \, (1 \leq i \leq n)$ does not change
the shape of the tableaux \cite{B} and 
$H(b_1 \ot b_2) = H(\dot{b}_1 \ot \dot{b}_2)$, Claim 3 
follows from  Proposition \ref{th:main2}.
\end{proof}


\section{\mathversion{bold}$U_q'(A_{2n-1}^{(2)})$ crystal case}
\label{sec:discuss}
\subsection{Definitions}
Given a non-negative integer $l$,
let us denote by $B_{l}$ the $U_q'(A_{2n-1}^{(2)})$ crystal
defined in \cite{KKM}.
(Their $B_l$ is identical to our $B_{l}$.)
As a set $B_{l}$ reads
$$
B_{l} = \left\{(
x_1,\ldots, x_n,\overline{x}_n,\ldots,\overline{x}_1) \Biggm|
x_i, \overline{x}_i \in \Z_{\ge 0},
\sum_{i=1}^n(x_i + \overline{x}_i) = l \right\}.
$$
$B_{l}$ is isomorphic to
$B(l \La_1)$ as a crystal for $U_q(C_n)$.
The crystal structure is given by
{\allowdisplaybreaks
\begin{eqnarray}
\et{0} b &=&
\begin{cases}
	(x_1,x_2-1,\ldots,\overline{x}_2,\overline{x}_1+1)
	& \mbox{ if }x_2 > \overline{x}_2, \\
	(x_1-1,x_2,\ldots,\overline{x}_2+1,\overline{x}_1)
	& \mbox{ if } x_2 \leq  \overline{x}_2,
\end{cases} \nonumber\\
\et{n} b &=& (x_1,\ldots,x_n+1,\overline{x}_n-1,\ldots,\overline{x}_1),
\nonumber\\
\et{i} b &=&
\begin{cases}
	(x_1,\ldots,x_i+1,x_{i+1}-1,\ldots,\overline{x}_1)
	& \mbox{ if } x_{i+1} > \overline{x}_{i+1},\\
	(x_1,\ldots,\overline{x}_{i+1}+1,\overline{x}_i-1,\ldots,\overline{x}_1)
	& \mbox{ if } x_{i+1} \le \overline{x}_{i+1},
\end{cases}\nonumber\\
\ft{0} b &=&
\begin{cases}
	(x_1,x_2+1,\ldots,\overline{x}_2,\overline{x}_1-1)
	& \mbox{ if }x_2 \ge \overline{x}_2, \\
	(x_1+1,x_2,\ldots,\overline{x}_2-1,\overline{x}_1)
	& \mbox{ if } x_2 < \overline{x}_2,
\end{cases}
\nonumber\\
\ft{n} b &=& (x_1,\ldots,x_n-1,\overline{x}_n+1,\ldots,\overline{x}_1),
\nonumber\\
\ft{i} b &=&
\begin{cases}
	(x_1,\ldots,x_i-1,x_{i+1}+1,\ldots,\overline{x}_1)
	& \mbox{ if } x_{i+1} \ge \overline{x}_{i+1},\\
	(x_1,\ldots,\overline{x}_{i+1}-1,\overline{x}_i+1,\ldots,\overline{x}_1)
	& \mbox{ if } x_{i+1} < \overline{x}_{i+1},
\end{cases}
\end{eqnarray}
}
where $b= (x_1, \ldots, x_n, \overline{x}_n,\ldots,\overline{x}_1)$
and $i = 1, \ldots, n-1$.
For this $b$ we have
\begin{eqnarray}
\vphi_0 (b) &=& \ol{x}_1 + (\ol{x}_2-x_2)_+ , \quad
\veps_0 (b) = x_1 + (x_2-\ol{x}_2)_+ ,\nonumber\\
\vphi_i (b) &=& x_i + (\ol{x}_{i+1}-x_{i+1})_+
\quad \mbox{for} \, i=1,\ldots,n-1, \nonumber\\
\veps_i (b) &=& \ol{x}_i + (x_{i+1}-\ol{x}_{i+1})_+
\quad \mbox{for} \, i=1,\ldots,n-1, \nonumber\\
\vphi_n (b) &=& x_n , \quad \veps_n (b) = \ol{x}_n.
\end{eqnarray}
We shall depict the element
$b= (x_1, \ldots, x_n, \overline{x}_n,\ldots,\overline{x}_1) \in
B_{l}$ with the tableau:
\begin{equation}
{\mathcal T} (b)=\overbrace{\fbox{$\vphantom{\ol{1}} 1 \cd 1$}}^{x_1}\!
\fbox{$\vphantom{\ol{1}}\cd$}\!
\overbrace{\fbox{$\vphantom{\ol{1}}n \cd n$}}^{x_n}\!
\overbrace{\fbox{$\vphantom{\ol{1}}\ol{n} \cd \ol{n}$}}^{\ol{x}_n}\!
\fbox{$\vphantom{\ol{1}}\cd$}\!
\overbrace{\fbox{$\ol{1} \cd \ol{1}$}}^{\ol{x}_1}.
\end{equation}
The length of this one-row tableau is equal to $l$, namely
$\sum_{i=1}^n(x_i + \overline{x}_i) =l$.

In this section we normalize the energy function as
\begin{equation}\label{eq:kikaku}
H_{B_l B_k}((l,0,\mbox{---},0)\ot (0,\mbox{---},0,k) ) = 0,
\end{equation}
irrespective of $l < k$ or $l \ge k$.
\subsection{\mathversion{bold}Main theorem : $A^{(2)}_{2n-1}$ case}
\label{subsec:isoruleaa}
The insertion symbol $\longrightarrow$ and the reverse bumping
in this section are the same ones as in Section \ref{sec:descri}.

Given $b_1 \otimes b_2 \in B_{l} \otimes B_{k}$,
we define the element
$b'_2 \otimes b'_1 \in B_{k} \otimes B_{l}$
and $l',k', m \in \Z_{\ge 0}$ by the following rule.

\begin{rules}\label{rule:A}
\hfill\par\noindent
Set $z = \min(\sharp\,\fbx{1} \text{ in }{\mathcal T}(b_1),\,
\sharp\,\fbx{\ol{1}} \text{ in }{\mathcal T}(b_2))$.
Remove \fbx{1}'s (resp. \fbx{\ol{1}}'s) from
${\mathcal T}(b_1)$ (resp. ${\mathcal T}(b_2)$) $z$ times and
call  the resulting tableaux $\acute{\mathcal T} (b_1)$
(resp.$\grave{\mathcal T} (b_2)\,$).
Let $l' = \acute{\mathcal T} (b_1)_1 = l-z$
and $k' = \grave{\mathcal T} (b_2)_1 = k-z$.
Operate the column insertion
and set
$\P (b_2 \stackrel{\ast}{\to} b_1) = \left(
\grave{\mathcal T} (b_2) \longrightarrow
\acute{\mathcal T} (b_1) \right)$.
(This $\P (b_2 \stackrel{\ast}{\to} b_1)$ coincides with
the column insertion ${\mathcal T} (b_2)
\stackrel{\ast}{\longrightarrow} {\mathcal T} (b_1)$, because of
$\left( \fbx{\ol{1}} \stackrel{\ast}{\longrightarrow}
\fbx{1} \right) = \emptyset$.)
$\P (b_2 \stackrel{\ast}{\to} b_1)$ has the form:

\setlength{\unitlength}{5mm}
\begin{picture}(20,4)
\put(8,1){\line(1,0){3.5}}
\put(8,2){\line(1,0){9}}
\put(8,3){\line(1,0){9}}
\put(8,1){\line(0,1){2}}
\put(11.5,1){\line(0,1){1}} 
\put(12.5,2){\line(0,1){1}} 
\put(17,2){\line(0,1){1}}
\put(12.5,2){\makebox(4.5,1){$i_{m+1} \;\cdots\; i_{l'}$}}
\put(8,1){\makebox(3,1){$\;\;i_1 \cdots i_m$}}
\put(8.5,2){\makebox(3,1){$\;\;j_1 \cdots\cdots j_{k'}$}}
\end{picture}

\noindent
where $m$ is  the length of the second row, hence that of the first
row is $l'+k'-m$. ($0 \le m \le k'$.)

Next we bump out  $l'$ letters from
the tableau $T^{(0)} = \P (b_2 \stackrel{\ast}{\to} b_1)$ by the reverse bumping
algorithm.
For the boxes containing $i_{l'}, i_{l'-1}, \ldots, i_1$ in the above
tableau, we do it first for $i_{l'}$ then $i_{l'-1}$ and so on.
Correspondingly, let $w_{1}$ be the first letter that is bumped out from
the leftmost column and $w_2$ be the second  and so on.
Denote by $T^{(i)}$  the resulting tableau when $w_i$ is bumped out
($1 \le i \le l'$).
Now $b'_1 \in B_l$ and $b'_2 \in B_k$ are uniquely specified by
\[
{\mathcal T}(b'_2) =
\setlength{\unitlength}{5mm}
\begin{picture}(6.5,1.4)(0,0.3)
\multiput(0,0)(0,1){2}{\line(1,0){6}}
\put(0,0){\line(0,1){1}}
\put(3,0){\line(0,1){1}}
\put(6,0){\line(0,1){1}}
\put(3,0){\makebox(3,1){$T^{(l')}$}}
\put(0,0){\makebox(3,1){$1\cdots 1$}}
\put(0,0.9){\makebox(3,1){$z$}}
\end{picture},\;\;
{\mathcal T}(b'_1) =
\setlength{\unitlength}{5mm}
\begin{picture}(7.5,1.4)(0,0.3)
\multiput(0,0)(0,1){2}{\line(1,0){7}}
\put(0,0){\line(0,1){1}}
\put(1.25,0){\line(0,1){1}}
\put(2.75,0){\line(0,1){1}}
\put(4,0){\line(0,1){1}}
\put(7,0){\line(0,1){1}}
\put(0,0){\makebox(1.25,1){$w_{1}$}}
\put(1.25,0){\makebox(1.5,1){$\cdots$}}
\put(2.75,0){\makebox(1.25,1){$w_{l'}$}}
\put(4,0){\makebox(3,1){$\ol{1}\cdots\ol{1}$}}
\put(4,0.9){\makebox(3,1){$z$}}
\end{picture}.
\]
\end{rules}
\vskip3ex

Our main result for $U'_q(A^{(2)}_{2n-1})$ is
\begin{theorem}\label{th:main3}
Given $b_1 \ot b_2 \in B_l \ot B_k$, specify $b'_2 \ot b'_1 \in B_k \ot B_l$
and $l', k', m$ by Rule \ref{rule:A}.
Let $\iota: B_l \ot B_k \stackrel{\sim}{\rightarrow} B_k \ot B_l$ be the isomorphism of
$U'_q(A^{(2)}_{2n-1})$ crystal.
Then we have
\begin{align*}
\iota(b_1\otimes b_2)& = b'_2 \otimes b'_1,\\
H_{B_l B_k}(b_1 \otimes b_2) &= 2\min(l',k')- m.
\end{align*}
\end{theorem}

\par\noindent
\begin{example}
If $1123 \ot 1\bar{1}\bar{1}$ is regarded as an element of
$U'_q(A_{2n-1}^{(2)})$ crystal $B_{4} \ot B_{3}$,
it is mapped to $113 \ot 12\bar{1}\bar{1} \in B_{3} \ot B_{4}$
under the isomorphism.
Here $\P(1\bar{1}\bar{1} \stackrel{\ast}{\to} 1123) =123$ and
$H(1123 \ot 1\bar{1}\bar{1})=2$.
If $1123 \ot 1\bar{1}\bar{1}$ is regarded as an element of
$U'_q(C_n^{(1)})$ crystal $B_{4} \ot B_{3}$,
it is mapped to $123 \ot 01\bar{1}\bar{0} \in B_{3} \ot B_{4}$
under the isomorphism.
Here $\hat{\P}(1\bar{1}\bar{1} \to 1123
 ) ={0 \atop 1}{1 \atop \bar{1}}{2 \atop \bar{0}}
{3 \atop \vphantom{1}}$ and
$H(1123 \ot 1\bar{1}\bar{1})=0$.
\end{example}
\subsection{\mathversion{bold}Proof : $A^{(2)}_{2n-1}$ case}
Given $b_1\ot b_2 \in B_k\ot B_k$, determine $b'_2\ot b'_1 \in B_k\ot B_l$
by Rule \ref{rule:A}.
To prove Theorem \ref{th:main3}, we are to show the following claims:
\begin{eqnarray}
1.&
\left(\grave{\mathcal T} (b_2) \longrightarrow
\acute{\mathcal T} (b_1) \right) =
\left(\grave{\mathcal T} (b'_1) \longrightarrow
\acute{\mathcal T} (b'_2)\right).
\label{eq:astcolumins}\\
2.&
H_{B_lB_k}(b_1 \ot b_2)=\left(\grave{\mathcal T} (b_2) \longrightarrow
\acute{\mathcal T} (b_1) \right)_1 - \vert l-k \vert.
\label{eq:hval}
\end{eqnarray}

\par\noindent
\begin{lemma}\label{lem:choohi2}
We have 
$$
\iota : (l,0,\mbox{---},0) \ot (k,0,\mbox{---},0) 
 \mapsto 
(k,0,\mbox{---},0) \ot (l,0,\mbox{---},0)
$$ 
under the isomorphism
$B_{l} \otimes B_{k} \stackrel{\sim}{\rightarrow} B_{k} \otimes B_{l}$.
\end{lemma}
\begin{proof}
They are the unique elements in $B_l \ot B_k$ and $B_k \ot
B_l$
respectively that do not vanish when $(\et{0})^{l+k}$ is
applied
and do not vanish when $(\ft{1})^{l+k}$ is applied.
\end{proof}

\begin{proof}[Proof of Theorem \ref{th:main3}]
Claim 1 is due to Proposition \ref{lem:2} and the fact that
the irreducible decomposition of the
$U_q(C_n)$ module $V(l\Lambda_1) \ot V(k\Lambda_1)$ is
multiplicity-free (for generic $q$).

We call an element $b$ of a $U'_q(A_{2n-1}^{(2)})$ crystal
a {\itshape $\mathit{U_q(C_n)}$ highest element} if
it satisfies
$\etd_i b = 0$ for $i=1,2,\ldots ,n$.
To show Claim 2, it suffices to check it for $U_q(C_n)$
highest elements.
Then the general case follows from  Proposition \ref{lem:2}
because $H(\tilde{f}_{i_1}\cdots \tilde{f}_{i_j}(b_1\ot b_2)) =
H(b_1 \ot b_2)$ for $i_1, \ldots, i_j \in \{1,\ldots, n\}$
if $\tilde{f}_{i_1}\cdots \tilde{f}_{i_j}(b_1\ot b_2) \neq 0$.
We assume  $l \geq k$ with no loss of generality.
Suppose that $b_1 \ot b_2 \simeq b_2' \ot b_1'$
is a $U_q(C_n)$ highest element.
In general  it has the form:
$$b_1\ot b_2 = (l,0,\mbox{---},0) \ot (x_1,x_2,\mbox{---},\ol{x}_1),$$
where $x_1, x_2$ and $\ol{x}_1$ are arbitrary
as long as $k = x_1 + x_2 + \ol{x}_1$.
Applying
\begin{displaymath}
\et{0}^{\ol{x}_1} \et{2}^{x_2+\ol{x}_1} \cdots
\et{n-1}^{x_2+\ol{x}_1} \et{n}^{x_2+\ol{x}_1}
\et{n-1}^{x_2+\ol{x}_1}
\cdots \et{2}^{x_2+\ol{x}_1} \et{0}^{x_2+\ol{x}_1}
\end{displaymath}
to the both sides of Lemma \ref{lem:choohi2}, 
we find
\begin{displaymath}
(l,0,\mbox{---},0) \ot (x_1,x_2,\mbox{---},\ol{x}_1)
\simeq
(k,0,\mbox{---},0) \ot (x_1',x_2,\mbox{---},\ol{x}_1).
\end{displaymath}
Here  $x_1' = l-x_2-\ol{x}_1$.
In the course of the application of $\tilde{e}_i$'s,
the value of the energy function has changed as
$$
H\left((l,0,\mbox{---},0) \ot (x_1,x_2,\mbox{---},\ol{x}_1)\right) =
H\left((l,0,\mbox{---},0) \ot (k,0,\mbox{---},0)\right) - x_2 - 2 \ol{x}_1.
$$
Thus according to our normalization (\ref{eq:kikaku}) we have
$H(b_1 \ot b_2)=2k-x_2 - 2 \ol{x}_1$.
On the other hand
for this highest element
the column insertions (\ref{eq:astcolumins}) lead to
a tableau
\setlength{\unitlength}{5mm}
\begin{picture}(4,2)
\put(0,2){\line(1,0){4}}
\put(0,1){\line(1,0){4}}
\put(0,0){\line(1,0){3}}
\put(0,0){\line(0,1){2}}
\put(3,0){\line(0,1){1}}
\put(4,1){\line(0,1){1}}
\put(0,1){\makebox(4,1){$1\cdots\cdots 1$}}
\put(0,0){\makebox(3,1){$2\cdots 2$}}
\end{picture}
whose first row has the length  $l+k-x_2-2\ol{x}_1$.
This completes the proof of Claim 2.
\end{proof}

\begin{remark}
For  $b  = (x_1, \ldots, \overline{x}_1) \in B_{l-1}$ ($l \ge 2$) define
\begin{align*}
\acute{\tau}^l_{l-1}(b) &= (x_1+1, x_2, \ldots, \overline{x}_1) \in B_l,\\
\grave{\tau}^l_{l-1}(b) &= (x_1, \ldots, \ol{x}_2, \overline{x}_1+1) \in B_l.
\end{align*}
Then we have
$\acute{\tau}^l_{l-1}(c_1) \ot \grave{\tau}^k_{k-1}(c_2) \simeq
\acute{\tau}^k_{k-1}(c_2') \ot \grave{\tau}^l_{l-1}(c_1')$
under the isomorphism $B_{l} \ot B_{k} \simeq B_{k} \ot B_{l}$,
if and only if
$c_1 \ot c_2 \simeq c_2' \ot c_1'$
under  $B_{l-1} \ot B_{k-1} \simeq B_{k-1}
\ot B_{l-1}$.
We also have
$H_{B_lB_k}(\acute{\tau}^l_{l-1}(c_1) \ot
\grave{\tau}^k_{k-1}(c_2)) =
H_{B_{l-1}B_{k-1}}(c_1 \ot c_2)$.
\end{remark}
%

\renewcommand{\theequation}{\Alph{section}.\arabic{equation}}
\appendix
\section{Proof of Proposition \ref{th:main2}}\label{ap:proofiso}
\subsection{\mathversion{bold}Column insertions of type I $U_q(C_n)$ highest elements}

Let us consider an element in $B_{l} \ot B_{k}$ depicted by
\par\noindent
\setlength{\unitlength}{5mm}
\begin{picture}(22,3)
\put(0,1){\makebox(4,1){$b_1 \ot b_2=$}}
\put(4,1){\line(1,0){11}}
\put(4,2){\line(1,0){11}}
\put(4,1){\line(0,1){1}}
\put(15,1){\line(0,1){1}}
\put(4,1){\makebox(11,1){$1$}}
\put(4,2){\makebox(11,1){$\scriptstyle{l}$}}
\put(15,1){\makebox(1,1){$\ot$}}
\put(16,1){\line(1,0){5}}
\put(16,2){\line(1,0){5}}
\put(16,1){\line(0,1){1}}
\put(17,1){\line(0,1){1}}
\put(18,1){\line(0,1){1}}
\put(19,1){\line(0,1){1}}
\put(20,1){\line(0,1){1}}
\put(21,1){\line(0,1){1}}
\put(16,1){\makebox(1,1){$0$}}
\put(17,1){\makebox(1,1){$1$}}
\put(18,1){\makebox(1,1){$2$}}
\put(19,1){\makebox(1,1){$\bar{1}$}}
\put(20,1){\makebox(1,1){$\bar{0}$}}
\put(16,0){\makebox(1,1){$\scriptstyle{x_0}$}}
\put(17,0){\makebox(1,1){$\scriptstyle{x_1}$}}
\put(18,0){\makebox(1,1){$\scriptstyle{x_2}$}}
\put(19,0){\makebox(1,1){$\scriptstyle{\overline{x}_1}$}}
\put(20,0){\makebox(1,1){$\scriptstyle{x_0}$}}
\put(21,1){\makebox(1,0){.}}
\end{picture}
\par\noindent
(In this appendix we denote $\T(b)$ simply by $b$.)
It is a $U_q(C_n)$ highest element.
We denote by $b_2' \ot b_1'$ the image of this element under the
isomorphism $\iota: \; B_{l} \ot B_{k} \to B_{k} \ot B_{l}$
.
\subsubsection{}
Let $\ol{x}_1 \leq x_1$.
Then $b_2' \ot b_1'$ is depicted by
\par\noindent
\setlength{\unitlength}{5mm}
\begin{picture}(22,3)
\put(0,1){\makebox(4,1){$b_2' \ot b_1'=$}}
\put(4,1){\line(1,0){5}}
\put(4,2){\line(1,0){5}}
\put(4,1){\line(0,1){1}}
\put(9,1){\line(0,1){1}}
\put(4,1){\makebox(5,1){$1$}}
\put(4,2){\makebox(5,1){$\scriptstyle{k}$}}
\put(9,1){\makebox(1,1){$\ot$}}
\put(10,1){\line(1,0){11}}
\put(10,2){\line(1,0){11}}
\put(10,1){\line(0,1){1}}
\put(12,1){\line(0,1){1}}
\put(15,1){\line(0,1){1}}
\put(17,1){\line(0,1){1}}
\put(19,1){\line(0,1){1}}
\put(21,1){\line(0,1){1}}
\put(10,1){\makebox(2,1){$0$}}
\put(12,1){\makebox(3,1){$1$}}
\put(15,1){\makebox(2,1){$2$}}
\put(17,1){\makebox(2,1){$\bar{1}$}}
\put(19,1){\makebox(2,1){$\bar{0}$}}
\put(10,0){\makebox(2,1){$\scriptstyle{x_0}$}}
\put(12,2){\makebox(3,1){$\scriptstyle{x_1+l-k}$}}
\put(15,0){\makebox(2,1){$\scriptstyle{x_2}$}}
\put(17,0){\makebox(2,1){$\scriptstyle{\overline{x}_1}$}}
\put(19,0){\makebox(2,1){$\scriptstyle{x_0}$}}
\end{picture}
\par\noindent
The column insertions $( b_2 \longrightarrow b_1)$ and
$( b_1' \longrightarrow b_2')$ lead to the same intermediate result;
\par\noindent
\setlength{\unitlength}{5mm}
\begin{picture}(20,4)
\put(1,1){\line(1,0){4}}
\put(1,2){\line(1,0){4}}
\put(1,1){\line(0,1){1}}
\put(3,1){\line(0,1){1}}
\put(5,1){\line(0,1){1}}
\put(1,1){\makebox(2,1){$0$}}
\put(3,1){\makebox(2,1){$1$}}
\put(1,0){\makebox(2,1){$\scriptstyle{x_0}$}}
\put(3,0){\makebox(2,1){$\scriptstyle{x_1-\overline{x}_1}$}}
\put(5,1){\makebox(3,1){$\longrightarrow $}}
\put(8,1){\line(1,0){8}}
\put(8,2){\line(1,0){11}}
\put(8,3){\line(1,0){11}}
\put(8,1){\line(0,1){2}}
\put(11,1){\line(0,1){2}}
\put(13,1){\line(0,1){1}}
\put(16,1){\line(0,1){1}}
\put(19,2){\line(0,1){1}}
\put(8,1){\makebox(3,1){$1$}}
\put(11,1){\makebox(2,1){$2$}}
\put(13,1){\makebox(3,1){$\overline{0}$}}
\put(8,2){\makebox(3,1){$0$}}
\put(11,2){\makebox(8,1){$1$}}
\put(8,0){\makebox(3,1){$\scriptstyle{\overline{x}_1}$}}
\put(11,0){\makebox(2,1){$\scriptstyle{x_2}$}}
\put(13,0){\makebox(3,1){$\scriptstyle{\overline{x}_1+x_0}$}}
\put(11,3){\makebox(8,1){$\scriptstyle{l-\overline{x}_1}$}}
\end{picture}
\par\noindent
The value of the energy function is $x_0+x_1-\ol{x}_1$.
\subsubsection{}
Let $\ol{x}_1 > x_1$.
Then $b_2' \ot b_1'$ is depicted by
\par\noindent
\setlength{\unitlength}{5mm}
\begin{picture}(22,3)
\put(0,1){\makebox(4,1){$b_2' \ot b_1'=$}}
\put(4,1){\line(1,0){5}}
\put(4,2){\line(1,0){5}}
\put(4,1){\line(0,1){1}}
\put(9,1){\line(0,1){1}}
\put(4,1){\makebox(5,1){$1$}}
\put(4,2){\makebox(5,1){$\scriptstyle{k}$}}
\put(9,1){\makebox(1,1){$\ot$}}
\put(10,1){\line(1,0){11}}
\put(10,2){\line(1,0){11}}
\put(10,1){\line(0,1){1}}
\put(12,1){\line(0,1){1}}
\put(15,1){\line(0,1){1}}
\put(17,1){\line(0,1){1}}
\put(19,1){\line(0,1){1}}
\put(21,1){\line(0,1){1}}
\put(10,1){\makebox(2,1){$0$}}
\put(12,1){\makebox(3,1){$1$}}
\put(15,1){\makebox(2,1){$2$}}
\put(17,1){\makebox(2,1){$\bar{1}$}}
\put(19,1){\makebox(2,1){$\bar{0}$}}
\put(10,0){\makebox(2,1){$\scriptstyle{x_0+y}$}}
\put(12,2){\makebox(3,1){$\scriptstyle{x_1+l-k-y}$}}
\put(15,0){\makebox(2,1){$\scriptstyle{x_2}$}}
\put(17,0){\makebox(2,1){$\scriptstyle{\overline{x}_1-y}$}}
\put(19,0){\makebox(2,1){$\scriptstyle{x_0+y}$}}
\end{picture}
\par\noindent
where
\begin{equation}
y=\min [l-k,\overline{x}_1-x_1].
\end{equation}
The column insertions $( b_2 \longrightarrow b_1)$ and
$( b_1' \longrightarrow b_2')$ lead to the same intermediate result;
\par\noindent
For $x_1 + x_2 > \ol{x}_1$,
\par\noindent
\setlength{\unitlength}{5mm}
\begin{picture}(20,4)
\put(1,1){\line(1,0){4}}
\put(1,2){\line(1,0){4}}
\put(1,1){\line(0,1){1}}
\put(5,1){\line(0,1){1}}
\put(1,1){\makebox(4,1){$0$}}
\put(1,0){\makebox(4,1){$\scriptstyle{x_0}$}}
\put(5,1){\makebox(3,1){$\longrightarrow $}}
\put(8,1){\line(1,0){8}}
\put(8,2){\line(1,0){11}}
\put(8,3){\line(1,0){11}}
\put(8,1){\line(0,1){2}}
\put(10,1){\line(0,1){1}}
\put(11,2){\line(0,1){1}}
\put(12.5,1){\line(0,1){1}}
\put(16,1){\line(0,1){1}}
\put(19,2){\line(0,1){1}}
\put(8,1){\makebox(2,1){$1$}}
\put(10,1){\makebox(2.5,1){$2$}}
\put(12.5,1){\makebox(3.5,1){$\overline{0}$}}
\put(8,2){\makebox(3,1){$0$}}
\put(11,2){\makebox(8,1){$1$}}
\put(8,0){\makebox(2,1){$\scriptstyle{x_1}$}}
\put(10,0){\makebox(2.5,1){$\scriptstyle{x_2}$}}
\put(12.5,0){\makebox(3.5,1){$\scriptstyle{\overline{x}_1+x_0}$}}
\put(8,3){\makebox(3,1){$\scriptstyle{\overline{x}_1}$}}
\put(11,3){\makebox(8,1){$\scriptstyle{l-\overline{x}_1}$}}
\end{picture}
\par\noindent
For $x_1 + x_2 \leq \ol{x}_1$,
\par\noindent
\setlength{\unitlength}{5mm}
\begin{picture}(20,4)
\put(1,1){\line(1,0){4}}
\put(1,2){\line(1,0){4}}
\put(1,1){\line(0,1){1}}
\put(5,1){\line(0,1){1}}
\put(1,1){\makebox(4,1){$0$}}
\put(1,0){\makebox(4,1){$\scriptstyle{x_0}$}}
\put(5,1){\makebox(3,1){$\longrightarrow $}}
\put(8,1){\line(1,0){8}}
\put(8,2){\line(1,0){11}}
\put(8,3){\line(1,0){11}}
\put(8,1){\line(0,1){2}}
\put(9.5,1){\line(0,1){1}}
\put(11,1){\line(0,1){2}}
\put(11.9,1){\line(0,1){2}}
\put(12.1,1){\line(0,1){2}}
\put(13,1){\line(0,1){2}}
\put(16,1){\line(0,1){1}}
\put(19,2){\line(0,1){1}}
\put(8,1){\makebox(1.5,1){$1$}}
\put(9.5,1){\makebox(1.5,1){$2$}}
\put(11,1){\makebox(0.9,1){$\overline{1}$}}
\put(12.1,1){\makebox(0.9,1){$\overline{0}$}}
\put(11,2){\makebox(0.9,1){$0$}}
\put(12.1,2){\makebox(0.9,1){$1$}}
\put(13,1){\makebox(3,1){$\overline{0}$}}
\put(8,2){\makebox(3,1){$0$}}
\put(13,2){\makebox(6,1){$1$}}
\put(8,0){\makebox(1.5,1){$\scriptstyle{x_1}$}}
\put(9.5,0){\makebox(1.5,1){$\scriptstyle{x_2}$}}
\put(11,3){\makebox(2,1){$\scriptstyle{\overline{x}_1-x_1-x_2}$}}
\put(13,0){\makebox(3,1){$\scriptstyle{x_0+x_1+x_2}$}}
\put(13,3){\makebox(6,1){$\scriptstyle{l-\overline{x}_1}$}}
\end{picture}
\par\noindent
The value of the energy function is $x_0$.
Here and in the following we use the notation:
\par\noindent
\setlength{\unitlength}{5mm}
\begin{picture}(20,4)
\put(0,1){\line(1,0){2}}
\put(0,2){\line(1,0){2}}
\put(0,3){\line(1,0){2}}
\put(0,1){\line(0,1){2}}
\put(0.9,1){\line(0,1){2}}
\put(1.1,1){\line(0,1){2}}
\put(2,1){\line(0,1){2}}
\put(0,1){\makebox(0.9,1){$\overline{1}$}}
\put(1.1,1){\makebox(0.9,1){$\overline{0}$}}
\put(0,2){\makebox(0.9,1){$0$}}
\put(1.1,2){\makebox(0.9,1){$1$}}
\put(0,3){\makebox(2,1){$\scriptstyle{m}$}}
\put(2,1.5){\makebox(4,1){$=$}}
\put(6,1){\line(1,0){2}}
\put(6,2){\line(1,0){2}}
\put(6,3){\line(1,0){2}}
\put(6,1){\line(0,1){2}}
\put(7,1){\line(0,1){2}}
\put(8,1){\line(0,1){2}}
\put(6,1){\makebox(1,1){$\overline{1}$}}
\put(7,1){\makebox(1,1){$\overline{0}$}}
\put(6,2){\makebox(1,1){$0$}}
\put(7,2){\makebox(1,1){$1$}}
\put(6,3){\makebox(1,1){$\scriptstyle{\frac{m}{2}}$}}
\put(7,3){\makebox(1,1){$\scriptstyle{\frac{m}{2}}$}}
\put(8,1.5){\makebox(5,1){($m$: even),}}
\put(13,1){\line(1,0){2}}
\put(13,2){\line(1,0){2}}
\put(13,3){\line(1,0){2}}
\put(13,1){\line(0,1){2}}
\put(13.9,2){\line(0,1){1}}
\put(14.1,1){\line(0,1){1}}
\put(15,1){\line(0,1){2}}
\put(13,1){\makebox(1.1,1){$\overline{1}$}}
\put(14.1,1){\makebox(0.9,1){$\overline{0}$}}
\put(13,2){\makebox(0.9,1){$0$}}
\put(13.9,2){\makebox(1.1,1){$1$}}
\put(12.9,3){\makebox(0.9,1){$\scriptstyle{\frac{m-1}{2}}$}}
\put(14.1,3){\makebox(1.1,1){$\scriptstyle{\frac{m+1}{2}}$}}
\put(12.9,0){\makebox(0.9,1){$\scriptstyle{\frac{m+1}{2}}$}}
\put(14.1,0){\makebox(1.1,1){$\scriptstyle{\frac{m-1}{2}}$}}
\put(15,1.5){\makebox(5,1){($m$: odd).}}
\end{picture}

\subsection{\mathversion{bold}Column insertions of type II $U_q(C_n)$ highest elements}

Let
\par\noindent
\setlength{\unitlength}{5mm}
\begin{picture}(22,3)
\put(0,1){\makebox(4,1){$b_1 \ot b_2=$}}
\put(4,1){\line(1,0){11}}
\put(4,2){\line(1,0){11}}
\put(4,1){\line(0,1){1}}
\put(6,1){\line(0,1){1}}
\put(13,1){\line(0,1){1}}
\put(15,1){\line(0,1){1}}
\put(4,1){\makebox(2,1){$0$}}
\put(6,1){\makebox(7,1){$1$}}
\put(13,1){\makebox(2,1){$\bar{0}$}}
\put(4,0){\makebox(2,1){$\scriptstyle{y_0}$}}
\put(6,2){\makebox(7,1){$\scriptstyle{l-2y_0}$}}
\put(13,0){\makebox(2,1){$\scriptstyle{y_0}$}}
\put(15,1){\makebox(1,1){$\ot$}}
\put(16,1){\line(1,0){5}}
\put(16,2){\line(1,0){5}}
\put(16,1){\line(0,1){1}}
\put(17.5,1){\line(0,1){1}}
\put(19.5,1){\line(0,1){1}}
\put(21,1){\line(0,1){1}}
\put(16,1){\makebox(1.5,1){$1$}}
\put(17.5,1){\makebox(2,1){$2$}}
\put(19.5,1){\makebox(1.5,1){$\bar{1}$}}
\put(16,0){\makebox(1.5,1){$\scriptstyle{x_1}$}}
\put(17.5,0){\makebox(2,1){$\scriptstyle{x_2}$}}
\put(19.5,0){\makebox(1.5,1){$\scriptstyle{\overline{x}_1}$}}
\end{picture}
\par\noindent
be a $U_q(C_n)$ highest element in $B_{l} \ot B_{k}$.
Thus we assume $l-2y_0 \geq x_2+\ol{x}_1$.
\subsubsection{}
Let $l-k > y_0 \geq x_1 - \ol{x}_1$.
Then $b_2' \ot b_1'$ is depicted by
\par\noindent
\setlength{\unitlength}{5mm}
\begin{picture}(22,3)
\put(0,1){\makebox(4,1){$b_2' \ot b_1'=$}}
\put(4,1){\line(1,0){5}}
\put(4,2){\line(1,0){5}}
\put(4,1){\line(0,1){1}}
\put(9,1){\line(0,1){1}}
\put(4,1){\makebox(5,1){$1$}}
\put(4,2){\makebox(5,1){$\scriptstyle{k}$}}
\put(9,1){\makebox(1,1){$\ot$}}
\put(10,1){\line(1,0){11}}
\put(10,2){\line(1,0){11}}
\put(10,1){\line(0,1){1}}
\put(12,1){\line(0,1){1}}
\put(15,1){\line(0,1){1}}
\put(17,1){\line(0,1){1}}
\put(19,1){\line(0,1){1}}
\put(21,1){\line(0,1){1}}
\put(10,1){\makebox(2,1){$0$}}
\put(12,1){\makebox(3,1){$1$}}
\put(15,1){\makebox(2,1){$2$}}
\put(17,1){\makebox(2,1){$\bar{1}$}}
\put(19,1){\makebox(2,1){$\bar{0}$}}
\put(10,0){\makebox(2,1){$\scriptstyle{z}$}}
\put(12,2){\makebox(3,1){$\scriptstyle{x_1+l-k-y_0-z}$}}
\put(15,0){\makebox(2,1){$\scriptstyle{x_2}$}}
\put(17,2){\makebox(2,1){$\scriptstyle{\overline{x}_1+y_0-z}$}}
\put(19,0){\makebox(2,1){$\scriptstyle{z}$}}
\end{picture}
\par\noindent
where
\begin{equation}
z=\min [y_0+\overline{x}_1-x_1,l-k-y_0].
\end{equation}
The column insertions $( b_2 \longrightarrow b_1)$ and
$( b_1' \longrightarrow b_2')$ give the same result;
\par\noindent
For $y_0 \geq k$,
\par\noindent
\setlength{\unitlength}{5mm}
\begin{picture}(20,4)
\put(8,1){\line(1,0){4.5}}
\put(8,2){\line(1,0){11}}
\put(8,3){\line(1,0){11}}
\put(8,1){\line(0,1){2}}
\put(9.5,1){\line(0,1){1}}
\put(11,1){\line(0,1){1}}
\put(12.5,1){\line(0,1){1}}
\put(13,2){\line(0,1){1}}
\put(14,2){\line(0,1){1}}
\put(19,2){\line(0,1){1}}
\put(8,1){\makebox(1.5,1){$1$}}
\put(9.5,1){\makebox(1.5,1){$2$}}
\put(11,1){\makebox(1.5,1){$\overline{1}$}}
\put(8,2){\makebox(5,1){$0$}}
\put(13,2){\makebox(1,1){$1$}}
\put(14,2){\makebox(5,1){$\overline{0}$}}
\put(8,0){\makebox(1.5,1){$\scriptstyle{x_1}$}}
\put(9.5,0){\makebox(1.5,1){$\scriptstyle{x_2}$}}
\put(11,0){\makebox(1.5,1){$\scriptstyle{\overline{x}_1}$}}
\put(8,3){\makebox(5,1){$\scriptstyle{y_0}$}}
\put(13,3){\makebox(1,1){$\scriptstyle{l-2y_0}$}}
\put(14,3){\makebox(5,1){$\scriptstyle{y_0}$}}
\end{picture}
\par\noindent
For $k > y_0 \geq x_1 + x_2$,
\par\noindent
\setlength{\unitlength}{5mm}
\begin{picture}(20,4)
\put(8,1){\line(1,0){5.5}}
\put(8,2){\line(1,0){11}}
\put(8,3){\line(1,0){11}}
\put(8,1){\line(0,1){2}}
\put(9,1){\line(0,1){1}}
\put(10,1){\line(0,1){1}}
\put(11.5,1){\line(0,1){2}}
\put(12.4,1){\line(0,1){2}}
\put(12.6,1){\line(0,1){2}}
\put(13.5,1){\line(0,1){2}}
\put(15.5,2){\line(0,1){1}}
\put(19,2){\line(0,1){1}}
\put(8,1){\makebox(1,1){$1$}}
\put(9,1){\makebox(1,1){$2$}}
\put(10,1){\makebox(1.5,1){$\overline{1}$}}
\put(8,2){\makebox(3.5,1){$0$}}
\put(11.5,1){\makebox(0.9,1){$\overline{1}$}}
\put(12.6,1){\makebox(0.9,1){$\overline{0}$}}
\put(11.5,2){\makebox(0.9,1){$0$}}
\put(12.6,2){\makebox(0.9,1){$1$}}
\put(13.5,2){\makebox(2,1){$1$}}
\put(15.5,2){\makebox(3.5,1){$\overline{0}$}}
\put(8,0){\makebox(1,1){$\scriptstyle{x_1}$}}
\put(9,0){\makebox(1,1){$\scriptstyle{x_2}$}}
\put(8,3){\makebox(3.5,1){$\scriptstyle{y_0}$}}
\put(11.5,0){\makebox(2,1){$\scriptstyle{k-y_0}$}}
\put(13.5,3){\makebox(2,1){$\scriptstyle{l-k-y_0}$}}
\put(15.5,3){\makebox(3.5,1){$\scriptstyle{y_0}$}}
\end{picture}
\par\noindent
For $x_1+x_2 > y_0 \geq x_1 + x_2-\ol{x}_1$,
\par\noindent
\setlength{\unitlength}{5mm}
\begin{picture}(20,4.5)
\put(8,1){\line(1,0){8.5}}
\put(8,2){\line(1,0){11}}
\put(8,3){\line(1,0){11}}
\put(8,1){\line(0,1){2}}
\put(10,1){\line(0,1){1}}
\put(11.5,1){\line(0,1){2}}
\put(12.4,1){\line(0,1){2}}
\put(12.6,1){\line(0,1){2}}
\put(13.5,1){\line(0,1){2}}
\put(16.5,1){\line(0,1){1}}
\put(17,2){\line(0,1){1}}
\put(19,2){\line(0,1){1}}
\put(8,1){\makebox(2,1){$1$}}
\put(10,1){\makebox(1.5,1){$2$}}
\put(8,2){\makebox(3.5,1){$0$}}
\put(11.5,1){\makebox(0.9,1){$\overline{1}$}}
\put(12.6,1){\makebox(0.9,1){$\overline{0}$}}
\put(11.5,2){\makebox(0.9,1){$0$}}
\put(12.6,2){\makebox(0.9,1){$1$}}
\put(13.5,1){\makebox(3,1){$\overline{0}$}}
\put(13.5,2){\makebox(3.5,1){$1$}}
\put(17,2){\makebox(2,1){$\overline{0}$}}
\put(8,0){\makebox(2,1){$\scriptstyle{x_1}$}}
\put(10,0){\makebox(1.5,1){$\scriptstyle{x_2}$}}
\put(11.5,3.1){\makebox(2,1){$\stackrel{\scriptstyle \ol{x}_1+y_0-x_1}
{\scriptstyle -x_2}
$}}
\put(13.5,0){\makebox(3,1){$\scriptstyle{x_1+x_2-y_0}$}}
\put(14,3){\makebox(3,1){$\scriptstyle{l-2y_0-\ol{x}_1}$}}
\put(17,3){\makebox(2,1){$\scriptstyle{y_0}$}}
\end{picture}
\par\noindent
For $x_1+x_2-\ol{x}_1 > y_0$,
\par\noindent
\setlength{\unitlength}{5mm}
\begin{picture}(20,4)
\put(8,1){\line(1,0){8.5}}
\put(8,2){\line(1,0){11}}
\put(8,3){\line(1,0){11}}
\put(8,1){\line(0,1){2}}
\put(10,1){\line(0,1){1}}
\put(12.5,2){\line(0,1){1}}
\put(13.5,1){\line(0,1){1}}
\put(16.5,1){\line(0,1){1}}
\put(17,2){\line(0,1){1}}
\put(19,2){\line(0,1){1}}
\put(8,1){\makebox(2,1){$1$}}
\put(10,1){\makebox(3.5,1){$2$}}
\put(8,2){\makebox(4.5,1){$0$}}
\put(13.5,1){\makebox(3,1){$\overline{0}$}}
\put(12.5,2){\makebox(4.5,1){$1$}}
\put(17,2){\makebox(2,1){$\overline{0}$}}
\put(8,0){\makebox(2,1){$\scriptstyle{x_1}$}}
\put(10,0){\makebox(3.5,1){$\scriptstyle{x_2}$}}
\put(13.5,0){\makebox(3,1){$\scriptstyle{\ol{x}_1}$}}
\put(8,3){\makebox(4.5,1){$\scriptstyle{y_0+\ol{x}_1}$}}
\put(12.5,3){\makebox(4.5,1){$\scriptstyle{l-2y_0-\ol{x}_1}$}}
\put(17,3){\makebox(2,1){$\scriptstyle{y_0}$}}
\end{picture}
\par\noindent
The value of the energy function is $0$.
\subsubsection{}
Let $l-k \leq y_0$ and $2y_0+k-l- x_1 + \ol{x}_1 > 0$.
Then $b_2' \ot b_1'$ is depicted by
\par\noindent
\setlength{\unitlength}{5mm}
\begin{picture}(22,3)
\put(0,1){\makebox(4,1){$b_2' \ot b_1'=$}}
\put(4,1){\line(1,0){5}}
\put(4,2){\line(1,0){5}}
\put(4,1){\line(0,1){1}}
\put(5.5,1){\line(0,1){1}}
\put(7.5,1){\line(0,1){1}}
\put(9,1){\line(0,1){1}}
\put(4,1){\makebox(1.5,1){$0$}}
\put(5.5,1){\makebox(2,1){$1$}}
\put(7.5,1){\makebox(1.5,1){$\ol{0}$}}
\put(4,0){\makebox(1.5,1){$\scriptstyle{y_0-l+k}$}}
\put(5.5,2){\makebox(2,1){$\scriptstyle{2l-k-2y_0}$}}
\put(7.5,0){\makebox(1.5,1){$\scriptstyle{y_0-l+k}$}}
\put(9,1){\makebox(1,1){$\ot$}}
\put(10,1){\line(1,0){11}}
\put(10,2){\line(1,0){11}}
\put(10,1){\line(0,1){1}}
\put(13,1){\line(0,1){1}}
\put(16,1){\line(0,1){1}}
\put(21,1){\line(0,1){1}}
\put(10,1){\makebox(3,1){$1$}}
\put(13,1){\makebox(3,1){$2$}}
\put(16,1){\makebox(5,1){$\bar{1}$}}
\put(10,0){\makebox(3,1){$\scriptstyle{x_1}$}}
\put(13,0){\makebox(3,1){$\scriptstyle{x_2}$}}
\put(16,0){\makebox(5,1){$\scriptstyle{\overline{x}_1+l-k}$}}
\end{picture}
\par\noindent
The column insertions $( b_2 \longrightarrow b_1)$ and
$( b_1' \longrightarrow b_2')$ lead to the same intermediate result;
\par\noindent
For $l-2y_0 \geq 2 \ol{x}_1$,
\par\noindent
\setlength{\unitlength}{5mm}
\begin{picture}(20,4)
\put(1,1){\line(1,0){4}}
\put(1,2){\line(1,0){4}}
\put(1,1){\line(0,1){1}}
\put(5,1){\line(0,1){1}}
\put(1,1){\makebox(4,1){$1$}}
\put(1,0){\makebox(4,1){$\scriptstyle{y_0-l+k}$}}
\put(5,1){\makebox(3,1){$\longrightarrow $}}
\put(8,1){\line(1,0){8}}
\put(8,2){\line(1,0){11}}
\put(8,3){\line(1,0){11}}
\put(8,1){\line(0,1){2}}
\put(10.5,1){\line(0,1){1}}
\put(11,2){\line(0,1){1}}
\put(13,1){\line(0,1){1}}
\put(16,1){\line(0,1){2}}
\put(19,2){\line(0,1){1}}
\put(8,1){\makebox(2.5,1){$1$}}
\put(10.5,1){\makebox(2.5,1){$2$}}
\put(13,1){\makebox(3.5,1){$\overline{0}$}}
\put(8,2){\makebox(3,1){$0$}}
\put(11,2){\makebox(5,1){$1$}}
\put(16,2){\makebox(3,1){$\ol{0}$}}
\put(8,0){\makebox(2.5,1){$\scriptstyle{l-y_0-\ol{x}_1-x_2}$}}
\put(10.5,0){\makebox(2.5,1){$\scriptstyle{x_2}$}}
\put(13,0){\makebox(3.5,1){$\scriptstyle{\overline{x}_1}$}}
\put(8,3){\makebox(3,1){$\scriptstyle{y_0+\overline{x}_1}$}}
\put(11,3){\makebox(5,1){$\scriptstyle{l-2y_0-\overline{x}_1}$}}
\put(16,3){\makebox(3,1){$\scriptstyle{y_0}$}}
\end{picture}
\par\noindent
For $l-2y_0 < 2\ol{x}_1$,
\par\noindent
\setlength{\unitlength}{5mm}
\begin{picture}(20,4)
\put(1,1){\line(1,0){4}}
\put(1,2){\line(1,0){4}}
\put(1,1){\line(0,1){1}}
\put(5,1){\line(0,1){1}}
\put(1,1){\makebox(4,1){$1$}}
\put(1,0){\makebox(4,1){$\scriptstyle{y_0-l+k}$}}
\put(5,1){\makebox(3,1){$\longrightarrow $}}
\put(8,1){\line(1,0){8}}
\put(8,2){\line(1,0){11}}
\put(8,3){\line(1,0){11}}
\put(8,1){\line(0,1){2}}
\put(10,1){\line(0,1){1}}
\put(11,1){\line(0,1){2}}
\put(11.9,1){\line(0,1){2}}
\put(12.1,1){\line(0,1){2}}
\put(13,1){\line(0,1){2}}
\put(16,1){\line(0,1){2}}
\put(19,2){\line(0,1){1}}
\put(8,1){\makebox(2,1){$1$}}
\put(10,1){\makebox(1,1){$2$}}
\put(11,1){\makebox(0.9,1){$\overline{1}$}}
\put(12.1,1){\makebox(0.9,1){$\overline{0}$}}
\put(11,2){\makebox(0.9,1){$0$}}
\put(12.1,2){\makebox(0.9,1){$1$}}
\put(13,1){\makebox(3,1){$\overline{0}$}}
\put(8,2){\makebox(3,1){$0$}}
\put(13,2){\makebox(3,1){$1$}}
\put(16,2){\makebox(3,1){$\ol{0}$}}
\put(7.8,0){\makebox(2,1){$\stackrel{\scriptstyle l-y_0-\ol{x}_1}
{\scriptstyle -x_2}$}}
\put(10,0){\makebox(1,1){$\scriptstyle{x_2}$}}
\put(11,3){\makebox(2,1){$\scriptstyle{2 \overline{x}_1-l+2y_0}$}}
\put(13,0){\makebox(3,1){$\scriptstyle{l-2y_0-\ol{x}_1}$}}
\put(16,3){\makebox(3,1){$\scriptstyle{y_0}$}}
\end{picture}
\par\noindent
The value of the energy function is $y_0-l+k$.
\subsubsection{}
Let $y_0 < x_1-\ol{x}_1$ and $2y_0+k-l- x_1 + \ol{x}_1 \leq 0$.
Then $b_2' \ot b_1'$ is depicted by
\par\noindent
\setlength{\unitlength}{5mm}
\begin{picture}(22,3)
\put(0,1){\makebox(4,1){$b_2' \ot b_1'=$}}
\put(4,1){\line(1,0){5}}
\put(4,2){\line(1,0){5}}
\put(4,1){\line(0,1){1}}
\put(5.5,1){\line(0,1){1}}
\put(7.5,1){\line(0,1){1}}
\put(9,1){\line(0,1){1}}
\put(4,1){\makebox(1.5,1){$0$}}
\put(5.5,1){\makebox(2,1){$1$}}
\put(7.5,1){\makebox(1.5,1){$\ol{0}$}}
\put(4,0){\makebox(1.5,1){$\scriptstyle{y_0-w}$}}
\put(5.5,2){\makebox(2,1){$\scriptstyle{k-2y_0+2w}$}}
\put(7.5,0){\makebox(1.5,1){$\scriptstyle{y_0-w}$}}
\put(9,1){\makebox(1,1){$\ot$}}
\put(10,1){\line(1,0){11}}
\put(10,2){\line(1,0){11}}
\put(10,1){\line(0,1){1}}
\put(13,1){\line(0,1){1}}
\put(16,1){\line(0,1){1}}
\put(21,1){\line(0,1){1}}
\put(10,1){\makebox(3,1){$1$}}
\put(13,1){\makebox(3,1){$2$}}
\put(16,1){\makebox(5,1){$\bar{1}$}}
\put(10,2){\makebox(3,1){$\scriptstyle{x_1+l-k-w}$}}
\put(13,0){\makebox(3,1){$\scriptstyle{x_2}$}}
\put(16,0){\makebox(5,1){$\scriptstyle{\overline{x}_1+w}$}}
\end{picture}
\par\noindent
where $w=(2y_0-x_1+\ol{x}_1)_+$.
The column insertions $( b_2 \longrightarrow b_1)$ and
$( b_1' \longrightarrow b_2')$ lead to the same intermediate result;
\par\noindent
\setlength{\unitlength}{5mm}
\begin{picture}(20,4)
\put(1,1){\line(1,0){4}}
\put(1,2){\line(1,0){4}}
\put(1,1){\line(0,1){1}}
\put(5,1){\line(0,1){1}}
\put(1,1){\makebox(4,1){$1$}}
\put(1,0){\makebox(4,1){$\scriptstyle{x_1-\ol{x}_1-y_0}$}}
\put(5,1){\makebox(3,1){$\longrightarrow $}}
\put(8,1){\line(1,0){7}}
\put(8,2){\line(1,0){11}}
\put(8,3){\line(1,0){11}}
\put(8,1){\line(0,1){2}}
\put(10.5,1){\line(0,1){2}}
\put(13,1){\line(0,1){1}}
\put(15,1){\line(0,1){1}}
\put(16,2){\line(0,1){1}}
\put(19,2){\line(0,1){1}}
\put(8,1){\makebox(2.5,1){$1$}}
\put(10.5,1){\makebox(2.5,1){$2$}}
\put(13,1){\makebox(2,1){$\overline{0}$}}
\put(8,2){\makebox(2.5,1){$0$}}
\put(10.5,2){\makebox(5.5,1){$1$}}
\put(16,2){\makebox(3,1){$\ol{0}$}}
\put(8,0){\makebox(2.5,1){$\scriptstyle{y_0+\ol{x}_1}$}}
\put(10.5,0){\makebox(2.5,1){$\scriptstyle{x_2}$}}
\put(13,0){\makebox(2,1){$\scriptstyle{\overline{x}_1}$}}
\put(10.5,3){\makebox(5.5,1){$\scriptstyle{l-2y_0-\overline{x}_1}$}}
\put(16,3){\makebox(3,1){$\scriptstyle{y_0}$}}
\end{picture}
\par\noindent
The value of the energy function is $x_1-\ol{x}_1-y_0$.

\section{\mathversion{bold}Alternative rule for $C^{(1)}_n$}
\label{sec:anoth}
\subsection{Algorithm for the isomorphism}
Let $b_1 = (x_1, \ldots, \overline{x}_1) \in B_{l},
b_2 = (y_1, \ldots, \overline{y}_1) \in B_{k}$.
We are going to show the rule of finding the image of $b_1 \otimes b_2$
under the isomorphism
\begin{eqnarray*}
\iota : B_{l} \otimes B_{k} &\stackrel{\sim}{\rightarrow}& B_{k} \otimes B_{l} \\
          b_1 \otimes b_2 &\mapsto& b'_2 \otimes b'_1.
\end{eqnarray*}
Let $x_0 = \overline{x}_0 = (l-\sum_{i=1}^n (x_i + \overline{x}_i))/2$
and $y_0 = \overline{y}_0 = (k-\sum_{i=1}^n (y_i + \overline{y}_i))/2$.
We assume $l \geq k$.
We start with the following initial diagram.

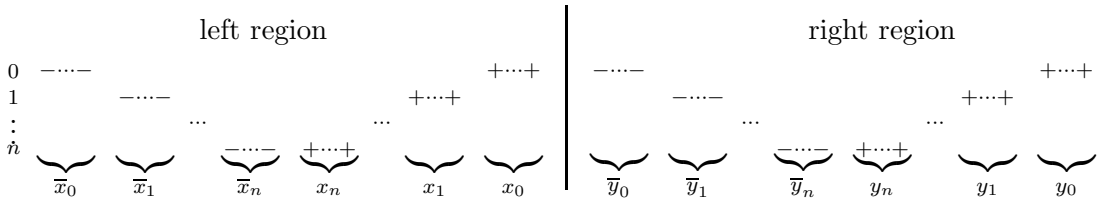
\setlength{\unitlength}{0.7mm}
\begin{picture}(220,40)(20,0)
\put(0,5){\makebox(5,5){$\scriptstyle{n}$}}
\put(0,10){\makebox(5,5){$\scriptstyle{\vdots}$}}
\put(0,15){\makebox(5,5){$\scriptstyle{1}$}}
\put(0,20){\makebox(5,5){$\scriptstyle{0}$}}
\put(5,20){\makebox(15,5){$\scriptstyle{- \cdots -}$}}
\put(20,15){\makebox(15,5){$\scriptstyle{- \cdots -}$}}
\put(35,10){\makebox(5,5){$\scriptstyle{\cdots}$}}
\put(40,5){\makebox(15,5){$\scriptstyle{- \cdots -}$}}
\put(90,20){\makebox(15,5){$\scriptstyle{+ \cdots +}$}}
\put(75,15){\makebox(15,5){$\scriptstyle{+ \cdots +}$}}
\put(70,10){\makebox(5,5){$\scriptstyle{\cdots}$}}
\put(55,5){\makebox(15,5){$\scriptstyle{+ \cdots +}$}}
\put(5,0){\makebox(15,5){$\underbrace{\hphantom{\scriptstyle{- \cdots -}}}_{\scriptstyle{\overline{x}_0}}$}}
\put(20,0){\makebox(15,5){$\underbrace{\hphantom{\scriptstyle{- \cdots -}}}_{\scriptstyle{\overline{x}_1}}$}}
\put(40,0){\makebox(15,5){$\underbrace{\hphantom{\scriptstyle{- \cdots -}}}_{\scriptstyle{\overline{x}_n}}$}}
\put(55,0){\makebox(15,5){$\underbrace{\hphantom{\scriptstyle{- \cdots -}}}_{\scriptstyle{x_n}}$}}
\put(75,0){\makebox(15,5){$\underbrace{\hphantom{\scriptstyle{- \cdots -}}}_{\scriptstyle{x_1}}$}}
\put(90,0){\makebox(15,5){$\underbrace{\hphantom{\scriptstyle{- \cdots -}}}_{\scriptstyle{x_0}}$}}
\put(107.5,0){\line(0,1){35}}
\put(50,30){\makebox(0,0){\small left region}}
\put(167.5,30){\makebox(0,0){\small right region}}
\put(110,20){\makebox(15,5){$\scriptstyle{- \cdots -}$}}
\put(125,15){\makebox(15,5){$\scriptstyle{- \cdots -}$}}
\put(140,10){\makebox(5,5){$\scriptstyle{\cdots}$}}
\put(145,5){\makebox(15,5){$\scriptstyle{- \cdots -}$}}
\put(195,20){\makebox(15,5){$\scriptstyle{+ \cdots +}$}}
\put(180,15){\makebox(15,5){$\scriptstyle{+ \cdots +}$}}
\put(175,10){\makebox(5,5){$\scriptstyle{\cdots}$}}
\put(160,5){\makebox(15,5){$\scriptstyle{+ \cdots +}$}}
\put(110,0){\makebox(15,5){$\underbrace{\hphantom{\scriptstyle{- \cdots -}}}_{\overline{y}_0}$}}
\put(125,0){\makebox(15,5){$\underbrace{\hphantom{\scriptstyle{- \cdots -}}}_{\overline{y}_1}$}}
\put(145,0){\makebox(15,5){$\underbrace{\hphantom{\scriptstyle{- \cdots -}}}_{\overline{y}_n}$}}
\put(160,0){\makebox(15,5){$\underbrace{\hphantom{\scriptstyle{- \cdots -}}}_{y_n}$}}
\put(180,0){\makebox(15,5){$\underbrace{\hphantom{\scriptstyle{- \cdots -}}}_{y_1}$}}
\put(195,0){\makebox(15,5){$\underbrace{\hphantom{\scriptstyle{- \cdots -}}}_{y_0}$}}

\end{picture}\\[4ex]
\noindent
By using Lemma \ref{lem:1} one can remove 
$(\, \fbx{0}\, , \fbx{\ol{0}}\, )$ pairs 
from $b_1$ and $b_2$ simultaneously as many times as possible.
Thus we assume in the following that either $x_0$ or $y_0$ is equal to $0$
throughout.
Then the general procedure to obtain the isomorphism and energy function is as
follows.
\begin{enumerate}
\setcounter{enumi}{-1}
    \item Each symbol $+$ or $-$ is marked or unmarked.
	In the initial diagram all the symbols are unmarked.
	\item There are three regions (left, right, and middle---the
	latter is empty in the initial diagram).
	Pick the leftmost symbol $a$ in the right region.
	Find $a$'s partner $b$ in the left region according to the rule 2-3.
	Apply (a), and repeat this procedure as many turns as possible, and
	then apply (b).
	During the procedure if a symbol named $a$ is a $+$ (resp. $-$)
	symbol we call it $+_a$ (resp. $-_a$).
	\begin{enumerate}
		\item If $a$ exists and
		there is the partner $b$, mark $b$
		according to the rule 4.
	Put a new line on the right of $a$
	which forms the new boundary between the middle region and the right region.
	(In the second turn or later, delete the old line on the left of $a$.)
	    \item If $a$ does not exist or
	    there is no partner of $a$, then stop.
	Enumerate the cardinality of the symbols in the right region
	and denote it by $h$.
	This $h$ is equal to the value of the energy function, which is
	so normalized as the minimal value is equal to $0$.
	Proceed to (c) or (d) according to the value of $h$.
	    \item If $h=0$, the procedure is finished. See (e).
	    \item If $h>0$, give up the diagram.
	    Go back to the initial diagram and
	    mark the leftmost $h$ symbols in the left region.
	    Then start again the procedure from the rule 1 in this new setting, and
	    stop it keeping the rightmost $h$ symbols in the right region
	    untouched. Then see (e).
	    \item The isomorphism $\iota$ is obtained as follows.
	    At the end of the procedure, the marked symbols signify the
	    contents of $b'_2$, and the unmarked symbols signify the
	    contents of $b'_1$.
	\end{enumerate}
	\item If $a$ is a $-$ symbol ($-_a$)
	in the $i$-th row, look at the $i$-th row in
	the left region.
	\begin{enumerate}
		\item If there are unmarked $+$ symbols
		in the
		$i$-th row in the left region,
		pick one of them and call it $+_c$.
	    \begin{enumerate}
	    	\item If there is no unmarked symbols
	    	(besides $+_c$) neither in the
	    	$i$-th row nor in the lower rows in the left region, then
	    	$+_c$ itself is identified with the partner $b(=+_b)$.
	    	\item If there are
	        unmarked symbols (besides $+_c$) either in the $i$-th
	        row or in
	        the lower rows in the left region, move $-_a$ and $+_c$ to
	        the $(i-1)$-th row.
	        Then apply the procedure (b).
	    \end{enumerate}
	    \item If there is no unmarked $+$ symbol in the
		$i$-th row in the left region, or one has already done the
		procedure (a)-ii,
	    \begin{enumerate}
	    	\item If there are unmarked $-$ symbols
	    	in the left
	    	region whose positions are lower than that of $-_a$,
	    	then the partner $b(=-_b)$ is chosen from one of those
	    	$-$ symbols that has the highest
	    	position.
	    	\item If there is no unmarked $-$ symbol in the left
	    	region whose positions are lower than that of $-_a$,
	    	then the partner $b(=+_b)$ is chosen from one of the unmarked $+$ symbols
	    	in the left region that has the lowest position.
	    \end{enumerate}
	\end{enumerate}
		\item If $a$ is a $+$ symbol ($+_a$),
		then the partner $b(=+_b)$ is chosen from one of the unmarked $+$ symbols
	    whose positions
	    are higher than that of $+_a$ but the lowest among them.
	    \item If the partner $b$ is a $-$ symbol ($-_b$),
	    mark it.
	    If $b$ is a $+$ symbol ($+_b$) in the $j$-th row,
	    look at the $j$-th row in the left and the middle region.
	    \begin{enumerate}
	    	\item If there are unmarked $-$ symbols in the
	    	$j$-th row
	    	either in the left region or in the middle region, pick
	    	the leftmost one of them and call it $-_d$.
	        Move $+_b$ and $-_d$ to the $(j+1)$-th row and
	        then mark the $+_b$.
	        \item If there is no unmarked $-$ symbol in the
	    	$j$-th row
	    	neither in the left region
	    	nor in the middle region, then mark the $+_b$.
	    \end{enumerate}
\end{enumerate}
This description of the rule is derived from the column insertion rule
in Section \ref{subsec:isorule} accompanied with the
{\itshape reverse row insertion procedure} for the $C$-tableaux.
We do not describe the latter procedure in this paper.

In the rule 4-(a), we have chosen $-_d$ to be the leftmost one.
However the final result of the procedure is actually the same 
for any choice of the $-$ symbols in the $j$-th row of the middle and left regions.

\subsection{Examples}
Let us present two examples.
We signify the marked symbols with  circles.
\subsubsection{Example 1}
Let us derive
\begin{equation}
\iota  : 1134\bar{3}\bar{2}\bar{1} \otimes
\bar{4} \bar{4} \bar{4} \bar{1} \bar{1}
\mapsto 144\bar{2}\bar{1} \otimes  0 \bar{4} \bar{4} \bar{4} \bar{4} \bar{1} \bar{0}
\end{equation}
under the isomorphism $B_{7} \otimes B_{5} \simeq B_{5} \otimes B_{7}$ of
the $U_q'(C_4^{(1)})$ crystals.
The value of the energy function is $0$ for this element.
The initial diagram is as follows.
\par\noindent
\setlength{\unitlength}{4mm}
\begin{picture}(15,5)(-5,0)
\put(0,0){\makebox(1,1){$\scriptstyle{4}$}}
\put(0,1){\makebox(1,1){$\scriptstyle{3}$}}
\put(0,2){\makebox(1,1){$\scriptstyle{2}$}}
\put(0,3){\makebox(1,1){$\scriptstyle{1}$}}
\put(0,4){\makebox(1,1){$\scriptstyle{0}$}}
\put(1,3){\makebox(1,1){$-$}}
\put(2,2){\makebox(1,1){$-$}}
\put(3,1){\makebox(1,1){$-$}}
\put(4,0){\makebox(1,1){$+$}}
\put(5,1){\makebox(1,1){$+$}}
\put(6,3){\makebox(1,1){$+$}}
\put(7,3){\makebox(1,1){$+$}}
\put(8,0){\line(0,1){5}}
\put(8,3){\makebox(1,1){$-$}}
\put(9,3){\makebox(1,1){$-$}}
\put(10,0){\makebox(1,1){$-$}}
\put(11,0){\makebox(1,1){$-$}}
\put(12,0){\makebox(1,1){$-$}}
\end{picture}
\par\noindent
We apply 2-(a)-ii, and then 2-(b)-i.
\par\noindent
\setlength{\unitlength}{4mm}
\begin{picture}(15,5)(-5,0)
\put(0,0){\makebox(1,1){$\scriptstyle{4}$}}
\put(0,1){\makebox(1,1){$\scriptstyle{3}$}}
\put(0,2){\makebox(1,1){$\scriptstyle{2}$}}
\put(0,3){\makebox(1,1){$\scriptstyle{1}$}}
\put(0,4){\makebox(1,1){$\scriptstyle{0}$}}
\put(1,3){\makebox(1,1){$-$}}
\put(1,3){\makebox(1,1){$\bigcirc$}}
\put(2,2){\makebox(1,1){$-$}}
\put(3,1){\makebox(1,1){$-$}}
\put(4,0){\makebox(1,1){$+$}}
\put(5,1){\makebox(1,1){$+$}}
\put(6,3){\makebox(1,1){$+$}}
\put(7,4){\makebox(1,1){$+$}}
\put(8,0){\line(0,1){5}}
\put(8,4){\makebox(1,1){$-$}}
\put(9,0){\line(0,1){5}}
\put(9,3){\makebox(1,1){$-$}}
\put(10,0){\makebox(1,1){$-$}}
\put(11,0){\makebox(1,1){$-$}}
\put(12,0){\makebox(1,1){$-$}}
\end{picture}
\par\noindent
Again, we apply 2-(a)-ii, and then 2-(b)-i.
\par\noindent
\setlength{\unitlength}{4mm}
\begin{picture}(15,5)(-5,0)
\put(0,0){\makebox(1,1){$\scriptstyle{4}$}}
\put(0,1){\makebox(1,1){$\scriptstyle{3}$}}
\put(0,2){\makebox(1,1){$\scriptstyle{2}$}}
\put(0,3){\makebox(1,1){$\scriptstyle{1}$}}
\put(0,4){\makebox(1,1){$\scriptstyle{0}$}}
\put(1,3){\makebox(1,1){$-$}}
\put(1,3){\makebox(1,1){$\bigcirc$}}
\put(2,2){\makebox(1,1){$-$}}
\put(2,2){\makebox(1,1){$\bigcirc$}}
\put(3,1){\makebox(1,1){$-$}}
\put(4,0){\makebox(1,1){$+$}}
\put(5,1){\makebox(1,1){$+$}}
\put(6,4){\makebox(1,1){$+$}}
\put(7,4){\makebox(1,1){$+$}}
\put(8,0){\line(0,1){5}}
\put(8,4){\makebox(1,1){$-$}}
\put(9,4){\makebox(1,1){$-$}}
\put(10,0){\line(0,1){5}}
\put(10,0){\makebox(1,1){$-$}}
\put(11,0){\makebox(1,1){$-$}}
\put(12,0){\makebox(1,1){$-$}}
\end{picture}
\par\noindent
We apply 2-(a)-i.
\par\noindent
\setlength{\unitlength}{4mm}
\begin{picture}(15,5)(-5,0)
\put(0,0){\makebox(1,1){$\scriptstyle{4}$}}
\put(0,1){\makebox(1,1){$\scriptstyle{3}$}}
\put(0,2){\makebox(1,1){$\scriptstyle{2}$}}
\put(0,3){\makebox(1,1){$\scriptstyle{1}$}}
\put(0,4){\makebox(1,1){$\scriptstyle{0}$}}
\put(1,3){\makebox(1,1){$-$}}
\put(1,3){\makebox(1,1){$\bigcirc$}}
\put(2,2){\makebox(1,1){$-$}}
\put(2,2){\makebox(1,1){$\bigcirc$}}
\put(3,1){\makebox(1,1){$-$}}
\put(4,0){\makebox(1,1){$+$}}
\put(4,0){\makebox(1,1){$\bigcirc$}}
\put(5,1){\makebox(1,1){$+$}}
\put(6,4){\makebox(1,1){$+$}}
\put(7,4){\makebox(1,1){$+$}}
\put(8,0){\line(0,1){5}}
\put(8,4){\makebox(1,1){$-$}}
\put(9,4){\makebox(1,1){$-$}}
\put(10,0){\makebox(1,1){$-$}}
\put(11,0){\line(0,1){5}}
\put(11,0){\makebox(1,1){$-$}}
\put(12,0){\makebox(1,1){$-$}}
\end{picture}
\par\noindent
We apply 2-(b)-ii to find the partner, and then apply 4-(a) to
mark the partner.
\par\noindent
\setlength{\unitlength}{4mm}
\begin{picture}(15,5)(-5,0)
\put(0,0){\makebox(1,1){$\scriptstyle{4}$}}
\put(0,1){\makebox(1,1){$\scriptstyle{3}$}}
\put(0,2){\makebox(1,1){$\scriptstyle{2}$}}
\put(0,3){\makebox(1,1){$\scriptstyle{1}$}}
\put(0,4){\makebox(1,1){$\scriptstyle{0}$}}
\put(1,3){\makebox(1,1){$-$}}
\put(1,3){\makebox(1,1){$\bigcirc$}}
\put(2,2){\makebox(1,1){$-$}}
\put(2,2){\makebox(1,1){$\bigcirc$}}
\put(3,0){\makebox(1,1){$-$}}
\put(4,0){\makebox(1,1){$+$}}
\put(4,0){\makebox(1,1){$\bigcirc$}}
\put(5,0){\makebox(1,1){$+$}}
\put(5,0){\makebox(1,1){$\bigcirc$}}
\put(6,4){\makebox(1,1){$+$}}
\put(7,4){\makebox(1,1){$+$}}
\put(8,0){\line(0,1){5}}
\put(8,4){\makebox(1,1){$-$}}
\put(9,4){\makebox(1,1){$-$}}
\put(10,0){\makebox(1,1){$-$}}
\put(11,0){\makebox(1,1){$-$}}
\put(12,0){\line(0,1){5}}
\put(12,0){\makebox(1,1){$-$}}
\end{picture}
\par\noindent
Again we apply 2-(b)-ii to find the partner, and then apply 4-(a) to
mark the partner.
\par\noindent
\setlength{\unitlength}{4mm}
\begin{picture}(15,5)(-5,0)
\put(0,0){\makebox(1,1){$\scriptstyle{4}$}}
\put(0,1){\makebox(1,1){$\scriptstyle{3}$}}
\put(0,2){\makebox(1,1){$\scriptstyle{2}$}}
\put(0,3){\makebox(1,1){$\scriptstyle{1}$}}
\put(0,4){\makebox(1,1){$\scriptstyle{0}$}}
\put(1,3){\makebox(1,1){$-$}}
\put(1,3){\makebox(1,1){$\bigcirc$}}
\put(2,2){\makebox(1,1){$-$}}
\put(2,2){\makebox(1,1){$\bigcirc$}}
\put(3,0){\makebox(1,1){$-$}}
\put(4,0){\makebox(1,1){$+$}}
\put(4,0){\makebox(1,1){$\bigcirc$}}
\put(5,0){\makebox(1,1){$+$}}
\put(5,0){\makebox(1,1){$\bigcirc$}}
\put(6,3){\makebox(1,1){$+$}}
\put(6,3){\makebox(1,1){$\bigcirc$}}
\put(7,4){\makebox(1,1){$+$}}
\put(8,0){\line(0,1){5}}
\put(8,4){\makebox(1,1){$-$}}
\put(9,3){\makebox(1,1){$-$}}
\put(10,0){\makebox(1,1){$-$}}
\put(11,0){\makebox(1,1){$-$}}
\put(12,0){\makebox(1,1){$-$}}
\put(13,0){\line(0,1){5}}
\end{picture}
\par\noindent
The procedure is finished.
Here the set of marked symbols stands for $144\bar{2}\bar{1} \in B_{5}$,
and the set of unmarked symbols stands for
$0\bar{4}\bar{4}\bar{4}\bar{4}\bar{1}\bar{0} \in B_{7}$.

Let us check the isomorphism by the definition.
\begin{equation}
	\begin{array}{ccc}
	1134\bar{3}\bar{2}\bar{1} \otimes
	\bar{4} \bar{4} \bar{4} \bar{1} \bar{1} &
	\mapright{\iota} &
	144\bar{2}\bar{1} \otimes
	0\bar{4} \bar{4} \bar{4} \bar{4} \bar{1} \bar{0} \\
	\mapdown{\psi_1} & & \mapdown{\psi_1} \\
	1111111 \otimes 122\bar{1}\bar{1} & \mapright{\iota} &
	11111 \otimes 01122\bar{1}\bar{0}
	\end{array}
\end{equation}
where
\begin{equation}
\psi_1 = (\et{2})^2 (\et{3})^2 (\et{1})^6 (\et{2})^6 (\et{3})^4
(\et{4})^6 (\et{3})^4 (\et{2})^2 \et{1}.
\end{equation}
We arrive at a $U_q(C_4)$ highest element.
\begin{equation}
	\begin{array}{ccc}
	1111111 \otimes 122\bar{1}\bar{1} & \mapright{\iota} &
	11111 \otimes 01122\bar{1}\bar{0} \\
	\mapdown{\psi_2} & & \mapdown{\psi_2} \\
	1111111 \otimes 11111 & \mapright{\iota} &
	11111 \otimes 1111111
	\end{array}
\end{equation}
where
\begin{equation}
\psi_2 = (\ftd_0)^8 (\ftd_1)^3 (\ftd_2)^3 (\ftd_3)^3
(\ftd_4)^3 (\ftd_3)^3 (\ftd_2)^3 \ftd_1 (\etd_0)^3.
\end{equation}
The energy was raised by $1$ when the third $\etd_0$ was applied, and
lowered by $1$ when the first $\ftd_0$ was applied.
Then it was raised by $5$ when the fourth to the
eighth $\ftd_0$'s were applied.
\subsubsection{Example 2}
Let us derive
\begin{equation}
\iota  : 0\bar{2}\bar{2}\bar{1}\bar{1}\bar{1}\bar{0} \otimes
1112\bar{2} \bar{1}
\mapsto
0\bar{2}\bar{1}\bar{1}\bar{1}\bar{0} \otimes
1112\bar{2}\bar{2} \bar{1}
\end{equation}
under the isomorphism  $B_{7} \otimes B_{6} \simeq B_{6} \otimes B_{7}$ of the
$U_q'(C_2^{(1)})$ crystals.
The value of the energy function is $4$ for this element.
The initial diagram is as follows.
\par\noindent
\setlength{\unitlength}{4mm}
\begin{picture}(15,3)(-5,0)
\put(0,0){\makebox(1,1){$\scriptstyle{2}$}}
\put(0,1){\makebox(1,1){$\scriptstyle{1}$}}
\put(0,2){\makebox(1,1){$\scriptstyle{0}$}}
\put(1,2){\makebox(1,1){$-$}}
\put(2,1){\makebox(1,1){$-$}}
\put(3,1){\makebox(1,1){$-$}}
\put(4,1){\makebox(1,1){$-$}}
\put(5,0){\makebox(1,1){$-$}}
\put(6,0){\makebox(1,1){$-$}}
\put(7,2){\makebox(1,1){$+$}}
\put(8,0){\line(0,1){3}}
\put(8,1){\makebox(1,1){$-$}}
\put(9,0){\makebox(1,1){$-$}}
\put(10,0){\makebox(1,1){$+$}}
\put(11,1){\makebox(1,1){$+$}}
\put(12,1){\makebox(1,1){$+$}}
\put(13,1){\makebox(1,1){$+$}}
\end{picture}
\par\noindent
We apply 2-(b)-i.
\par\noindent
\setlength{\unitlength}{4mm}
\begin{picture}(15,3)(-5,0)
\put(0,0){\makebox(1,1){$\scriptstyle{2}$}}
\put(0,1){\makebox(1,1){$\scriptstyle{1}$}}
\put(0,2){\makebox(1,1){$\scriptstyle{0}$}}
\put(1,2){\makebox(1,1){$-$}}
\put(2,1){\makebox(1,1){$-$}}
\put(3,1){\makebox(1,1){$-$}}
\put(4,1){\makebox(1,1){$-$}}
\put(5,0){\makebox(1,1){$-$}}
\put(6,0){\makebox(1,1){$-$}}
\put(6,0){\makebox(1,1){$\bigcirc$}}
\put(7,2){\makebox(1,1){$+$}}
\put(8,0){\line(0,1){3}}
\put(8,1){\makebox(1,1){$-$}}
\put(9,0){\line(0,1){3}}
\put(9,0){\makebox(1,1){$-$}}
\put(10,0){\makebox(1,1){$+$}}
\put(11,1){\makebox(1,1){$+$}}
\put(12,1){\makebox(1,1){$+$}}
\put(13,1){\makebox(1,1){$+$}}
\end{picture}
\par\noindent
We apply 2-(b)-ii to find the partner, and then apply 4-(a) to
mark the partner.
\par\noindent
\setlength{\unitlength}{4mm}
\begin{picture}(15,3)(-5,0)
\put(0,0){\makebox(1,1){$\scriptstyle{2}$}}
\put(0,1){\makebox(1,1){$\scriptstyle{1}$}}
\put(0,2){\makebox(1,1){$\scriptstyle{0}$}}
\put(1,1){\makebox(1,1){$-$}}
\put(2,1){\makebox(1,1){$-$}}
\put(3,1){\makebox(1,1){$-$}}
\put(4,1){\makebox(1,1){$-$}}
\put(5,0){\makebox(1,1){$-$}}
\put(6,0){\makebox(1,1){$-$}}
\put(6,0){\makebox(1,1){$\bigcirc$}}
\put(7,1){\makebox(1,1){$+$}}
\put(7,1){\makebox(1,1){$\bigcirc$}}
\put(8,0){\line(0,1){3}}
\put(8,1){\makebox(1,1){$-$}}
\put(9,0){\makebox(1,1){$-$}}
\put(10,0){\line(0,1){3}}
\put(10,0){\makebox(1,1){$+$}}
\put(11,1){\makebox(1,1){$+$}}
\put(12,1){\makebox(1,1){$+$}}
\put(13,1){\makebox(1,1){$+$}}
\end{picture}
\par\noindent
This time we find that there is no partner in the left region for
the leftmost $+$ symbol in the right region.
We interrupt the procedure here according to 1-(b).
Since there are four symbols in the right region,
we find that the value of the energy function is equal to $4$.
Following 1-(d) we give up this diagram and go to the initial diagram with
four marked symbols.
\par\noindent
\setlength{\unitlength}{4mm}
\begin{picture}(15,3)(-5,0)
\put(0,0){\makebox(1,1){$\scriptstyle{2}$}}
\put(0,1){\makebox(1,1){$\scriptstyle{1}$}}
\put(0,2){\makebox(1,1){$\scriptstyle{0}$}}
\put(1,2){\makebox(1,1){$-$}}
\put(1,2){\makebox(1,1){$\bigcirc$}}
\put(2,1){\makebox(1,1){$-$}}
\put(2,1){\makebox(1,1){$\bigcirc$}}
\put(3,1){\makebox(1,1){$-$}}
\put(3,1){\makebox(1,1){$\bigcirc$}}
\put(4,1){\makebox(1,1){$-$}}
\put(4,1){\makebox(1,1){$\bigcirc$}}
\put(5,0){\makebox(1,1){$-$}}
\put(6,0){\makebox(1,1){$-$}}
\put(7,2){\makebox(1,1){$+$}}
\put(8,0){\line(0,1){3}}
\put(8,1){\makebox(1,1){$-$}}
\put(9,0){\makebox(1,1){$-$}}
\put(10,0){\makebox(1,1){$+$}}
\put(11,1){\makebox(1,1){$+$}}
\put(12,1){\makebox(1,1){$+$}}
\put(13,1){\makebox(1,1){$+$}}
\end{picture}
\par\noindent
We apply 2-(b)-i.
\par\noindent
\setlength{\unitlength}{4mm}
\begin{picture}(15,3)(-5,0)
\put(0,0){\makebox(1,1){$\scriptstyle{2}$}}
\put(0,1){\makebox(1,1){$\scriptstyle{1}$}}
\put(0,2){\makebox(1,1){$\scriptstyle{0}$}}
\put(1,2){\makebox(1,1){$-$}}
\put(1,2){\makebox(1,1){$\bigcirc$}}
\put(2,1){\makebox(1,1){$-$}}
\put(2,1){\makebox(1,1){$\bigcirc$}}
\put(3,1){\makebox(1,1){$-$}}
\put(3,1){\makebox(1,1){$\bigcirc$}}
\put(4,1){\makebox(1,1){$-$}}
\put(4,1){\makebox(1,1){$\bigcirc$}}
\put(5,0){\makebox(1,1){$-$}}
\put(6,0){\makebox(1,1){$-$}}
\put(6,0){\makebox(1,1){$\bigcirc$}}
\put(7,2){\makebox(1,1){$+$}}
\put(8,0){\line(0,1){3}}
\put(8,1){\makebox(1,1){$-$}}
\put(9,0){\line(0,1){3}}
\put(9,0){\makebox(1,1){$-$}}
\put(10,0){\makebox(1,1){$+$}}
\put(11,1){\makebox(1,1){$+$}}
\put(12,1){\makebox(1,1){$+$}}
\put(13,1){\makebox(1,1){$+$}}
\end{picture}
\par\noindent
We apply 2-(b)-ii to find the partner, and then apply 4-(b) to
mark the partner.
\par\noindent
\setlength{\unitlength}{4mm}
\begin{picture}(15,3)(-5,0)
\put(0,0){\makebox(1,1){$\scriptstyle{2}$}}
\put(0,1){\makebox(1,1){$\scriptstyle{1}$}}
\put(0,2){\makebox(1,1){$\scriptstyle{0}$}}
\put(1,2){\makebox(1,1){$-$}}
\put(1,2){\makebox(1,1){$\bigcirc$}}
\put(2,1){\makebox(1,1){$-$}}
\put(2,1){\makebox(1,1){$\bigcirc$}}
\put(3,1){\makebox(1,1){$-$}}
\put(3,1){\makebox(1,1){$\bigcirc$}}
\put(4,1){\makebox(1,1){$-$}}
\put(4,1){\makebox(1,1){$\bigcirc$}}
\put(5,0){\makebox(1,1){$-$}}
\put(6,0){\makebox(1,1){$-$}}
\put(6,0){\makebox(1,1){$\bigcirc$}}
\put(7,2){\makebox(1,1){$+$}}
\put(7,2){\makebox(1,1){$\bigcirc$}}
\put(8,0){\line(0,1){3}}
\put(8,1){\makebox(1,1){$-$}}
\put(9,0){\makebox(1,1){$-$}}
\put(10,0){\line(0,1){3}}
\put(10,0){\makebox(1,1){$+$}}
\put(11,1){\makebox(1,1){$+$}}
\put(12,1){\makebox(1,1){$+$}}
\put(13,1){\makebox(1,1){$+$}}
\end{picture}
\par\noindent
The procedure is finished.
Here the set of marked symbols stands for
$0\bar{2}\bar{1}\bar{1}\bar{1}\bar{0} \in B_{6}$,
and the set of unmarked symbols stands for
$1112\bar{2}\bar{2} \bar{1} \in B_{7}$.
\section{\mathversion{bold}$C_n^{(1)}$ Kostka polynomials}
\label{sec:kostka}
Let $\mu_1 \geq \cdots \geq \mu_L (\geq 1)$ be a set of integers.
We set $\mu =(\mu_1, \mu_2, \ldots, \mu_L)$.
Consider the tensor product of $U_q'(C_n^{(1)})$ crystals
$B_{\mu_1}\ot\cd\ot B_{\mu_L}$.
Let
$\lambda = (\lambda_1, \ldots, \lambda_n)$ be another partition
satisfying $|\lambda| \leq |\mu| \,\mbox{and}\,|\lambda| \equiv |\mu| \pmod{2}$.
We define a classically restricted 1dsum
(cf. \cite{HKOTY}):
\begin{equation}\label{eq:Xdef}
X_{\lambda,\mu}(t)
=\mathop{{\sum}^*}t^{\sum_{0\le i<j\le L}H(b_i\ot b^{(i+1)}_j)},
\end{equation}
where the sum $\sum^*$ is taken over all
$b_1\ot\cd\ot b_L\in B_{\mu_1}\ot\cd\ot B_{\mu_L}$
satisfying
\begin{equation*}
\tilde{e}_i(b_1 \ot \cdots \ot b_L) = 0,\;
\varphi_i(b_1 \ot \cdots \ot b_L) = \lambda_i - \lambda_{i+1},
 \quad 1 \le i \le n  \;(\lambda_{n+1}=0).
\end{equation*}
This condition is equivalent to
\begin{equation}
\label{eq:hwcond}
(b_L \stackrel{\ast}{\longrightarrow}\cdots \stackrel{\ast}{\longrightarrow} ( b_3
\stackrel{\ast}{\longrightarrow} (b_2
\stackrel{\ast}{\longrightarrow} b_1 ) ) \cdots ) =
T(\lambda).
\end{equation}
$T(\lambda)$ is the unique tableau of both
shape and weight $\lambda$.
Namely all the letters in the first row are $1$ and
those in the second row are $2$, and so on.
(We distinguished $b_j$ from $\T (b_j)$.)
In the summation we set $b^{(i)}_i=b_i$,
and $b^{(i)}_j$ ($i<j$) are defined by
successive use of the crystal isomorphism,
\begin{eqnarray*}
&&
\begin{array}{ccccc}\hspace{-5mm}
B_{\mu_i}\ot\cd\ot B_{\mu_{j-1}}\ot B_{\mu_j}&\simeq&
B_{\mu_i}\ot\cd\ot B_{\mu_j}\ot B_{\mu_{j-1}}&\simeq&\cd\\
b_i\ot\cd\ot b_{j-1}\ot b_j&\mapsto&
b_i\ot\cd\ot b^{(j-1)}_j\ot b'_{j-1}&\mapsto&\cd
\end{array}\\
&&\hspace{5cm}
\begin{array}{ccc}
\cd&\simeq&B_{\mu_j}\ot B_{\mu_i}\ot\cd\ot B_{\mu_{j-1}}\\
\cd&\mapsto&b^{(i)}_j\ot b'_i\ot\cd\ot b'_{j-1}.
\end{array}
\end{eqnarray*}
The above condition (\ref{eq:hwcond}) implies that
$b^{(1)}_j$
in the tableau presentation should have the form
$\T (b^{(1)}_j)=\fbx{0 \cd 0}\fbx{1 \cd 1}\fbx{\ol{0} \cd \ol{0}}$.
$b_0$ is chosen so that $H(b_0 \ot b^{(1)}_j) = - \sharp (\fbx{0} \,
\mbox{in} \, \T (b^{(1)}_j))$.
Up to additive constant this agrees with
the choice of $b_0$ in \cite{HKOTY}.
Up to an overall power of $t$, this is a polynomial which may be viewed as
a $C^{(1)}_n$-analogue of the Kostka polynomial.
In fact if  $|\lambda| = |\mu|$, $X_{\lambda,\mu}(t)$ coincides with
the ordinary Kostka polynomial $K_{\lambda,\mu}(t)$.

Following the tables in pp.~239-240 of \cite{Ma}
we give a list of $X_{\lambda,\mu}(t)$ or the matrices
$X(t) := \{ X_{\lambda,\mu}(t) \}$ for $|\mu| \leq 6$ and
$|\lambda| = |\mu|-2,|\mu|-4,\ldots$ with $n \ge L$.
$X_{\lambda, \mu}(t)$ is independent of $n$ if $n \ge L$.
If $n< L$ it is  $n$-dependent in general.
For instance, let $\lambda = (1^3)$ and $\mu = (1^5)$.
The element $\fbx{1} \ot \fbx{2} \ot \fbx{3} \ot \fbx{4} \ot \fbx{\ol{4}}
\in (B_{1})^{\ot 5}$ contributes to $X_{(1^3),(1^5)}(t)$ for $n \geq 4$,
but does not for $n=2,3$.
We have checked that all the data in the table agrees
with the fermionic formula in \cite{HKOTY}.
In the tables, a row (resp. column) specifies $\lambda$ (resp. $\mu$) 
in $X_{\lambda, \mu}(t)$.

\begin{displaymath}
X_{\emptyset,(2)}(t) = t^{-1}, \quad X_{\emptyset,(1^2)}(t) = 1.
\end{displaymath}

\begin{displaymath}
X_{(1),(3)}(t) = t^{-1}, \quad X_{(1),(21)}(t) = t^{-1}+1, \quad
X_{(1),(1^3)}(t) = 1+t+t^2.
\end{displaymath}

\par\noindent
$
\footnotesize
\begin{array}{|c||c|c|c|c|c|}\hline
	  & ( 4 ) & ( 3 1 ) & ( 2^{2} ) & ( 2 1^{2} ) & ( 1^{4} )  \\\hline\hline
	 \emptyset &t^{\!-\!2} &t^{\!-\!1} &t^{\!-\!2}\!+\!1 &t^{\!-\!1}\!+\!t &1\!+\!{t^2}\!+\!{t^4}  \\\hline
	 ( 2 ) &t^{\!-\!1} &t^{\!-\!1}\!+\!1 &t^{\!-\!1}\!+\!1\!+\!t &2\!+\!t\!+\!{t^2} &t\!+\!{t^2}\!+\!2{t^3}\!+\!{t^4}\!+\!{t^5}  \\\hline
	 ( 1^{2} ) & &t^{\!-\!1} &1&t^{\!-\!1}\!+\!1\!+\!t &1\!+\!t\!+\!2{t^2}\!+\!{t^3}\!+\!{t^4}  \\\hline
\end{array}
$
\vskip3ex

\par\noindent $\footnotesize
\begin{array}{|c||c|c|c|c|c|}\hline
	  & ( 5 ) & ( 4 1 ) & ( 3 2 ) & ( 3 1^{2} ) & ( 2^{2} 1 )  \\\hline\hline
	 ( 1 ) &t^{\!-\!2} &t^{\!-\!2}\!+\!t^{\!-\!1} &t^{\!-\!2}\!+\!t^{\!-\!1}\!+\!1 &2t^{\!-\!1}\!+\!1\!+\!t &t^{\!-\!2}\!+\!t^{\!-\!1}\!+\!2\!+\!t\!+\!{t^2}  \\\hline
	 ( 3 ) &t^{\!-\!1} &t^{\!-\!1}\!+\!1 &t^{\!-\!1}\!+\!1\!+\!t &2\!+\!t\!+\!{t^2} &1\!+\!2t\!+\!{t^2}\!+\!{t^3}  \\\hline
	 ( 2 1 ) & &t^{\!-\!1} &t^{\!-\!1}\!+\!1 &t^{\!-\!1}\!+\!2\!+\!t &t^{\!-\!1}\!+\!2\!+\!2t\!+\!{t^2}  \\\hline
	 ( 1^{3} ) & & & &t^{\!-\!1} &1 \\\hline
\end{array}
$
\vskip3ex

\par\noindent $\footnotesize
\begin{array}{|c||c|c|}\hline
	  & ( 2 1^{3} ) & ( 1^{5} )  \\\hline\hline
	 ( 1 ) &t^{\!-\!1}\!+\!2\!+\!2t\!+\!2{t^2}\!+\!{t^3}\!+\!{t^4} &1\!+\!t\!+\!2{t^2}\!+\!2{t^3}\!+\!3{t^4}\!+\!2{t^5}\!+\!2{t^6}\!+\!{t^7}\!+\!{t^8}  \\\hline
	 ( 3 ) &t\!+\!2{t^2}\!+\!2{t^3}\!+\!{t^4}\!+\!{t^5} &{t^3}\!+\!{t^4}\!+\!2{t^5}\!+\!2{t^6}\!+\!2{t^7}\!+\!{t^8}\!+\!{t^9}  \\\hline
	 ( 2 1 ) &2\!+\!3t\!+\!3{t^2}\!+\!2{t^3}\!+\!{t^4} &t\!+\!2{t^2}\!+\!3{t^3}\!+\!4{t^4}\!+\!4{t^5}\!+\!3{t^6}\!+\!2{t^7}\!+\!{t^8}  \\\hline
	 ( 1^{3} ) &t^{\!-\!1}\!+\!1\!+\!t\!+\!{t^2} &1\!+\!t\!+\!2{t^2}\!+\!2{t^3}\!+\!2{t^4}\!+\!{t^5}\!+\!{t^6}  \\\hline
\end{array}
$
\vskip3ex

\par\noindent $\footnotesize
\begin{array}{|c||c|c|c|c|c|c|}\hline
	  & ( 6 ) & ( 5 1 ) & ( 4 2 ) & ( 4 1^{2} ) & ( 3^{2} ) & ( 3 2 1 )  \\\hline\hline
	 \emptyset &t^{\!-\!3} &t^{\!-\!2} &t^{\!-\!3}\!+\!t^{\!-\!1} &t^{\!-\!2}\!+\!1 &t^{\!-\!2}\!+\!1 &t^{\!-\!2}\!+\!t^{\!-\!1}\!+\!t  \\\hline
	 ( 2 ) &t^{\!-\!2} &t^{\!-\!2}\!+\!t^{\!-\!1} &2t^{\!-\!2}\!+\!t^{\!-\!1}\!+\!1 &3t^{\!-\!1}\!+\!1\!+\!t &2t^{\!-\!1}\!+\!1\!+\!t &t^{\!-\!2}\!+\!2t^{\!-\!1}\!+\!3\!+\!t\!+\!{t^2}  \\\hline
	 ( 1^{2} ) & &t^{\!-\!2} &t^{\!-\!1} &t^{\!-\!2}\!+\!t^{\!-\!1}\!+\!1 &t^{\!-\!2}\!+\!1 &t^{\!-\!2}\!+\!2t^{\!-\!1}\!+\!1\!+\!t  \\\hline
	 ( 4 ) &t^{\!-\!1} &t^{\!-\!1}\!+\!1 &t^{\!-\!1}\!+\!1\!+\!t &2\!+\!t\!+\!{t^2} &1\!+\!t\!+\!{t^2} &1\!+\!2t\!+\!{t^2}\!+\!{t^3}  \\\hline
	 ( 3 1 ) & &t^{\!-\!1} &t^{\!-\!1}\!+\!1 &t^{\!-\!1}\!+\!2\!+\!t &t^{\!-\!1}\!+\!1\!+\!t &t^{\!-\!1}\!+\!3\!+\!2t\!+\!{t^2}  \\\hline
	 ( 2 1^{2} ) & & & &t^{\!-\!1} & &t^{\!-\!1}\!+\!1  \\\hline
	 ( 2^{2} ) & & &t^{\!-\!1} &1&1&t^{\!-\!1}\!+\!1\!+\!t  \\\hline
	 ( 1^{4} ) & & & & & &  \\\hline
\end{array}
$
\vskip3ex

\par\noindent $\footnotesize
\begin{array}{|c||c|c|}\hline
	  & ( 3 1^{3} ) & ( 2^{3} ) \\\hline\hline
	 \emptyset &t^{\!-\!1}\!+\!1\!+\!t\!+\!{t^3} &t^{\!-\!3}\!+\!t^{\!-\!1}\!+\!1\!+\!t\!+\!{t^3}   \\\hline
	 ( 2 ) &t^{\!-\!1}\!+\!3\!+\!3t\!+\!3{t^2}\!+\!{t^3}\!+\!{t^4} &t^{\!-\!2}\!+\!t^{\!-\!1}\!+\!3\!+\!2t\!+\!3{t^2}\!+\!{t^3}\!+\!{t^4}  \\\hline
	 ( 1^{2} ) &2t^{\!-\!1}\!+\!2\!+\!3t\!+\!{t^2}\!+\!{t^3} &t^{\!-\!1}\!+\!1\!+\!2t\!+\!{t^2}\!+\!{t^3}  \\\hline
	 ( 4 ) &t\!+\!2{t^2}\!+\!2{t^3}\!+\!{t^4}\!+\!{t^5} &t\!+\!{t^2}\!+\!2{t^3}\!+\!{t^4}\!+\!{t^5}  \\\hline
	 ( 3 1 ) &2\!+\!3t\!+\!4{t^2}\!+\!2{t^3}\!+\!{t^4} &1\!+\!2t\!+\!3{t^2}\!+\!2{t^3}\!+\!{t^4}  \\\hline
	 ( 2 1^{2} ) &t^{\!-\!1}\!+\!2\!+\!2t\!+\!{t^2} &1\!+\!t\!+\!{t^2}  \\\hline
	 ( 2^{2} ) &1\!+\!2t\!+\!{t^2}\!+\!{t^3} &t^{\!-\!1}\!+\!1\!+\!2t\!+\!{t^2}\!+\!{t^3}  \\\hline
	 ( 1^{4} ) &t^{\!-\!1} &  \\\hline
\end{array}
$
\vskip3ex

\par\noindent $\footnotesize
\begin{array}{|c||c|c|}\hline
            & ( 2^{2} 1^{2} )                                & ( 2 1^{4} )                                                      \\\hline\hline
\emptyset       & t^{\!-\!2}\!+\!2\!+\!t\!+\!{t^2}\!+\!{t^4}                         & t^{\!-\!1}\!+\!2t\!+\!{t^2}\!+\!2{t^3}\!+\!{t^4}\!+\!{t^5}\!+\!{t^7}                     \\\hline
( 2 )       & 2t^{\!-\!1}\!+\!2\!+\!5t\!+\!3{t^2}\!+\!3{t^3}\!+\!{t^4}\!+\!{t^5} & 2\!+\!2t\!+\!5{t^2}\!+\!4{t^3}\!+\!6{t^4}\!+\!3{t^5}\!+\!3{t^6}\!+\!{t^7}\!+\!{t^8}  \\\hline
( 1^{2} )   & t^{\!-\!2}\!+\!t^{\!-\!1}\!+\!4\!+\!2t\!+\!3{t^2}\!+\!{t^3}\!+\!{t^4}      & t^{\!-\!1}\!+\!2\!+\!4t\!+\!4{t^2}\!+\!5{t^3}\!+\!3{t^4}\!+\!3{t^5}\!+\!{t^6}\!+\!{t^7}    \\\hline
( 4 )       & 2{t^2}\!+\!2{t^3}\!+\!2{t^4}\!+\!{t^5}\!+\!{t^6}         & {t^3}\!+\!{t^4}\!+\!3{t^5}\!+\!2{t^6}\!+\!2{t^7}\!+\!{t^8}\!+\!{t^9}               \\\hline
( 3 1 )     & 1\!+\!4t\!+\!4{t^2}\!+\!4{t^3}\!+\!2{t^4}\!+\!{t^5}        & t\!+\!3{t^2}\!+\!5{t^3}\!+\!6{t^4}\!+\!5{t^5}\!+\!4{t^6}\!+\!2{t^7}\!+\!{t^8}    \\\hline
( 2 1^{2} ) & t^{\!-\!1}\!+\!2\!+\!3t\!+\!2{t^2}\!+\!{t^3}                   & 2\!+\!3t\!+\!5{t^2}\!+\!4{t^3}\!+\!4{t^4}\!+\!2{t^5}\!+\!{t^6}                 \\\hline
( 2^{2} )   & 2\!+\!2t\!+\!3{t^2}\!+\!{t^3}\!+\!{t^4}                    & 2t\!+\!2{t^2}\!+\!4{t^3}\!+\!3{t^4}\!+\!3{t^5}\!+\!{t^6}\!+\!{t^7}             \\\hline
( 1^{4} )   & 1                                              & t^{\!-\!1}\!+\!1\!+\!t\!+\!{t^2}\!+\!{t^3}                                           \\\hline
\end{array}
$
\vskip3ex

\par\noindent $\footnotesize
\begin{array}{|c||c|}\hline
	        & ( 1^{6} )  \\\hline\hline
\emptyset       & 1\!+\!{t^2}\!+\!{t^3}\!+\!2{t^4}\!+\!{t^5}\!+\!3{t^6}\!+\!{t^7}\!+\!2{t^8}\!+\!{t^9}\!+\!{t^{10}}\!+\!{t^{12}}  \\\hline
( 2 )       & t\!+\!{t^2}\!+\!3{t^3}\!+\!3{t^4}\!+\!6{t^5}\!+\!5{t^6}\!+\!7{t^7}\!+\!5{t^8}\!+\!6{t^9}\!+\!3{t^{10}}\!+\!3{t^{11}}\!+\!{t^{12}}\!+\!{t^{13}}  \\\hline
( 1^{2} )   & 1\!+\!t\!+\!3{t^2}\!+\!3{t^3}\!+\!6{t^4}\!+\!5{t^5}\!+\!7{t^6}\!+\!5{t^7}\!+\!6{t^8}\!+\!3{t^9}\!+\!3{t^{10}}\!+\!{t^{11}}\!+\!{t^{12}}  \\\hline
( 4 )       &{t^6}\!+\!{t^7}\!+\!2{t^8}\!+\!2{t^9}\!+\!3{t^{10}}\!+\!2{t^{11}}\!+\!2{t^{12}}\!+\!{t^{13}}\!+\!{t^{14}}  \\\hline
( 3 1 )     & {t^3}\!+\!2{t^4}\!+\!4{t^5}\!+\!5{t^6}\!+\!7{t^7}\!+\!7{t^8}\!+\!7{t^9}\!+\!5{t^{10}}\!+\!4{t^{11}}\!+\!2{t^{12}}\!+\!{t^{13}}  \\\hline
( 2 1^{2} ) & t\!+\!2{t^2}\!+\!4{t^3}\!+\!5{t^4}\!+\!7{t^5}\!+\!7{t^6}\!+\!7{t^7}\!+\!5{t^8}\!+\!4{t^9}\!+\!2{t^{10}}\!+\!{t^{11}}  \\\hline
( 2^{2} )   & {t^2}\!+\!{t^3}\!+\!3{t^4}\!+\!3{t^5}\!+\!5{t^6}\!+\!4{t^7}\!+\!5{t^8}\!+\!3{t^9}\!+\!3{t^{10}}\!+\!{t^{11}}\!+\!{t^{12}}  \\\hline
( 1^{4} )   & 1\!+\!t\!+\!2{t^2}\!+\!2{t^3}\!+\!3{t^4}\!+\!2{t^5}\!+\!2{t^6}\!+\!{t^7}\!+\!{t^8}  \\\hline
\end{array}$
\vskip3ex


\end{document}